\documentclass[aps,pra,superscriptaddress]{revtex4}

\pdfoutput=1
\usepackage[intlimits]{amsmath}
\usepackage{amsfonts}
\usepackage{verbatim}
\usepackage{enumitem}
\usepackage{psfrag}

\usepackage{amsthm}
\usepackage{float}
\usepackage{tikz}

\advance\voffset by -1cm
\advance\hoffset by -0.5cm
\newcommand\beq            {\begin{equation}}
\newcommand\be            {\begin{equation}}
\newcommand\bea           {\begin{equation}\begin{array}l\displaystyle}
\newcommand\ee            {\end{equation}}
\newcommand\eeq            {\end{equation}}
\newcommand\bes           {\begin{subequations}}
\newcommand\esu           {\end{subequations}}

\renewcommand{\(}{\left(}
\renewcommand{\)}{\right)}
\renewcommand{\[}{\left[}
\renewcommand{\]}{\right]}

\newcommand{\bigx}[1]{\bBigg@{#1}}

\newcommand\red[1]{{\color{red}{#1}} }

\def\3pt#1#2#3{{\langle{#1}\vert{#2}\vert{#3}\rangle}}

\newcommand\doi[2]        {\href{http://dx.doi.org/#1}{#2}}

\newcommand{\EQ}{\begin{equation}}
\newcommand{\EN}{\end{equation}}
\usepackage{epsfig}
\usepackage{psfrag}
\usepackage{amsmath}
\usepackage{graphicx}
\usepackage{amsfonts}
\usepackage{amssymb}

\def\tilde{\widetilde}

\def\hat{\widehat}
\def\*{\star}
\def\[{\left[}
\def\]{\right]}
\def\({\left(}      

\def\){\right)}

\def\frac#1#2{\dfrac{#1}{#2}}
\def\inv#1{\dfrac{1}{#1}}
\def\half{\tfrac{1}{2}}

\def\2pi{\hbox{$2\pi i$}}

\def\dsl{\raise.15ex\hbox{/}\kern-.57em\partial}
\def\Dsl{\,\raise.15ex\hbox{/}\mkern-.13.5mu D}
\def\be{\beta}

   \def\CE{{\cal E}}

   \def\CN{{\cal N}}   
\def\CP{{\cal P}}

\def\2pi{\hbox{$2\pi i$}}

\def\dsl{\raise.15ex\hbox{/}\kern-.57em\partial}
\def\Dsl{\,\raise.15ex\hbox{/}\mkern-.13.5mu D}

\def\barray{\begin{eqnarray}}
\def\earray{\end{eqnarray}}
\def\beq{\begin{equation}}
\def\eeq{\end{equation}}

\def\AA{\leavevmode\setbox0=\hbox{h}
\dimen0=\ht0 \advance\dimen0 by-1ex\rlap{\raise.67\dimen0\hbox{\char'27}}A}

\def\blue#1{{\color{blue}{#1}}}
\def\red#1{{\color{red}{#1}}}

\def\blue#1{{\color{blue}{#1}}}
\def\red#1{{\color{red}{#1}}}
\def\half{\tfrac{1}{2}}

\def\prob{{ \bf Pr}}
\def\Prob2{\CP}

\def\Ex{{\bf E}}

\def\Primes{\mathbb{P}}
\def\Primesp{\mathbb{P}'}
\def\Pensemble{{\bf \CE}}

\def\half{\tfrac{1}{2}}

\def\dist{{\overset{d}{\longrightarrow}}}

\def\CBN{B_N}

\newtheorem{theorem}{Theorem}

\newtheorem{conjecture}{Conjecture}

\begin{document}

\title{{\Large {\bf Generalized Riemann Hypothesis, Time Series \\
and Normal Distributions}}}

\author{ Andr\'e  LeClair}
\affiliation{Cornell University, Physics Department, Ithaca, NY 14850} 
\author{Giuseppe Mussardo}
\affiliation{SISSA and INFN, Sezione di Trieste, via Bonomea 265, I-34136, 
Trieste, Italy}

\begin{abstract}
$L$ functions based on Dirichlet characters are natural  generalizations of  the Riemann $\zeta(s)$ function:  they both have series representations and satisfy an Euler product representation, i.e. an infinite product taken over prime numbers. In this paper we address the Generalized Riemann Hypothesis relative to the non-trivial complex zeros of the 
Dirichlet $L$ functions by studying the possibility to enlarge the original domain of convergence of their Euler product. The feasibility of this analytic continuation is ruled by the asymptotic behavior in $N$ of the series  $B_N =  \sum_{n=1}^N \cos \( t \log p_n - \arg \chi (p_n) \)$ involving Dirichlet characters $\chi$ modulo $q$ on primes $p_n$. Although deterministic, these series have pronounced stochastic features which make them analogous to random time series.  We show that the $B_N$'s satisfy various normal law probability distributions. The study of their large asymptotic behavior poses an interesting problem of statistical physics equivalent to the Single Brownian Trajectory Problem, here addressed by defining an appropriate ensemble $\CE$ involving intervals of primes. For non-principal characters, we show that the series $B_N$ present a universal diffusive random walk behavior  $B_N = O(\sqrt{N})$ in view of the Dirichlet theorem on the equidistribution of reduced residue classes modulo $q$ and the Lemke Oliver-Soundararajan conjecture on the distribution of pairs of residues on consecutive primes. This purely diffusive behavior of $B_N$ implies that the domain of convergence of the infinite product representation of the Dirichlet $L$-functions for non-principal characters can be extended from $\Re(s) > 1$ down to $\Re (s)  = \half$, without encountering any zeros before reaching this critical line.    

\vspace{5mm}
{\em Dedicated to Giorgio Parisi on the occasion of his 70th birthday.}

\end{abstract}

\maketitle

\section{Introduction}

The real world confronts the mathematician with events that are not strictly predictable, and the methods developed to deal with them have opened new domains of pure mathematics. Nowadays the concept of probability plays a vital role in many fields of physics, chemistry or mathematics and appears as well in a wide range of many other phenomena, including computer science, finance or biology (see, for instance \cite{Jaynes,Kampen,Parisi}). A key object is the normal probability distribution together with the associated random walk: 
the ubiquity and robustness of the normal distribution comes of course from a key concept in probability theory, namely the central limit theorem, a result which holds under quite general conditions. There is indeed a large degree of universality and simplicity behind this law. Consider, for instance, the random walk: the rules at the back of this process are quite simple but their consequences can be far from elementary. This is particularly true in a subject as Number Theory, a field usually seen as highly stiff and deterministic in view of the rigidity of the discrete laws of arithmetics. However, in the course of the last decades, much progress has  been achieved in this field by exploring the interplay between randomness and determinism, two aspects which coexist in particular in the realm of prime numbers  
(see, for instance, 
\cite{Dyson,Montgomery,Odlyzko,Rudnick-Sarnack,Conrey,EPFchi,Kac,Cramer,Erdos-Kac,Billingsley,GrosswaldSchnitzer,Chernoff,Schroeder,Tao,SarnakMoebius,Wolf,Torquato,shortpaper}). Mark Kac, for instance, unveiled a wide spectrum of aspects of Number Theory ruled by normal distributions (see \cite{Kac}): a famous example is the Erd\"os and Kac result concerning the number of prime factors of the integers, which indeed obeys a normal distribution \cite{Erdos-Kac}. It is worth stressing that probabilistic arguments may  also be a source of inspiration in identifying hidden properties of Number Theory, as illustrated by the famous random model of primes proposed by Cram\'er in which he was able to prove concise statements with probability equal to one: within his random model of primes, an example of those statements is given by this inequality about the gaps of these \textquotedblleft random prime numbers\textquotedblright \cite{Cramer}   
\beq
\label{gaps} 
p_{n+1} - p_n < \log^2 p_n\,\,\,, 
\eeq
and, as a matter of fact,  no violation of this inequality has found so far in the actual set of prime numbers with $p_n > 7$. 
 
In this paper we are going to show that a random walk approach provides a key to establish the validity of the so-called Generalized Riemann Hypothesis (GRH) for the Dirichlet $L$-functions of non-principal characters. While the original arguments were presented in a previous publication by us \cite{shortpaper}, this paper not only provides their thorough and detailed discussion but also embeds such a discussion in a broader analysis involving several other probabilistic aspects relative to the Dirichlet $L$-functions.
In this introduction we shall give a brief account of the problem and the central idea of our approach, skipping many technical details which however will be addressed later in the paper. 

The main concern of this paper is the location of the non-trivial zeros of the Dirichlet $L$-functions $L(s,\chi)$ of the complex variable $s = \sigma + i t$ based on a Dirichlet character $\chi$. A detailed discussions of these quantities and  the relative $L$ functions can be found, for instance, in \cite{Apostol,Iwaniec,Bombieri,Steuding,Sarnak,Davenport}. For  $\Re(s) > 1$ these functions admit two equivalent representations, one given in terms of an infinite series on the natural numbers $m$, the other in terms of an infinite product over the sequence of primes $p_n$ (hereafter labelled in ascending order) 
\beq
L(s,\chi) \,=\,
\sum_{m=1}^\infty \frac{\chi(m)}{m^s}\,=\, 
 \prod_{n=1}^\infty  \( 1 -  \frac{\chi (p_n)}{p_n^s} \)^{-1} \,\,\,. 
\label{Euleridentity}
\eeq
The infinite product representation is known as the {\em Euler product formula}.  As shown below, when the characters are non-zero, they are pure phases and therefore expressed in terms of some angles $\theta_m$ defined as 
\beq
\label{thetan}
\chi (m) \,=\,  e^{i \theta_m }, ~~~~~\forall ~\chi(m) \neq 0 \,\,\,.
\eeq   
The characters $\chi(m)$ are completely multiplicative arithmetic functions, a property which is at the origin of the Euler product formula together with the unique decomposition of any integer in terms of primes. Notice that the $L$-functions provide a generalization of the Riemann $\zeta$-function \cite{Riemann1,Riemann2,Riemann3}, which is obtained taking all
 $\chi(m) = 1$. 

%\vspace{3mm}
Following the original Riemann Hypothesis for the Riemann $\zeta$-function \cite{Riemann1,Riemann2,Riemann3} but, at the same time, widely enlarging the perspective and the foundation of such a conjecture,  the Generalized Riemann Hypothesis states that the non-trivial zeros of {\em all} the infinitely many $L$- functions lie along the critical line $\Re(s) = \half$. According to Davenport \cite{Davenport}, this conjecture seems to have been first formulated by the German mathematician Adolf Piltz in 1884 and since then, a large number of papers have been dealing with this hypothesis, too large to do justice to all the many authors who contributed to the development of the subject. Here it may be enough to mention a few basic results about $L$-functions particularly important for our purposes. Selberg \cite{Selberg1} was the first to obtain the counting formula $N(T,\chi)$ for the number of zeros up to height $T$ within the entire critical strip $0 \le \Re(s) \le 1$. Fujii \cite{Fujii} later refine this result providing a formula for the number of zeros in the critical strip with the ordinate between $[T+H, T]$. For the low lying  zeros near and at the critical line, their distribution was analyzed by Iwaniec, Luo and Sarnak \cite{IwaniecS}, assuming however the validity of the GRH. As for the original Montgomery-Odlyzko conjecture relative to the zeros of the Riemann $\zeta$-function and their relation to random matrix theory \cite{Montgomery,Odlyzko} (see also \cite{Rudnick-Sarnack}), the statistics of the zeros of the $L$-functions were studied by Conrey and, in a separate paper, by Hughes and Rudnick \cite{Conrey,Hughes}. Interestingly enough, Conrey, Iwaniec, and Soundararajan have estimated that more than $56\%$  of the non-trivial zeros are on the critical line \cite{Conrey2}. These mathematical results are also accompanied by some interesting interpretations of these functions that come from two different fields of Physics.

\vspace{3mm}
\noindent
{\bf Statistical Physics Interpretation}. From a Statistical Physics point of view, the $L$-functions can be naturally interpreted as generalized quantum partition functions of free systems of particles: these particles have energies given by $E_n = \log p_n$ \cite{Julia,Spector} and a certain assignment of their electric charges (see Appendix \ref{AA}). From this point of view, the identity (\ref{Euleridentity}) can be seen as the formula which expresses the equivalence between the canonical and grand-canonical statistical ensembles of these free systems, and therefore the zeros of $L(s,\chi)$ are nothing else but the Fisher zeros of these statistical systems \cite{Fisher}. 

\vspace{3mm}
\noindent
{\bf Quantum Physics Interpretation}. Even more interesting is the profound interplay between the spectral theory of quantum mechanics and the zeros of the Dirichlet $L$-function, in particular those of the $\zeta$-Riemann function: originally stated by P\'olya, this viewpoint has given rise to an important series of works on the $\zeta$-Riemann function 
by Berry, Keating, Connes, Sierra, Srednicki, Bender and many others \cite{BK1,BK2,BK3,BK4,Connes,Sierra1,Sierra2,Srednicki,Bender} (for a more complete list of references, see the review \cite{reviewRiemann}). In a nutshell, these approaches to the Riemann Hypothesis for the $\zeta$-Riemann function have the aim to identify a quantum mechanical Hamiltonian whose spectrum coincides with the imaginary parts of the Riemann zeros along the line $\Re(s) = \half$, so that to argue that the alignment of all these zeros along this axis can be seen as a consequence of the spectral theory of quantum Hamiltonians. However, in spite of all these interesting works, it is probably fair to say that the sought after hamiltonian has remained so far elusive. 
 
\vspace{3mm}
\noindent
{\bf Random Walk Approach}.  As shown originally in the paper \cite{shortpaper}, the GRH can be approached in a completely different way. The starting point of this new approach 
comes from a simple remark: if all the infinitely many $L$-Dirichlet functions have their non-trivial zeros along the axis $\sigma = \half$, behind this fact there should be some universal and robust reason which transcends the details of the characters entering their definition, relying instead on some of the general properties of these quantities.  For the $L$-functions associated to the non principal characters, such a reason can be nailed down to the existence of a random walk, in the sense that the value $\sigma = \half$ can be identified with the critical value of a random walk process which exists for all these functions. 

Where does  this random walk process come from? As discussed in Section \ref{Probabilistic}, a way to establish the validity of the GRH for the $L$-Dirichlet functions consists of showing that their infinite product representation can be extended from $ \Re(s) >1 $ to the half-plane $\Re(s) > \half+ \epsilon$ for any $\epsilon >0$ arbitrarily small. In turn, this analytical continuation of the infinite product representation inside the critical strip $0< \Re(s) <1$ is controlled by the following series on the primes \cite{EPFchi} (see Theorem \ref{BNtheorem} below)
\beq
\label{Bprimes}
B_N( t, \chi) = \sum_{n=1}^N \cos (t \log p_n  - \theta_{p_n}) \,\,\,.
\eeq
For every character $\chi$ there is an associated $B_N$ series\footnote{We will not always display this $\chi$ dependence and simply write $B_N(t)$ for these series.}. For non-principal characters, the phases $\theta_{p_n}$ are different from zero while for principal characters all of them vanish. The fact that {\em all} the non-trivial zeros of the $L$-functions lie in the complex plane along the critical line $\Re(s) = \half$ will be guaranteed if, for all values of the variable $t$, at large $N$ the series $B_N(t)$ behaves as 
\beq
B_N(t, \chi) \simeq  {\cal A}_\chi(t) \, N^{1/2 +\epsilon}
\,\,\,\,\,\,\,\,
,
\,\,\,\,\,\,\,\,
N \rightarrow \infty \,\,\,
\label{importantbehaviour}
\eeq
for any positive $\epsilon >0$, where the prefactor ${\cal A}_\chi(t)$ may depend on the character $\chi$ and possibly also on $t$. Sums with a power law behavior such as $N^{\alpha}$ commonly occur in the displacement of random walks: the value of the exponent equal to $\alpha = \half $ corresponds to the purely {\em diffusive} brownian motion, those with values of $\alpha $ in the interval $ 0 < \alpha < \half$ correspond instead to {\em sub-diffusive} motion, while those in the interval $\half < \alpha < 1$ to {\em super-diffusive} motion, as for instance happens in Levy flights (see \cite{Mazo,Rudnick,Levy,diffusion}). 

\begin{figure}[t]
\centering
\includegraphics[width=0.60\textwidth]{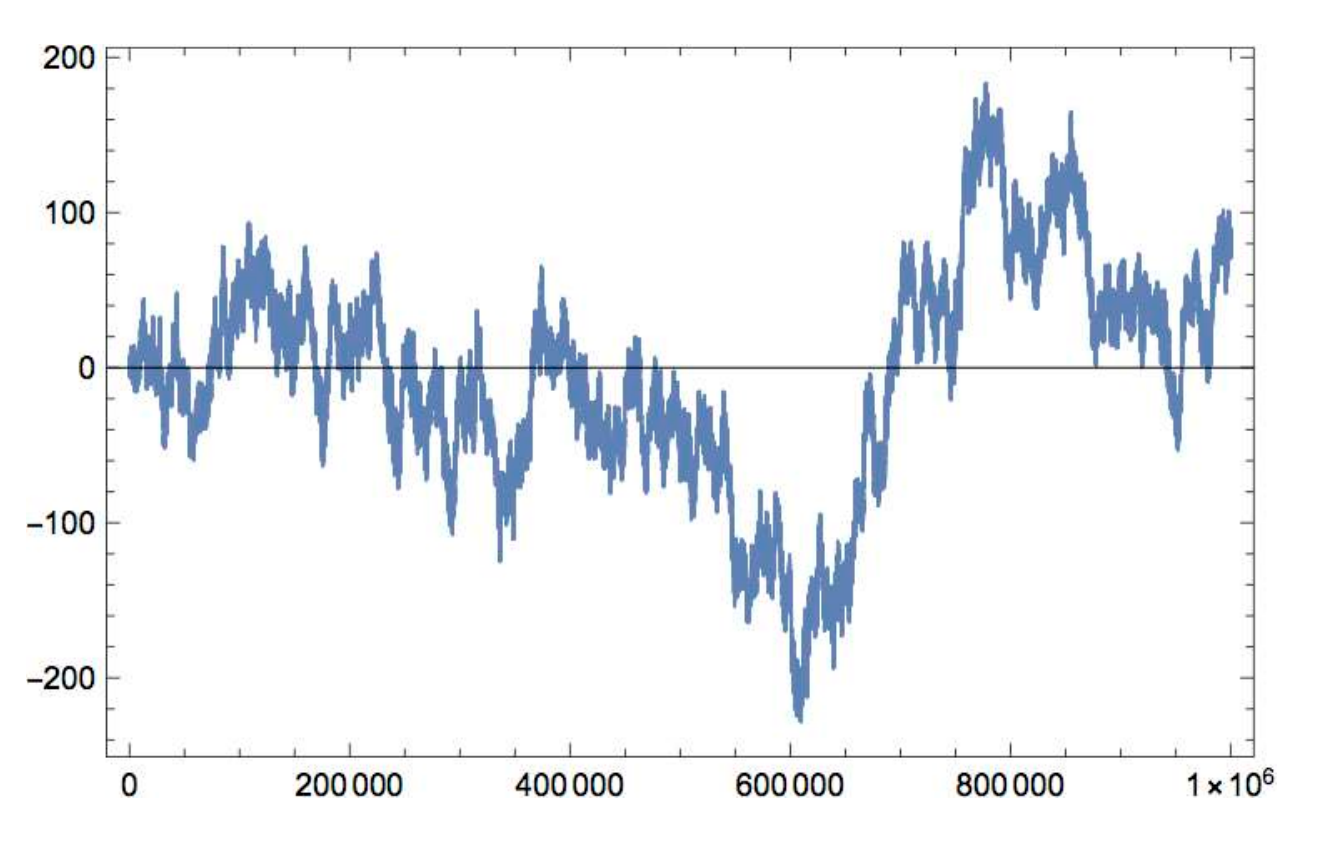}
\caption{Plot of $B_N(t)$ versus $N$ for $t=0$, with the angles $\theta_n$ given by the character $\chi_2$ mod 7. (See Table I for the values of this character).}
%\ref{tablech} .) }
\label{tipical}
\end{figure}

The functions $B_N(t, \chi)$ in eq.\,(\ref{Bprimes}) are of course purely deterministic series but, as we are going to show below, for all purposes they behave as random time series: their typical behavior varying $N$ (at a fixed value of $t$) is shown in Figure \ref{tipical} and one can see that their curve erratically fluctuates between positive and negative values with amplitudes which indeed increase (up to logarithmic corrections) as $N^{1/2}$. Random time series are ubiquitous quantities in science and we have taken advantage of several powerful methods which have been developed for their study (see, for instance, \cite{timeseries1,timeseries2}) to extract some relevant information on our series $B_N(t, \chi)$. 

As shown below, for  the aim of establishing the GRH  for non-principal characters, it is sufficient to study the behavior of these series at $t=0$, hereafter denoted as 
\beq
C_N \,=\, \sum_{n=1}^N \cos \theta_{p_{n}} \,\,\,.
\label{seriescn}
\eeq
These quantities only involve the sequence of angles $\theta_{p_n}$ of the characters computed at the primes $p_n$. A behavior as ${\cal O}(N^{1/2})$ of these series will guarantee the validity of the GRH for the corresponding Dirichlet functions. As originally shown in \cite{shortpaper}, such a universal diffusive behavior of the series $C_N$ comes from the Dirichlet theorem on the equidistribution of reduced residue classes modulo $q$ \cite{Diric,SelbergD} and the Lemke Oliver-Soundararajan conjecture on the distribution of pairs of residues on consecutive primes \cite{OliverSoundararajan}, and it is interesting to show how methods and statistical analysis related to time series help in corroborating this result. 

It is worth noticing (and we will comment more on this point later) that such a random walk process does not exist however for the $L$-functions associated to the principal characters, simply because all relative angles $\theta$'s for these characters vanish: this is what makes these functions a special case and it seems necessary to use a strategy different from the one presented in this paper in order to prove that their non-trivial zeros are all on the axis $\Re(s) = \half$. There is however a notable fact to take into account: thanks to the identity shown in the forthcoming eq.\,(\ref{identitylprincipal}), all $L$-functions based on principal characters share exactly the {\em same non-trivial zeros} of the Riemann $\zeta$ function. This means that the proof of the GRH for all $L$ functions relative to principal characters simply reduces to the proof of the original Riemann conjecture for the Riemann $\zeta$ function. An approach to this problem also based on random walk will be discussed somewhere else.

In the following the reader will find firstly the definition of the main quantities presented in this introduction and secondly the detailed discussion of the stochastic arguments which lead us to establish the asymptotic behavior (\ref{importantbehaviour}). In addition to some rigorous results, this paper also contains extensive numerical analysis as well as some heuristic arguments. In more detail, the paper is organized as follows.

\vspace{3mm}
\noindent
{\bf Contents of the paper}. Section \ref{CHH} is devoted to reviewing the main properties of  the Dirichlet $L$-functions, defined in terms of their infinite series or product representations, and to studying their analytic structure; we also discuss in certain detail the properties of the characters $\chi(n)$ which enter their definition. This Section can be skipped by readers already familiar with elementary aspects of analytic number theory.  Section \ref{surprising} presents two surprising results about the zeros of functions strictly related to the $L$-functions whose main outcome is to put in perspective the role of the primes in relation to the GRH. In Section \ref{extending} we introduce the quantities $B_N(t)$ that determine whether  one can extend the region of  convergence of the infinite product representation of the $L$-function. In Sections \ref{randomprimeensemble} and \ref{centrallimittheorem} we present the emergence of a normal law distribution for an analog of the $B_N(t, \chi)$ constructed on a set of ``random primes", whose notion is specified in the same Sections. Although this result points out the existence of a normal law distribution for the quantities we are interested in, for the aim of establishing the validity of the GRH this result is however inconclusive. Indeed, for that purpose, one needs also to control the asymptotic behavior of the mean entering this normal law: this point is the main content of the next Sections. In Section \ref{mean} we present a series of results that narrow down the behavior of the series $B_N(t, \chi)$: in particular, we show that the large $N$ behavior of the $B_N(t, \chi)$ at any $t$ is dictated by the large $N$ behavior of the series at $t=0$ and therefore only by the angles $\theta_{p_n}$. In Section \ref{timeseriess} we present some insights on the angles $\theta_{p_n}$ and the relative series $C_N$ defined in eq.\,(\ref{seriescn}) which come from adopting the point of view that $C_N$  is a random time series. Such an empirical study will find its theoretical framework in Section \ref{stprth}, where the statistical properties of the angles $\theta_{p_n}$ are nailed down on the basis of the Dirichlet theorem and the Lemke Oliver-Soundararajan conjecture. The growth of the series $C_N$  is the subject of Section \ref{Timeseries} where, based on the Dirichlet theorem and the Lemke Oliver-Soundararajan conjecture, we show that $C_N$ has a purely diffusive random walk behavior. 
Our conclusions are then discussed in Section \ref{conclusions}. 

The paper has also several Appendices. Appendix A discusses the statistical interpretation of the $L$-functions in terms of partition functions of free particles. Appendix B analyses the zeros and pole of the $L$-function in terms of the density of states of the equivalent statistical physics systems,  showing a connection with Fisher zeros.  Appendix C presents the proofs of two theorems concerning the $L$-functions, one due to Grosswald and Schnitzer, the other to Chernoff. Appendix D discusses the Kac's theorem relative to sums of trigonometric series with incommensurate frequencies and why this result cannot be used to prove the GRH. Appendix E shows how the large $N$ behavior of the series $B_N(t)$ is dictated by its value at $t=0$.

\vspace{3mm}
\noindent
{\bf Notation}. In the following the index $n$ is used {\em exclusively} in association with primes and denotes either the $n$th prime $p_n$ itself or quantities which depend on $p_n$, such as the angle $\theta_{p_n}$ defined by $\chi(p_n) \equiv e^{i \theta_{p_n}}$. Analogously, indices as $N$  stand for the upper limit of a sequence on primes. Hence any sum 
$\sum_{n=1}^N [...]$ 
on the index $n$ up to $N$ stands  for a sum involving the first $N$ primes. For this reason, we will use other indices, such as $m$ or $k$, etc. to denote the natural numbers and sums thereof.

\section{Dirichlet Characters and L-functions}\label{CHH}
This section collects all the basic results about L-functions \cite{Apostol,Iwaniec,Bombieri,Steuding,Sarnak}  we will need in the following 
and is designed to help the reader to follow the analysis we perform later. An informed reader can easily skip this section, since what  is presented here are well-known mathematical properties of the $L$-functions.    

\vspace{3mm}
\noindent
{\bf Arithmetical Progressions}. As the infinite series of odd numbers $1, 3, 5, \ldots (2 m+1)$ contains infinitely many primes, an interesting question to settle  is whether this property is also shared by other arithmetic progressions such as 
\beq
S_m \,=\, q \, m + h \,\,\,\,\,,\,\,\,\, m=0, 1, 2, \ldots 
\,\,\,\,\,\,\,\,\,\,\,\, 
q, h \in \mathbb{N} 
\label{sequence}
\eeq
In such progressions, the number $q$ is known as the {\em modulus} while the number $h$ as the {\em residue}. It is easy to see that a necessary condition to find a prime among the values of $S_m$ is that the two natural numbers $q$ and $h$ have no common divisors, namely they are {\em coprime}, a condition expressed as $(q,h) = 1$. In 1837 Dirichlet proved that this condition is also sufficient, that is if $(q,h) = 1$ then  the sequence  $S_m$ contains infinitely many primes. His ingenious proof involves some identities satisfied by functions defined by series expressions, known nowadays as Dirichlet $L$-functions which  generalize the more familiar Riemann $\zeta$(s) function \cite{Riemann1,Riemann2,Riemann3} with whom they share most of their analytic properties.

\vspace{3mm}
\noindent
{\bf L-functions: Infinite series definition}.
The Dirichlet L-functions of the complex variable $s = \sigma + i t$ are special cases of Dirichlet series: they are  
given by 
\beq
L(s,\chi) \,=\,\sum_{m=1}^\infty \frac{\chi(m)}{m^s} 
\,\,\,, 
\,\,\,\,\,\,\,\,
\Re (s) > 1 \,\,\,,
\label{Lfunctions}
\eeq 
where the arithmetic functions $\chi(m)$ are known as Dirichlet characters.  There are an infinite number of distinct Dirichlet characters, mainly characterized by their modulus $q$ which also determines their periodicity. As shown below, the non-zero characters are  complex numbers of modulus equal to $1$: hence, as any other Dirichlet series, the $L$ functions (\ref{Lfunctions}) are defined in an half-plane, here $\mathbb{R}(s) > 1$, where they converge absolutely.  These  $L$-functions can be then analytically continued  to $\Re (s) < 1$ using the functional relations presented in eq.\,(\ref{FELambda}) below.  Also, in  particular they can also be analytically continued into the so called {\em critical strip} $0 < \Re (s) < 1$ by certain integral representations. Thus they are analytic functions in the whole complex plane, except for the possibility of a pole at $s=1$.    
 
\vspace{3mm}
\noindent
{\bf Characters}. Let's discuss in more detail the characters $\chi$ entering the definition of the $L$-functions. To set the notation, given an integer $q$ we will denote by  the symbol $(q,a)$ the greatest common divisor of the two integers $m$ and $q$. If $(m,q) =1$, the two integers are said to be {\em coprime}. Given a modulus $q$, the prime residue classes modulo $q$ form an abelian group, denoted as   
\beq
(\mathbb{Z}/q\mathbb{Z})^* := \{m\, {\rm mod} \,q \,:\, (m, q) = 1\} \,\,\,. 
\label{groupab}
\eeq
The dimension of this group is given by the Euler totient arithmetic function $\varphi(q)$.  The latter  is defined to be the number of positive integers less than $q$ that are coprime to $q$.   Its value is given by 
\beq
\varphi(q) \,=\, q \, \prod_{p | q} \left(1 - \frac{1}{p}\right) \,\,\,,
\label{Eulertotient}
\eeq
where the product is over the distinct prime numbers dividing $q$. Notice that $\varphi(q)$ is an even integer number for $q \geq 3$. 

With  these definitions, a character $\chi$ of modulus $q$ is an arithmetic function from the finite abelian 
group $(\mathbb{Z}/q\mathbb{Z})^*$ onto $\mathbb{C}$ satisfying the following properties: 
\begin{enumerate}
\item 
$\chi(m+q) \,=\, \chi(m) $. 
\item 
$\chi(1) =1 $ and $\chi(0) = 0$. 
\item 
$\chi( n \, m ) \,=\, \chi(n) \, \chi(m)$. 
\item 
$\chi(m) = 0$ if $(m,q) > 1$ and $\chi(m) \neq 0$ if $(m,q) =1$. 
\item 
If $(m,q) =1$ then $(\chi(m))^{\varphi(q)} =1$, namely $\chi(m)$ have to be $\varphi(q)$-roots of unity. 
\item 
If $\chi$ is a Dirichlet character so is its  complex conjugate $\overline\chi$. 
\end{enumerate}
From property $5$, it follows that for a given modulus $q$ there are $\varphi(q)$ distinct Dirichlet characters that can be labeled as $\chi_{j}$ where $j = 1, 2, . . . , \varphi(q)$ denotes an arbitrary ordering.   We will not  display the arbitrary index $j$ in $\chi_j$,  except  for  explicit examples.  
Moreover, the characters satisfy the following orthogonality conditions 
\begin{eqnarray}
\sum_{r=1}^{\varphi(q)} \chi_r(k) \overline{\chi}_r(l) &\,=\, & 
\left\{ 
\begin{array}{cll}
\varphi(q) & & {\rm if} \,\,k\equiv l \,\,\, ({\rm mod} \, q) \\
0 & & {\rm if} \,\, k\not\equiv l \,\,\, ({\rm mod} \, q) 
\end{array}
\right. \\
\sum_{m=1}^{q} \chi_r(m) \overline{\chi}_s(m) & \,=\,& 
\varphi(q) \, \delta_{r,s} \,\,\,.
\end{eqnarray}

For a generic $q$, the {\em principal} character, usually denoted $\chi_1$, is defined as 
\beq
\chi_1(m) \,=\, 
\left\{ 
\begin{array}{cl}
 1 & \, {\rm if} \,\,(m, q) = 1 \\
 0 & \, \, {\rm otherwise} 
 \end{array}
 \right.
 \label{principalcharacter}
 \eeq
When $q = 1$, we have only the {\em trivial}  principal character $\chi(m) = 1$ for every $m$, and in this case the corresponding $L$-function reduces to the Riemann $\zeta$-function given by 
\beq
\zeta(s) \,=\,\sum_{m=1}^\infty \frac{1}{m^s} 
\,\,\,, 
\,\,\,\,\,\,\,\,
\Re (s) > 1 \,\,\,. 
\label{zetafunctions}
\eeq 
There is an important difference between principal versus non-principal characters. The principal characters, being only $1$ or $0$, satisfy 
\beq
\sum_{m=1}^{q-1} \chi_1(m) \,=\,\varphi(q) \,\neq \,0 \,\,\,,
\eeq
whereas  the non-principal characters  satisfy
\beq
\sum_{m=1}^{q-1} \chi(m) \,=\,0 \,\,\,.
\label{sumtozero}
\eeq
We will see below that these conditions determine the analytic structure of the $L$-functions.

\vspace{3mm}
\noindent
{\bf Parametrization of the angles}. Posing 
\beq
\label{thetan2}
\chi (m) \,=\,  e^{i \theta_m }, ~~~~~\forall ~\chi(m) \neq 0 \,\,\,, 
\eeq   
eq.\,(\ref{sumtozero}) shows that the angles $\theta_m$ of the non-principal characters defined in eq.\,(\ref{thetan}) are equally spaced over the unit circle. Since they are associated to the $\varphi(q)$ roots of unity, their possible values can be labelled as 
\beq
\label{notationangles}
\alpha_k \,=\, \frac{\pi (2 k - \varphi(q))}{\varphi(q)} \,\,\,\,\,\,\,\,\,\,\,\,\,, 
\,\,\,\,\,\,\,\,\,\,\,\, k = 1,\ldots, 
\varphi(q) \,\,\,.
\eeq 
In this parameterization, the angles $\alpha_k$, which are negative for $k=1,\ldots\,\varphi(q)/2$ and positive for $k=\varphi(q)/2+1,\ldots,\varphi(q)$, are related pairwise as 
\beq
\alpha_{\varphi/2 +k} \,=\,\alpha_k + \pi \,\,\,\,\,\,
,\,\,\,\, 
k=1,\ldots ,\varphi(q)/2 \,\,\,.
\label{pairwiseangles}
\eeq
Notice that the actual roots of unity entering the expression of the characters may be a smaller set of the $\varphi(q)$-roots of unity, 
\beq
\label{PhiSet}
\theta_m \in \Phi = \{ \phi_1,\phi_2,\ldots,\phi_r\} ~ {\rm with} ~  r \leq \varphi(q)
\eeq
 and  
$\phi_i$ equal to one of the angles of the set (\ref{notationangles}) (see the examples below). 
The integer $r$ is referred to as the {\em order} of the particular character and it is an integer that divides $\varphi (q)$.

\vspace{3mm}
\noindent
{\bf Primitive and non-primitive characters}. 
For values of $m$ coprime with $q$, the character $\chi(m)$ mod $q$ may have a period less than $q$. If this is the case, $\hat\chi$ will be called a {\em non-primitive} character, otherwise  $\chi$ is {\em primitive}. Obviously if $q$ is a prime number, then every character mod $q$ is primitive. If $\chi$ is a primitive character of modulus $ q$ and $q$  divides $\hat q$, then we can construct a character $\hat\chi(m)$ mod $\hat  q$ in  the following manner:
\beq
\hat\chi(m) \,=\, \left\{
\begin{array}{cll} 
\chi(m) & {\rm if} & (m,\hat q) = 1 \\
0 & {\rm if} & (m,\hat q) > 1 
\end{array}
\right.
\eeq
In this case, we say that the character $\hat\chi$ mod $q$ is induced by $ \chi$.\footnote{${\chi}$ is called the {\it conductor} of $\hat\chi$.}  It is important to stress that every non-primitive character is induced by a primitive one. For this reason from now on we focus our attention on the primitive characters only. 

\begin{figure}[t]
\begin{center}
%\vspace{8mm}
\includegraphics[width=6.5cm]{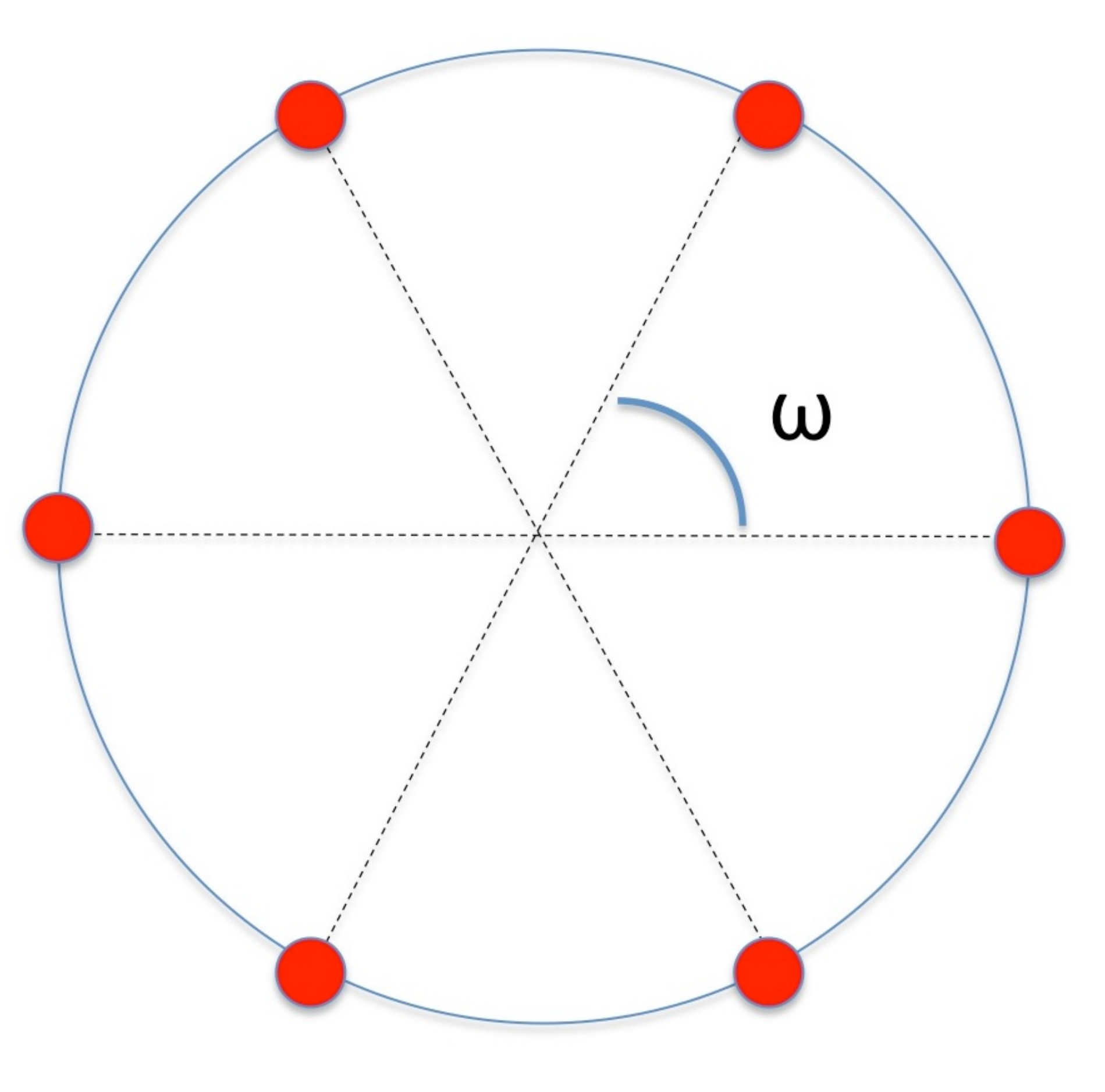}
\caption{$6$-th roots of unity associated to all characters mod 7.}
\label{character7}
\end{center}
\end{figure}

\vspace{3mm}
\noindent
{\bf Explicit examples of primitive characters}. As explicit examples, in the following we will mainly use those associated to the modules $q=3$, $q=5$ and $q=7$. Since these are prime numbers, all their relative characters are primitive. For $q=3$, they are expressed in terms of the square roots of $1$, for $q=5$ in terms of the $4$-th roots of unity, while for $q=7$ in terms of the $6$-th roots of unity shown in Figure \ref{character7}, with $\omega = e^{i\pi/3}$. The set of characters for these modules are shown in Table \ref{tablech}.  Notice that for 
$q=3$ all characters are real, for $q=5$ $\chi_1$ and $\chi_3$ are real (and the corresponding angles belong to a smaller set of the $4$-roots of unit) while the terms of the pair $(\chi_2,\chi_4)$ are complex conjugates of each other. In the case of $q=7$, $\chi_1$ and $\chi_4$ are real (and the corresponding angles belong to a smaller set of the $6$-roots of unit) while the terms of the pairs $(\chi_2,\chi_6)$ and $(\chi_3,\chi_5)$ are complex conjugates of each other. The characters $(\chi_3,\chi_5)$ are composed of $r=3$ subsets of the angles (\ref{notationangles}) (i.e. the angles $0, \pm 2\pi/3$), while the characters $(\chi_2,\chi_6)$ employ the full set of angles (\ref{notationangles}), with $\varphi(7) = 6$. 

\begin{table}[t]
\hspace{-15mm}
\begin{minipage}{.22\textwidth}
\centering
\begin{tabular}{|c | c | c | c |}
\hline\hline
n & 1 & 2 & 3 \\ [0.5ex] % inserts table %heading
\hline
$\chi_1(n)$ 
& \,\, 1 \,\, & \,\,\,\,1 \, & \,\, 0 \,\,\\
\hline
$\chi_2(n)$ & 1 & - 1 & \,\, 0 \,\,\\
\hline
\end{tabular}
\end{minipage}
\begin{minipage}{.36\textwidth}
\centering
\begin{tabular}{|c | c | c | c | c | c |}
\hline\hline
n & 1 & 2 & 3 & 4 & 5 \\ [0.5ex] % inserts table %heading
\hline
$\chi_1(n)$ 
& \,\, 1 \,\, & \,\,\,\,1 \, &  \,\,\,\,1 \, &  \,\,\,\,1 \,& \,\, 0 \,\,\\
\hline
$\chi_2(n)$ & 1 & i  & - i & -1 & \,\, 0 \,\,\\
\hline
$\chi_3(n)$ & 1 & -1 & -1  & 1 & \,\, 0 \,\,\\
\hline
$\chi_4(n)$ & 1 & -i  &  i & -1 & \,\, 0 \,\,\\
\hline
\end{tabular}
\end{minipage}
\begin{minipage}{.31\textwidth}
\centering
\begin{tabular}{| c | c | c | c | c | c | c | c |}
\hline\hline
n & 1 & 2 & 3 & 4 & 5 & 6 & 7 \\ [0.5ex] % inserts table %heading
\hline
$\chi_1(n)$ 
& \,\, 1 \,\,& 1 & 1 & 1 & 1 & \,\,1 \,\,& \,\,\, 0 \,\,\\
\hline
$\chi_2(n)$ & 1 & $\omega^2$ & $\omega$ & $-\omega$ & $-\omega^2$ & \,\,-1 \,\,& \,\,\, 0 
\,\,\\
\hline
$\chi_3(n)$ & 1 & $- \omega$ & $\omega^2$ & $\omega^2$ & $-\omega$ & \,\, 1 \,\,& \,\,\, 0 \,\,\\
\hline
$\chi_4(n)$ & 1 & 1 & -1 & 1 & -1 & \,\,-1 \,\,& \,\,\, 0 \,\,\\
\hline
$\chi_5(n)$ & 1 & $\omega^2$ & $-\omega$ & $-\omega$ & $\omega^2$ &\,\, 1 \,\, & \,\,\, 0 \,\,\\
\hline
$\chi_6(n)$ & 1 & $-\omega$ & $-\omega^2$ & $\omega^2$ & $\omega$ & \,\,-1 \,\, & \,\,\, 0 \,\, \\
\hline
\end{tabular}
\end{minipage}
\caption{Left hand side: characters mod $q=3$. Center: characters mod $q=5$. Right hand side: characters mod $q=7$, where 
$\omega = e^{i \pi/3}$. }
\label{tablech}
%\label{table:nonlin}
\end{table}

\vspace{3mm}
\noindent
{\bf Infinite product representation}. Due to the completely multiplicative property of the characters, the $L$-functions can also be  expressed in terms of the infinite product representation recalled in the Introduction  
\beq
\label{EPF}  
L(s, \chi)   \,=\, 
 \prod_{n=1}^\infty  \( 1 -  \frac{\chi (p_n)}{p_n^s} \)^{-1} \,\,\,,     ~~~~ \Re (s) > 1,
\eeq
where $p_n$ is the $n$-th prime in ascending order. This infinite product is certainly convergent for $\Re (s) >1$ (and it coincides with the series representation of the $L$-function which also converges in this domain), but it may have a larger domain of convergence. Recall the main goal of this paper is indeed to show that its abscissa of convergence can be safely extended down to $\Re (s) = \half$ 
for non-principal characters.   

Notice that if $ \chi$ is a primitive character mod $ q$ which induces another character $\hat\chi$ mod $\hat q$, we have 
\beq
L(s,\hat\chi) \,=\,L(s, \chi) \, \prod_{p | {\hat q} } \left(1 - \frac{\chi(p)}{p^s}\right) 
\,\,\,,
\label{primitivevsnonLfunctions}
\eeq
where the product extends to the {\em finite} set of primes $p$ which divide $\hat q$. Hence, every $L$-function is equal to the $L$-function $L(s,\chi)$ of a primitive character multiplied by a finite number of terms.  In fact,  the above formula shows that $L(s, \chi)$ and $L(s, \hat{\chi})$ share the same non-trivial zeros.  Therefore, for the purpose of  establishing  whether the infinite product representation (\ref{EPF}) of the $L$-function can be extended to $\Re(s)> \half$,  it is sufficient to focus our  attention only on the $L$-functions based on  primitive characters. 

\vspace{3mm}
\noindent
{\bf $L$-functions of principal characters and Riemann $\zeta$ function}. Notice that the principal character of modulus $q$ satisfies 
eq.\,(\ref{principalcharacter}) and therefore the relative $L$-functions can be expressed as 
\beq
L(s,\chi_1) \,=\,\prod_{p \nmid q} \left(1 - \frac{1}{p^s}\right)^{-1} \,=\, \zeta(s) \, \prod_{p\mid q} 
\left(1 - \frac{1}{p^s}\right)^{-1}\,\,\,,
\label{identitylprincipal}
\eeq
where $d | n$ denotes the integer $d$ which divides the integer $n$, and $d \nmid n$ otherwise. Since the finite product involving the primes which divide $q$ in the the right hand side never vanish, the zeros of the Dirichlet $L$-functions of principal characters coincide exactly with the zeros of the Riemann $\zeta$ function. Hence, establishing the GRH for these functions is equivalent to 
prove the original hypothesis by Riemann for the $\zeta$ function.

\vspace{3mm}
\noindent
{\bf Functional equation}. The $L$-functions associated to the primitive characters satisfy a functional equation similar to  that of the Riemann $\zeta$-function. To express such a functional equation, let's define the index $a$ as 
\beq\label{order}
a \equiv \begin{cases}1 \qquad
\mbox{if $\chi(-1)= -1$ ~~\,\,\,\,(odd)} \\
0 \qquad \mbox{if $\chi(-1) = \,\,\,\,\,1$  ~~~~(even)}
\end{cases}
\eeq
Moreover let's also introduce the Gauss sum 
\beq
G(\chi) \,=\,\sum_{m=1}^q \chi(m) \, e^{2 \pi i m/q}
\,\,\,,
\label{Gausssum}
\eeq
which satisfies $|G(\chi)|^2 = q$ if and only if the character $\chi$ is primitive. With these definitions, the functional equation for the primitive $L$-functions can be written as 
\beq
\label{FELambda}
L(1-s, \chi) \,=\, i^{-a} \,\frac{q^{s-1} \, \Gamma(s)}{(2\pi)^s} 
\, G(\chi) \, \left\{\begin{array}{c}\cos(\pi s/2)\\ \sin(\pi s/2) 
\end{array}\right\} \, 
L(s, \overline\chi)\,\,\,. 
\eeq
where the choice of cosine or sine depends upon the sign of $\chi(-1) = \pm 1$.  An equivalent but a more symmetric version of the functional equation (\ref{FELambda}) can be given in terms of the so-called {\em completed $L$-function} $\hat L(s,\chi)$ defined by  
\beq
\hat L(s,\chi) \equiv \left(\frac{q}{\pi}\right)^{(s+\delta)/2} \, 
\Gamma\left(\frac{s+\delta}{2}\right) \, L(s,\chi)\,\,, 
\label{completedLfunction}
\eeq
where $\delta = \frac{1}{2} (1 - \chi(-1))$. The completed $L$-function satisfies the functional equation 
\beq
\hat L(s,\chi) \,=\, \epsilon(\chi) \, \hat L(1-s, \overline\chi) \,\,\,,
\label{equivFunctional}
\eeq
where the quantity $\epsilon(\chi)$ 
\beq  
\epsilon(\chi) \,=\, \frac{G(\chi)}{i^\delta \sqrt{q}} \,\,\,,
\label{epsilonchi}
\eeq
is a constant of absolute value $1$.

\vspace{3mm}
\noindent
{\bf Analytic structure of the $L$-functions}.  As previously mentioned, there is an important distinction between the $L$-functions based on non-principal verses  principal characters which will be very
important for our purposes.  
\begin{itemize}
\item
The $L$ functions for  non-principal characters are {\em entire} functions,  i.e. analytic everywhere in the complex plane with no poles. 
\item The $L$-functions  $L(s,\chi_1)$ for  principal characters, on the contrary, are analytic everywhere except for a {\em simple pole} at $s=1$ with residue $\varphi(q)/q$. 
\end{itemize}
To show this result, let us  first  express any $L$-function in terms of a {\em finite} linear combination of the Hurwitz zeta function defined by the series  
\beq
\zeta(s,a) \,=\, \sum_{m=0}^\infty \frac{1}{(m+a)^s} \,\,\,,
\label{Hurwitz}
\eeq
whose domain of convergence is $\Re (s) >1$. Since we can split any integer $m$  as  
$$
m \,=\, q\, k +r \,\,\,\,\,\,\,,\,\,\,\,\,\,\,\,
{\rm where} \,\, 1 \leq r \leq q \,\,\,
{\rm and} \,\,\, k=0,1,2, \ldots 
$$
we have 
\begin{eqnarray}
L(s,\chi) &\,=\,& \sum_{m=1}^\infty \frac{\chi(m)}{m^s} \,=\, 
\sum_{r=1}^q \sum_{k=0}^\infty \frac{\chi(q k + r)}{(q k +r)^s} \,=\, 
\frac{1}{q^s} \, \sum_{r=1}^q \chi(r) \, \sum_{k=0}^\infty 
\frac{1}{\left(k+\frac{r}{q}\right)^s} \\
&=& \frac{1}{q^s} \, \sum_{r=1}^q \chi(r) \, \zeta\left(s, \frac{r}{q}\right) 
\,\,\,.\nonumber
\label{LHurwitz}
\end{eqnarray} 
The Hurwitz $\zeta$-function has a simple pole at $s=1$ with residue 1 and therefore the residue at this pole of the $L$-function is 
\beq
{\rm Res} \, L(s,\chi) \,=\, \frac{1}{q} \sum_{r=0}^q \chi(r)\,=\, 
\left\{
\begin{array}{cll}
\frac{\varphi(q)}{q} & & {\rm if }\, \, \chi = \chi_1 \\
0 & & {\rm if }\,\, \chi \neq \chi_1 \,\,\,.
\end{array}
\right.
\eeq

\vspace{3mm}
\noindent
{\bf Trivial Zeros}. Using the Euler product representation of the $L$-function it is easy to see that these functions have no zeros in the half-plane $\Re (s) > 1$, in particular $\log L(s,\chi)$ is finite in this region since the series converges there.  Examining the functional equation (\ref{FELambda}) one sees that, analogously to the Riemann $\zeta$-function, the trivial zeros of the $L$-functions are those in correspondence with the zeros of the trigonometric functions present in the expression. Therefore 
\begin{enumerate}
\item If $\chi(-1) = 1$, then the trivial zeros are along the negative real axis located at $\sigma = - 2 k$, with $k=0,1,2,\ldots$. 
\item If $\chi(-1) = -1$, then the trivial zeros are along the negative real axis  but  now located at $\sigma = - 2 k -1$, with $k=0,1,2,\ldots$.
\end{enumerate} 

\vspace{3mm}
\noindent
{\bf Non-trivial Zeros and Generalized Riemann Hypothesis}.
All other non-trivial zeros must lie in the critical strip $ 0 < \sigma < 1$. When the character is real, if $\rho = \sigma + i t$ is a zero of $L(s,\chi)$ then $\hat\rho = (1- \sigma) - i t$ is also a zero of the same $L$-function and, if $\sigma = 1/2$, the two zeros are then complex conjugates of  each other.  When the character $\chi$ is instead complex, a zero $\rho = \sigma + i t$ of $L(s,\chi)$ corresponds to a zero $\hat\rho = (1 - \sigma) - i t$ of $L(s,\overline\chi)$: in this case, if $\sigma =1/2$, the zeros of the $L$-functions associated to complex characters are not necessarily complex conjugates.   

According to the Generalized Riemann Hypothesis, all non-trivial zeros of the primitive\footnote{It is important to refer to primitive characters in order to exclude the zeros of the factors $\prod_{p | {\hat q}} \left(1 - \chi(p) \,p^{-s}\right)$ present in the non-primitive characters, see eq.\,(\ref{primitivevsnonLfunctions}), which are all along the line $\sigma = 0$.} $L$-functions lie on the critical line $ \sigma = \half$.   An explicit formula for the $n-th$ zero  as the solution of a transcendental equation was proposed in \cite{Transcendental}.

\vspace{3mm}
\noindent
{\bf  Our approach to the GRH}. 
Having completed in this section the review of known facts of the $L$-functions, it is worthwhile restating the approach to the GRH that we are pursuing here. This is based on the following observation: if  the Euler product formula were valid for $\Re (s) > \half $, i.e. a domain larger than the original one $\Re (s) > 1 $ stated in eq.\,(\ref{EPF}), then the GRH would follow by very simple arguments. Namely,  it would establish that there are no zeros with $\Re (s) > \half$. Combined with the functional equation,  this implies  there are no zeros with $\Re (s) < \half$. Thus, all non-trivial zeros have to be on the critical line $\Re (s) = \half$. It is known that they are infinite in number  since the number of them with ordinate $0< t< T$ is known to leading order as\footnote{This result holds for $L$-functions relative to 
primitive characters mod $q$.}   
\beq
N(T,\chi) \,=\,\frac{T}{\pi} \,\log\frac{q T}{2 \pi e} + {\mathcal O}(\log q T) 
\,\,\,.
\label{asymptoticzeross}
\eeq
Hence, in the next sections we are going to study the possibility to extend the infinite product representation of the $L$-functions from the original region $\sigma > 1$ to the new region $\sigma > \half$. Before doing this, it is however interesting to present in the meantime two remarkable results which unveil the crucial role played by the fluctuations of the primes in determining the zeros of the $L$-functions.

\section{Two surprising results}\label{surprising}
There are two surprising, but in a sense opposing, results concerning the zeros of both the Riemann $\zeta$ function and all other Dirichlet  $L$-functions. These results are the content of the following two theorems.

\begin{theorem} (Grosswald and Schnitzer) 
\label{GrosswaldSchnitzer} \cite{GrosswaldSchnitzer}. 
%{\bf Theorem 2} ({\em Grosswald and Schnitzer}) \cite{GrosswaldSchnitzer}.   
Let $L(s,\chi)$ be the Dirichlet $L$ function based on any  Dirichlet character of modulus 
$q$.  Let $\Primes = \{p_1, p_2, \ldots \}$  denote the set of primes while $\Primesp = \{ p'_1, p'_2, \ldots \}$ a set of  integers $p'_n$ satisfying 
\beq
 \label{ppn}
 p_n \leq p'_n < p_n + K, ~~~~~~p'_n = p_n ~~ ({\rm mod} ~ q)
\eeq
where $K \geq q$ is an arbitrary integer,  and define the modified $L$-function according to the infinite product 
\beq
\label{Lp}
L' (s, \chi) \,=\,  
\prod_{n=1}^\infty  \( 1 -  \frac{\chi(p'_n)}{(p'_n)^s} \)^{-1} \,\,\,.
\eeq
Then $L'(s,\chi) $ can be analytically continued to the half plane $\sigma > 0$ and in this domain it has the same zeros as the Dirichlet $L$-function $L(s,\chi)$. Moreover, if $\chi$ is a non-principal character then $L'(s, \chi)$ has no poles for $\sigma>0$, as does $L(s,\chi)$. 
\end{theorem}
\begin{theorem} ({\em Chernoff}) \label{Chernoff} \cite{Chernoff}.  
%{\bf Theorem 1} ({\em Chernoff}) \cite{Chernoff}.  
Consider the Euler infinite product representation of the Riemann $\zeta$-function. Substitute the primes $p_n$ in such a formula with their smooth approximation $p_n \sim n \log n$ and define the modified  function $\zeta'(s)$ 
according to the infinite product 
\beq 
\zeta'(s) \,=\,\prod_{n=1}^\infty  \( 1 -  \frac{1}{(n\,\log n)^s} \)^{-1}\,\,\,.
\label{CsC}
\eeq 
The  function $\zeta'(s)$can be analytically continued into the half-plane $\Re(s) > 0$ except for an isolated singularity at $s=0$.  Furthermore it no longer has any zeros in this region. 
\end{theorem}
%\item
\vspace{1mm}

\begin{figure}[t]
\centering\includegraphics[width=.5\textwidth]{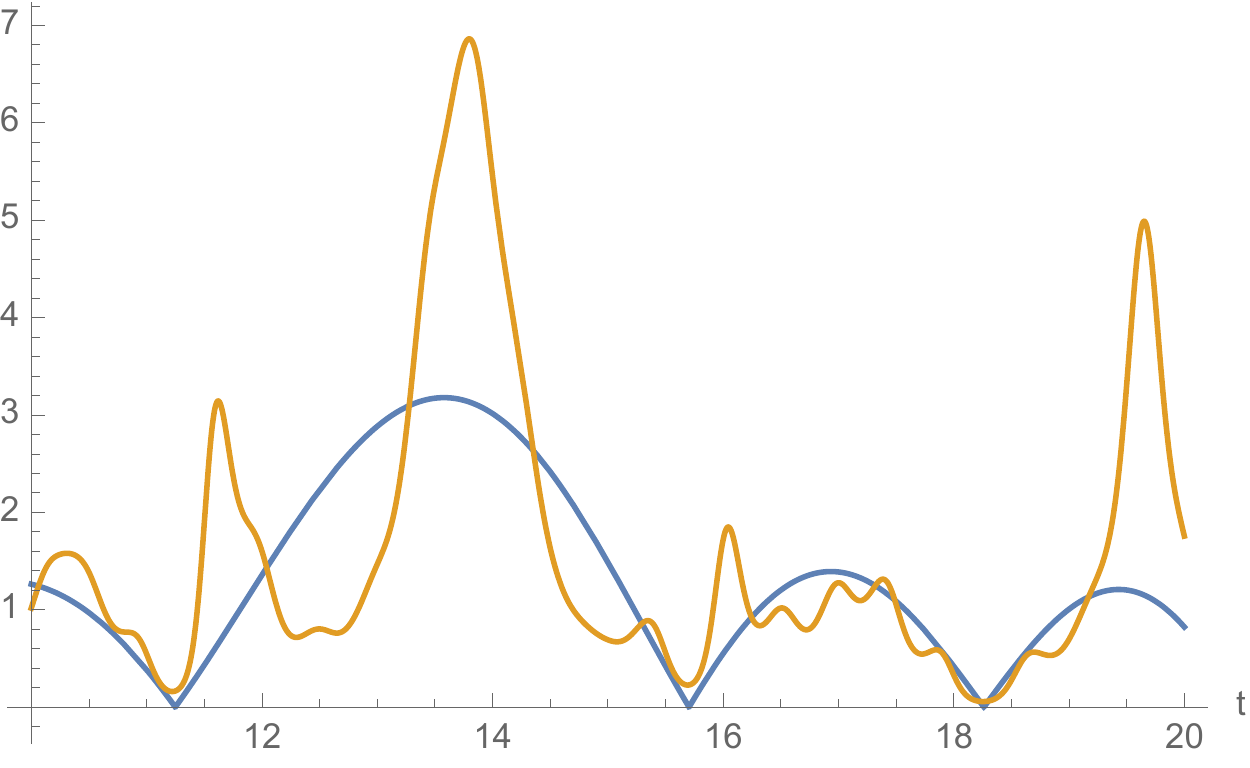}
\caption{Numerical illustration of Theorem \ref{GrosswaldSchnitzer} for the non-principal character $\chi_2$ mod $3$. The blue line is $|L(\half + i t)|$ where $t$ is the x-axis.  
The erratic orange line is $|L'(\half + i t)|$ for some randomly chosen state $\Primesp$.  We chose $N =10^6$ and $K=Mq=6$  for $M =2$.}
\label{LpDir}
\end{figure}

\noindent
What is surprising is the emergence of the following scenario: if in the infinite Euler product we use the smooth approximation of the prime numbers $p_n \simeq n \,\log n$, then all the non-trivial zeros of the Riemann $\zeta$-function in the critical strip completely disappear. On the contrary, if in the  Euler product we use another set of random numbers which shares with the primes the same modularity property and the same rate of growth, then all the non-trivial zeros of the original $L$-function remain  exactly at  the same location  in the critical strip!  In particular, Theorem \ref{GrosswaldSchnitzer} suggests that the validity of the Generalized Riemann Hypothesis may not depend on the detailed properties of the primes and this further justifies the probabilistic considerations presented later in this paper. A numerical check of Theorem \ref{GrosswaldSchnitzer} can be found in Figure \ref{LpDir}, where we have chosen the non-principal character $\chi_2$ mod $3$, whose values in the first period are given by 
\beq
\label{char}
\{ \chi (1), \chi(2),  \chi (3) \} = \{ 1, -1, 0 \} \,\,\,, 
\eeq
to plot $|L(\half+it)|$ and $|L'(\half+i t)|$  for a randomly chosen set of the integers $p_n'$ as a function of $t$ in the region of  the first 3 zeros. Whereas $|L'(\half+i t)|$ is erratic due to the randomness of the integers $p'_n$ and changes its shape if we change the set of these random numbers, the  validity of Theorem \ref{GrosswaldSchnitzer} is nevertheless clear, i.e.  the two functions share the same  zeros.  An interesting aspect of this plot concerns how we calculated $|L'(\half+i t)|$: we did not formally analytically continue it into the entire critical strip, since it is unknown how to do so numerically.      Rather   we used  the Euler product to continue it only to the right of the critical line, which is sufficient for our purposes. In short,  Figure \ref{LpDir} provides numerical evidence that the Euler product converges to the right of the critical line for $L'$, which is the key idea  we are going to address in the next sections. A short proof of both theorems is presented  in Appendix \ref{proofstheorems}, while we refer the reader to the original references for a more detailed discussion.

\section{Infinite product into the critical strip}\label{extending}
\label{Probabilistic}

The aim of this Section is to present a criterion which allows us to extend the region of convergence of the Euler product of the $L$-functions and to constrain the location of their zeros. The main result was proven in \cite{EPFchi}, and here we summarize it also providing additional relevant remarks. From now on, unless stated explicitly, we focus our attention only on $L$-functions corresponding to non-principal characters.  Consider the infinite product representation of the $L$-functions 
\beq 
L(s,\chi) \,=\, 
 \prod_{n=1}^\infty  \( 1 -  \frac{\chi (p_n)}{p_n^s} \)^{-1} \,\,\,,
\eeq
and take the formal logarithm on both sides of this equation, so that 
\beq 
\log L(s,\chi) \,=\,  X(s,\chi) + R(s,\chi)\,\,\,,
\eeq
where
\beq\label{PDir}
X(s,\chi) = \sum_{n=1}^{\infty} \dfrac{\chi(p_n)}{p_n^{\, s}}\,\,\,\,\,\,\, , \qquad
R(s,\chi) = \sum_{n=1}^{\infty} \sum_{m=2}^{\infty}
\dfrac{\chi(p_n)^{m}}{m p_n^{\, ms}}\,\,\,.
\eeq
Now $R(s,\chi)$ absolutely converges for $\sigma > \half$, so we can write
\beq\label{logsum2}
\log L(s,\chi) = X(s,\chi) + O(1)
\eeq
which indicates that the convergence of the Euler product to the right of the critical line depends only on properties of $X(s,\chi)$. The singularities of $ \log L(s,\chi)$ are determined by the zeros of $L(s,\chi)$ and, if present, also by the pole $s=1$. For what concerns the GRH, the main emphasis is of course in locating the {\em zeros} of these functions and the eventual presence of the pole at $s=1$ is a simple, though significant, complication\footnote{See the discussion below, after eq.\,(\ref{noshiftfunction}), for the effect induced by the pole at $s=1$. This is relevant for the Riemann Hypothesis relative to the $\zeta$ function.}. The advantage of the $L$-functions of non-principal characters is that they do not have a pole at $s=1$ and therefore for all of them we have a very concise  mathematical statement: $X(s,\chi)$ {\em is the diagnostic quantity which directly locates their non-trivial zeros}. Taking now the real part\footnote{Analogous arguments apply to the imaginary part of $X(s, \chi)$.}
of $X(s, \chi)$ in \eqref{PDir}, we have 
\begin{equation}
\label{SDef}
S(\sigma,t,  \chi) \,=\, \sum_{n=1}^\infty \dfrac{\cos(t \log p_n -  \theta_{p_n})}{p_n^{\, \sigma}} \,\,\,. 
\end{equation}
A further elaboration of this expression goes as follows. 
Defining
\begin{equation}
\label{CxDef}
B(x ; t, \chi) \,=\, \sum_{p \le x} \cos\( t \log p -  \theta_p\)\,\,\,,  
\end{equation}
we have  
\beq
B(p_n ; t, \chi) - B(p_{n-1} ; t, \chi) \,=\, \cos (t \log p_n - \theta_{p_n})\,\,\,,
\eeq
and then  
\begin{equation}
S(\sigma, t, \chi) = \sum_{n=1}^{\infty} B(p_n ; t,  \chi) \left( \dfrac{1}{p_n^{\,\sigma}} - 
\dfrac{1}{p_{n+1}^{\,\sigma}}\right) = \sigma \sum_{n=1}^{\infty} 
B(p_n ; t, \chi) \int_{p_n}^{p_{n+1}} \dfrac{1}{u^{\sigma+1}} du\,\,\,.
\end{equation}
Given that $B(x ; t,  \chi) = B(p_n ; t,  \chi)$ is a constant for $x \in (p_n, p_{n+1})$, we finally arrive to 
\begin{equation}
\label{importantintegral}
S(\sigma, t, \chi) = \sigma \int_{2}^{\infty} \dfrac{B(x ; t,  \chi)}{x^{\sigma+1}} dx\,\,\,. 
\end{equation}
Hence, the convergence of the integral is fixed by the behavior of the function $B(x ; t,  \chi)$ at $x \rightarrow \infty$: if $B(x ; t,  \chi) = O(x^{\alpha})$ for $x \rightarrow \infty$ and for any $t$, then the integral converges for $ \sigma > \alpha$ and diverges precisely at $\sigma =\alpha$. All this is the content the following theorem:
\begin{theorem}
\label{BNtheorem}
(Fran\c ca-LeClair)~\cite{EPFchi}. Defining 
\beq
\label{BNt} 
\CBN  (t, \chi) \, =\, \sum_{\substack{n=1 \\ p_n \nmid q}}^{N}    \cos \(   t \log p_n  - \theta_{p_n} \)  \,\,\,, 
\eeq
if, for all $t$, $\CBN = O(N^{1/2+\epsilon})$  (up to logarithms, see below), then the Euler product formula is valid for $\Re (s) > \half  +\epsilon$ because it converges in this region. This implies there are no zeros with $\Re (s) > \half + \epsilon$.   
 \end{theorem}
   
 \vspace{3mm}
 
The above theorem implies that if $\CBN = O(N^{1/2+\epsilon})$ for all  $\epsilon > 0$,  up to logarithms, then the Generalized Riemann Hypothesis is true. By ``up to logarithms" we are referring to  factors involving $\log N$ or $\log \log N$ etc,  which do not spoil the convergence argument.  For instance behaviors as $\CBN =O(\sqrt{N \log \log N})$ (suggested by the law of iterated logs)  or $\CBN = O( \sqrt{N} \log^a N )$ for any positive power $a$ will be sufficient.  For the latter, assuming
$B(x) = O\left(\sqrt{x} \log^a x \right)$ yields
\begin{equation}
|S(\sigma, t,  \chi)| \le  K \int_{1}^\infty 
\dfrac{\log^a x}{x^{\sigma+1/2}} dx = 
K \dfrac{\Gamma(a+1)}{(\sigma - 1/2)^{a+1}}
\end{equation}
{\it Henceforth in places we  will simply write $\CBN = O(\sqrt{N})$ without always displaying the $\epsilon$, and it is implicit that this can be relaxed with such logarithmic factors. }

\vspace{1mm}

We can now see immediately what the problem is with the principal characters $\chi_{1}$: in this case, in fact, all the angles $\theta_{p_n}$ in the expression of $B_N(t)$ are zero and therefore we have
\beq 
\CBN  (t, \chi_{1}) \, =\, \sum_{\substack{n=1 \\ p_n \nmid q}}^{N}    \cos \(   t \log p_n  \)  \,\,\,. 
\label{noshiftfunction}
\eeq
Clearly, for this series, there is one special value of $t$, i.e. $t=0$, for which in the limit $N \rightarrow \infty$ this series diverges {\em linearly} in $N$ since  $\CBN  (0, \chi_{1}) = N$. The Mellin transform (\ref{importantintegral}) now diverges at $\sigma =1$ but, in this case,  this divergence just signals the pole at $s=1$ of the corresponding $L$-functions and unfortunately gives no information on their zeros. It is for this reason that the GRH for $L$-functions relative to principal characters needs an approach different from the one presented here for all the other $L$-functions relative to non-principal characters, as for instance truncating the Euler product representation of these functions in a well-prescribed manner \cite{Gonek, ALZeta}. 
 We would like to recall that, according to eq.\,(\ref{identitylprincipal}), the GRH for the $L$-functions relative to principal characters is equivalent to the original Riemann Hypothesis for the Riemann $\zeta$-function. 
 
It is worth mentioning  that the behavior of trigonometric series as the one in eq.\,(\ref{noshiftfunction}) is the topic of a famous Kac's central limit theorem  \cite{Kac}. We present such a theorem in Appendix \ref{Kactheor} where we also discuss why this result by Kac does not help in establishing the valididy of the GRH (on the Kac's theorem, see also \cite{Kiessling}).

\section{An ensemble of random primes and its associated $L$-functions}\label{randomprimeensemble}

\begin{figure}[t]
\centering
\includegraphics[width=0.45\textwidth]{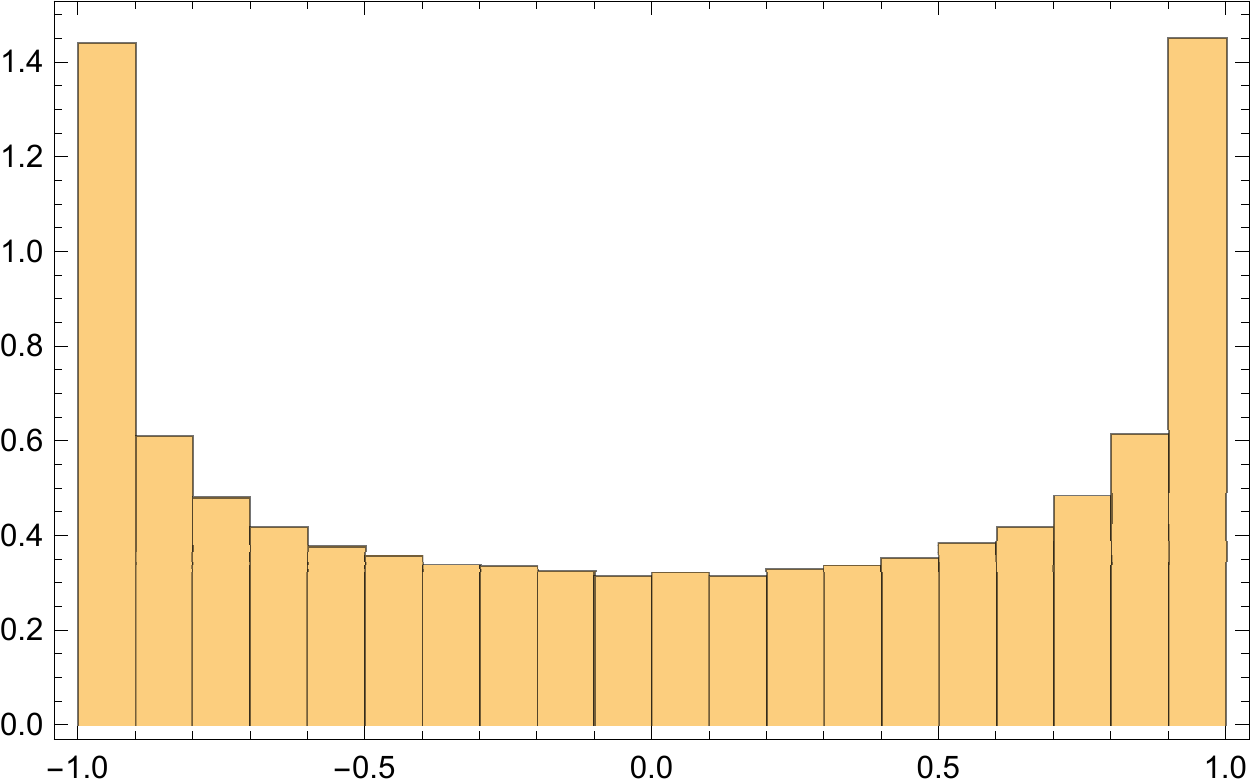}
\caption{Histogram relative to the first $10^5$ values of $b_n'(t=100)$ for $\chi_2$ mod q=7.}
\label{symmcurve}
\end{figure}

\begin{figure}[b]
\centering\includegraphics[width=.5\textwidth]{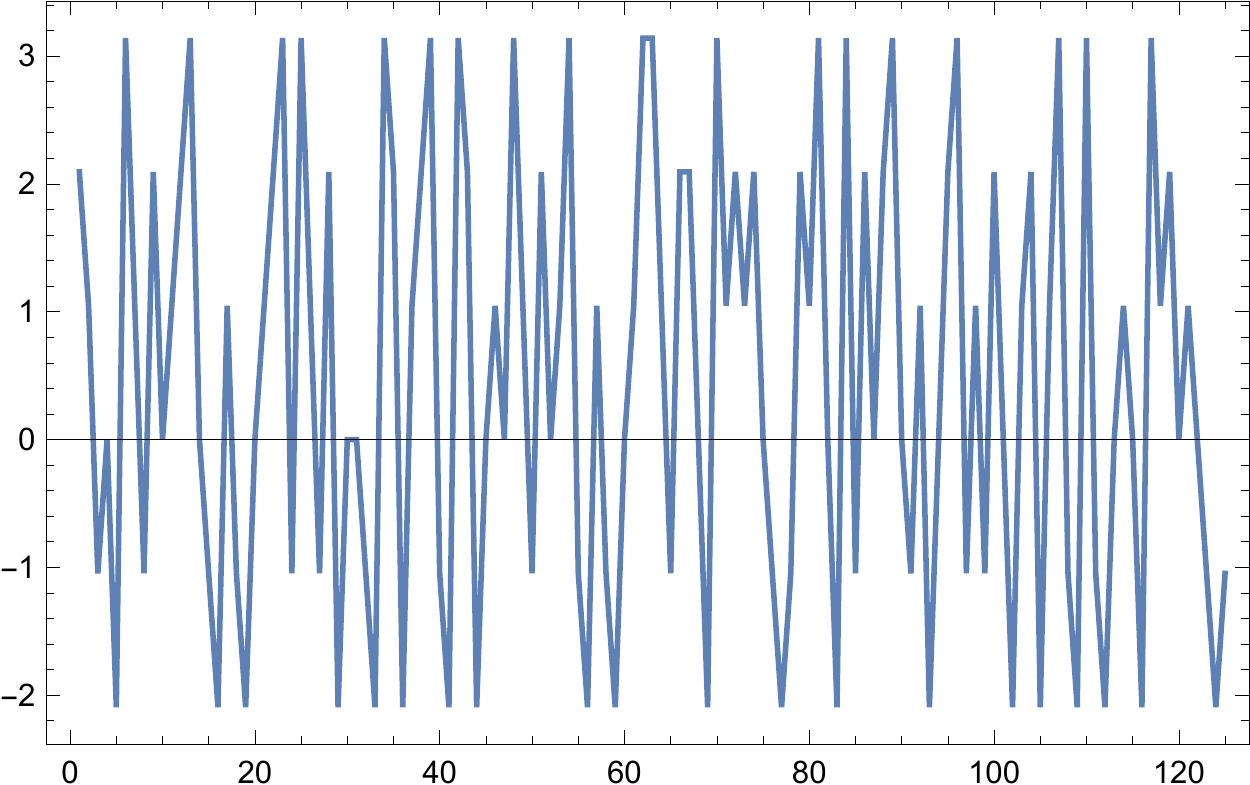}
\caption{Sequence of the angles $\theta_{p_n}$ (in unit of $\pi/3$) for the first $80$ primes relative to 
$\chi_2$ mod 7 given in Table 1.}
\label{jumphase}
\end{figure}

As we saw in the previous section, the validity of the GRH could be established if the series $\CBN (t,\chi)$ for large values of $N$ and any value of $t$ presents the purely diffusive behavior $\CBN = O(\sqrt{N})$. Such a square-root power law behavior is typically encountered in the study of the displacements 
\beq 
X_N \,=\, \sum_n^N x_n\,\,\,
\label{sumrandomvars}
\eeq 
relative to random walks, namely for sums involving independent and uncorrelated random variables $x_n$ 
of a finite variance. Following this analogy, the role of the random variables $x_n$ in our case is played by the quantities 
\beq
b_n(t) \,\equiv\,  \cos (t \log p_n  - \theta_{p_n}) \,\,\,,
\label{littlebn}
\eeq
which however are not random but completely deterministic. Yet, making an histogram of the $b_n(t)$ for a generic non zero value of $t$, one typically finds a curve as in Figure \ref{symmcurve}, which is quite close to the familiar probability distribution 
\beq
P(x) \,=\,\frac{1}{2\pi \sqrt{1-x^2}}\,\,\, 
\label{probcos}
\eeq
for a random variable $x = \cos \psi$, when the angle $\psi$ has a uniform distribution. Moreover, plotting the sequence of the values assumed by the angles $\theta_{p_n}$ varying the index $n$ (see Figure \ref{jumphase}), the apparent erratic motion of these quantities becomes immediately clear. Could this be the key for an effective random walk behavior for the deterministic series $\CBN (t,\chi)$? This would not be of course the first example of such a phenomenon: as mentioned above, we will review in Appendix \ref{Kactheor} the famous example of Mark Kac of deterministic series ruled by probabilistic normal law behavior \cite{Kac}.

These considerations suggest  exploring  a probabilistic treatment of the $b_n(t)$. Incidentally, there is a natural way to promote these quantities to be full fledged stochastic variables: this  is provided by the Grosswald and Schnitzer theorem previously mentioned. In fact, according to this theorem, we can replace the primes $p_n$ with a set of random integers $p_n'$ in the $b_n(t)$ without altering the position of the zeros of the $L$-functions. This allows us to define a set of stochastic variables $b_n'(t)$ and correspondingly an analogous series $\CBN' (t,\chi)$ for the infinite product representation of the random functions $L'(s,\chi)$ given in eq.\,(\ref{Lp}). As discussed in the following, we will see that  there holds a central limit theorem for the quantity $\CBN' (t,\chi)$! Encouraging as this may seem, it is however important to stress that this result is inconclusive towards establishing the validity of the GRH although it turns out to be useful anyway since it points to a way to sharpen our analysis, in particular to nail down the key properties which ultimately may give rise to the $O(\sqrt{N})$ growth of the original $\CBN (t,\chi)$ series. Let's now see in more detail all these steps. 

\vspace{3mm}
\noindent
{\bf A  probabilistic model  of  the primes}. Let us first define our probabilistic model which will be used to define  random $L$-functions $L'(s,\chi)$.   Let $\Primes = \{p_1, p_2, \ldots \}$  denote the set of primes,  where $p_1 =2,  p_2 = 3, $ and so forth. We will consider replacing $\Primes$ with the set $\Primesp = \{ p'_1, p'_2, \ldots \}$ where $p'_n$ is a randomly chosen integer satisfying 
 \beq
 \label{ppnn}
 p_n \leq p'_n < p_n+K, ~~~~~~p'_n = p_n ~~ ({\rm mod} ~ q)
 \eeq
where $q$ is the modulus of the Dirichlet character and $K\geq q$ is an arbitrary integer.   
To simplify the analysis, we can take $K = M \,q$ for some positive integer $M$,  such that 
\beq
\label{pprange}
 p'_n \in  \{p_n, p_n + q, p_n+2 q, \ldots, p_n + M q \} \,\,\,.
\eeq
Hence our stochastic model can be simply viewed  as follows: it consists in randomly chosing an integer $m_n \in [0,M]$   with equal probability,  and, at any $n$-th step of this process, we assign as output the integer  
\beq 
p'_n \,=\, p_n + m_n q
\,\,\,.
\label{fake}
\eeq
Therefore we are dealing with a sequence of independent and random integers $0 \leq m_n \leq M$ which are  superimposed onto a  ``ramp"  given by the primes $p_n$. Notice that $M$ can be {\em any} integer, in particular it can be arbitrarily large, so that the values assumed by the random variables $m_n$ can be spread out on  an arbitrarily  large interval of integers. Moreover, the ramp dictated by the primes $p_n$ does not effect either the independence of the $m_n$ nor their equal probability;  it simply implies that the random numbers $p'_n$ grow as the primes $p_n$ when $n \rightarrow \infty$. Since the $p'_n$ are random variables, we are led to  consider $\Pensemble =  \{ \Primesp \}$   which is the ensemble of all possible  $\Primesp$,  i.e the set of sets $\Primesp$.  We will refer to $\Pensemble$  as the {\em random-prime  ensemble},  and a specific element $\Primesp \in \Pensemble$ as a {\it state} of this  ensemble.  The actual primes $\Primes$ are then simply one state in this random-prime ensemble, more precisely the state in which $m_n = 0$ for {\em all} $n$. 

Given a state $\Primesp$, we can now define a modified $L$ function 
\beq
\label{Lppp}
L' (s, \chi) = 
\prod_{n=1}^\infty  \( 1 -  \frac{\chi(p'_n)}{(p'_n)^s} \)^{-1} \,\,\,,
\eeq
which is now a random function; yet, according to the Grosswald-Schnitzer theorem, it has exactly the same zeros as  $L(s,\chi)$ inside the critical strip. This suggests that a possible approach to proving the GRH consists of studying the convergence properties of the infinite product (\ref{Lppp}): if we were able to show that at least a  {\it single}  state $\Primesp$ leads to a $L'$ function with no zeros to the right of the critical line, then this implies the  validity of the GRH, since for a given $\chi$,  all the $L'(s, \chi)$ have the same zeros. For this reason, let's then focus our attention on the series 
\beq
\label{BprimesGS}
B'_N (t)  \,=\, \sum_{n=1}^N  b_n'(t)\,\,\,\,\,\,\,\,\ , ~~~~~~~~~b_n' (t) = \cos (t \log p'_n - \theta_{p_n})
 \,\,\,. 
\eeq
If $B'_N$ obeys  an appropriate central limit theorem, then an arbitrarily large fraction of the $B'_N$ are $O(N^{1/2 + \epsilon})$ for arbitrarily small positive $\epsilon$ and Theorem \ref{BNtheorem} would then imply there are no zeros to the right of the critical line for at least one state. In other words,  Theorem \ref{GrosswaldSchnitzer} would then promote {\it almost surely true} statements,  i.e. statements that are true with probability $1$, to {\it surely true}.   For clarity of presentation,  let us state this as a theorem:

\begin{theorem}
\label{BNpthm}
Given a character $\chi$, define the series on the random numbers $p_n'$
\beq
\label{BprimesGGS}
B'_N (t)  = \sum_{\substack{n=1 \\ p_n \nmid k}}^{N}     b'_n(t) \,\,\,\,\,\,\,\, , ~~~~~~~~~b'_n (t) = \cos (t \log p'_n - \theta_{p'_n})
\eeq
where $\theta_{p'_n} = \theta_{p_n}$.  
If  $B'_N (t) = O(N^{1/2+\epsilon})$ with probability equal to $1$  for any $\epsilon > 0$,  then the RH is true for the $L$-function based on this character.
\end{theorem}

\vspace{1mm}

\noindent
{\bf Remark}. {\it Notice that power law of $B'_N $ cannot be less that $1/2$. Indeed, if $B'_N = O(N^{1/2+ \epsilon})$ for $\epsilon < 0$,  then this would rule out zeros on the critical line,  which instead we know exist
\cite{Conrey2}.}

\bigskip
\noindent
 As we are going to show below,   $B'_N$ does obey a central limit theorem but, as we will explain, this result is not decisive to establish the validity of the GRH.

\begin{comment}
\blue{I think there are too many figures in the rest of this section,  several of which show the same thing.   Furthermore,  they don't convey a great deal.    
The reader has to decipher that we are not plotting jumps,  etc., so they can't figure out what is going on with a quick glance.   We should reserve the proliferation of figures for our main results later.   
I recommend dropping Figures 6, 8 and 9,   or in any case keep  only one of Figures 6, 7,8, or 9,  since they all show the same thing. I prefer Fig. 7. }

\blue{There more comments on figures below.    The most important figures are in the appendices!   Whereas we have all these tangential figures 6,7,8,9!}

\end{comment}

\vspace{3mm}
\noindent
{\bf Some properties of the random sequence $\{{\bf b_n'(t)}\}$.} Even though $B'_N(t)$ is a sum of random variables, there are however some differences with the standard sum of a random walk: for instance, the $b_n'(t)$ are not identically distributed and this may lead to a non-zero drift. It is useful to get more familiar with the properties of the sequence of $b_n'(t)$ for $n=1,2,\ldots N$. First of all, we express them as $b_n'(t) = \cos\psi_n'(t)$, where the angles $\psi_n'(t)$ are defined according to   
\beq
\psi_n'(t) \,\equiv\, t \log p'_n - \theta_{p_n}\,\,\,.
\eeq
In the following, when we say \textquotedblleft {\em short scale}\textquotedblright \, we mean looking exactly at the jumps, i.e. the strong fluctuations, which occur in the sequence of $\psi_n'(t)$ in passing from $n$ to $n+1$. On the other hand, disregarding their short scale jumps, the $b_n$ fall into a set of values whose envelopes, varying $n$, form a set of continuous curves: in the following, when we say \textquotedblleft {\em large scale}\textquotedblright  \,we mean looking at this continuous and smooth pattern of the $b_n$'s which emerges not taking into account the erratic jumps of the sequence in passing from $n$ to $n+1$ but  considering instead large intervals of the index $n$.  We emphasize that \textquotedblleft {\em long scale}\textquotedblright  \,does not signify at all a different kind of behavior:  if we zoom on any region of a plot of the long scale behavior and pay attention to the local jumps,  there  appears of course the erratic short scale behavior.  Let's consider, for simplicity, the case when the cardinality $r$ of the set of angles  in \eqref{PhiSet} of the characters coincides with $\varphi(q)$. There are essentially two situations to consider, according whether $t=0$ or $t \neq 0$. 

\begin{itemize}
\item When $t=0$,  as is evident from Figure \ref{jumphase}, the angles $\theta_{p_n}$ jump erratically (although deterministically) among all possible roots of unity relative to  the modulus $q$ of the characters: there are $\varphi(q)$  such roots of unity,   but in view of the identity $\cos(a) = \cos (-a)$, the sequence of $b_n'(0)$ consists of $1 + \varphi(q)/2$ values only. On  short scales, the sequence of $b_n'(0)$ jumps discontinuously from one value to another, as shown in left hand side of Figure \ref{plotsequence1}, while on large scales (i.e. disregarding the individual jumps) one observes a set of $1 + \varphi(q)/2$ flat values, as those shown in the plot on the right-hand side of Figure \ref{plotsequence1}.

\begin{figure}[t]
\centering
\includegraphics[width=0.45\textwidth]{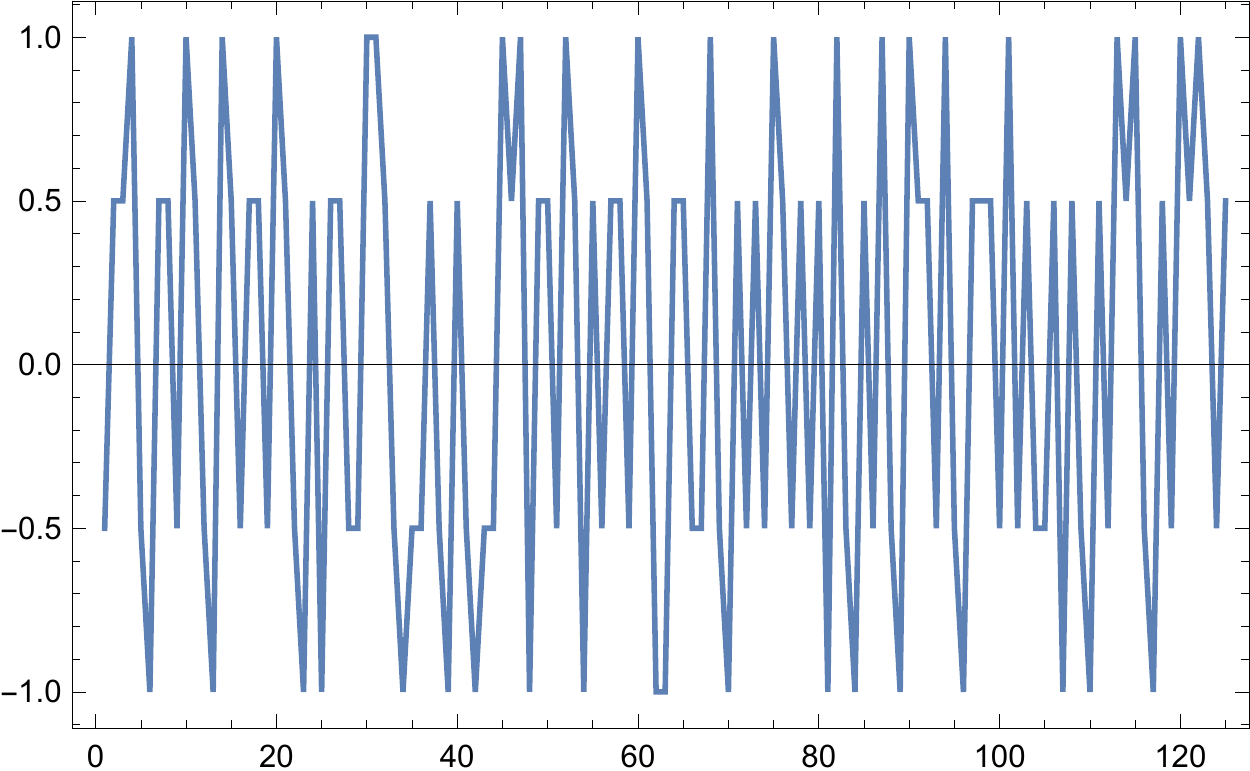}
\,\,\,
\includegraphics[width=0.45\textwidth]{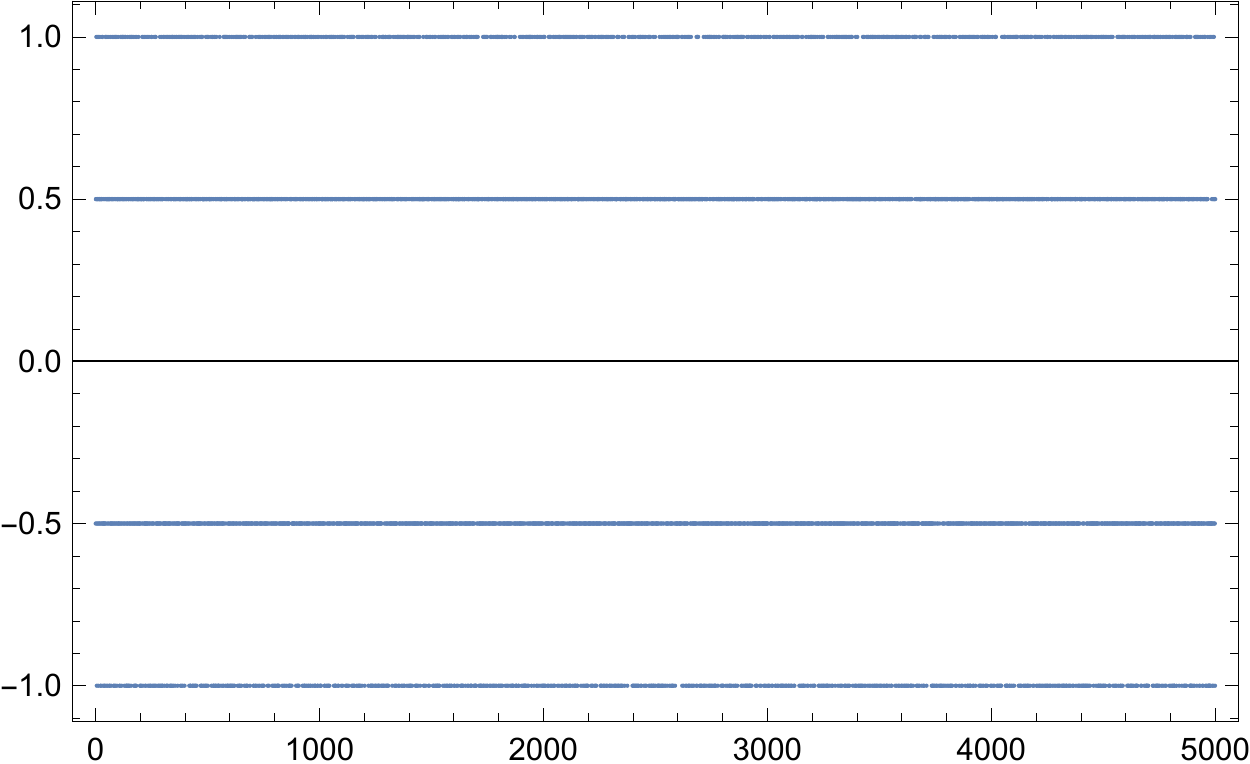}
\caption{$b_n'(0)$  for the character $\chi_2$ mod q=7. (a) Left-hand side: short scale behavior of the sequence,  first 125 values. Successive points are jointed to emphasize their jumps.
b) Right-hand side: large scale behavior of the sequence, first 5000 values. Successive points are not jointed in order to show their large scale smoothness.  
} 
\label{plotsequence1}
\end{figure}

\begin{figure}[b]
\centering
\includegraphics[width=0.45\textwidth]{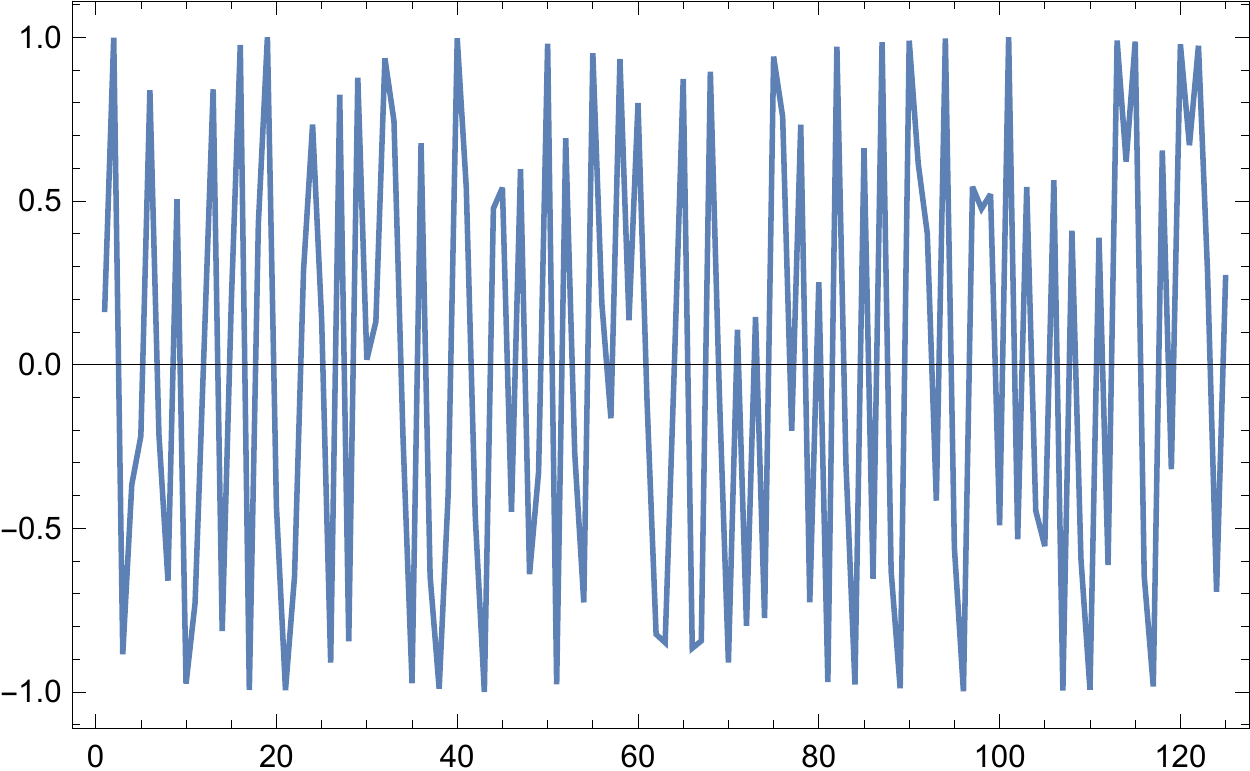}
\,\,\,
\includegraphics[width=0.45\textwidth]{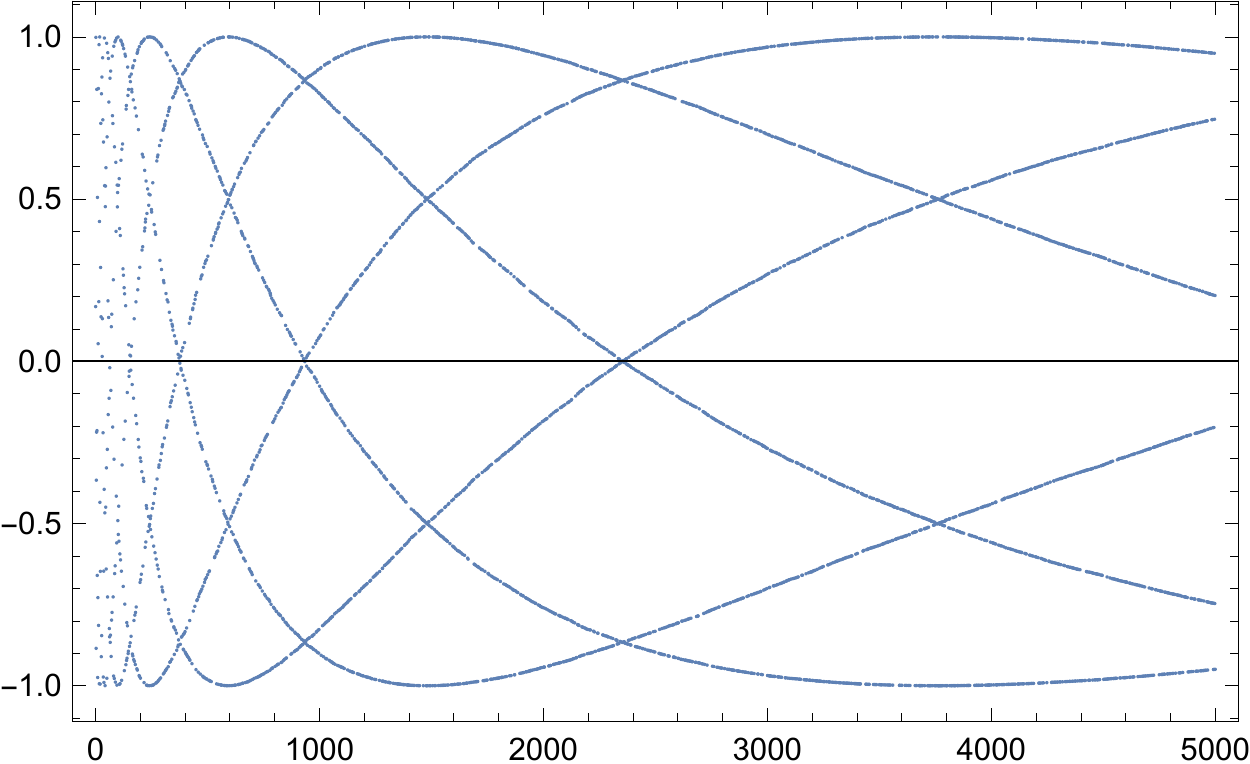}
\caption{$b_n'(1)$  for the character $\chi_2$ mod q=7. (a) Left-hand side: short scale behavior of the sequence,  first 125 values. b) Righ-hand side: large scale behavior of the sequence, first 5000 values.} 
\label{plotsequence2}
\end{figure}

\begin{figure}[t]
\centering
\includegraphics[width=0.45\textwidth]{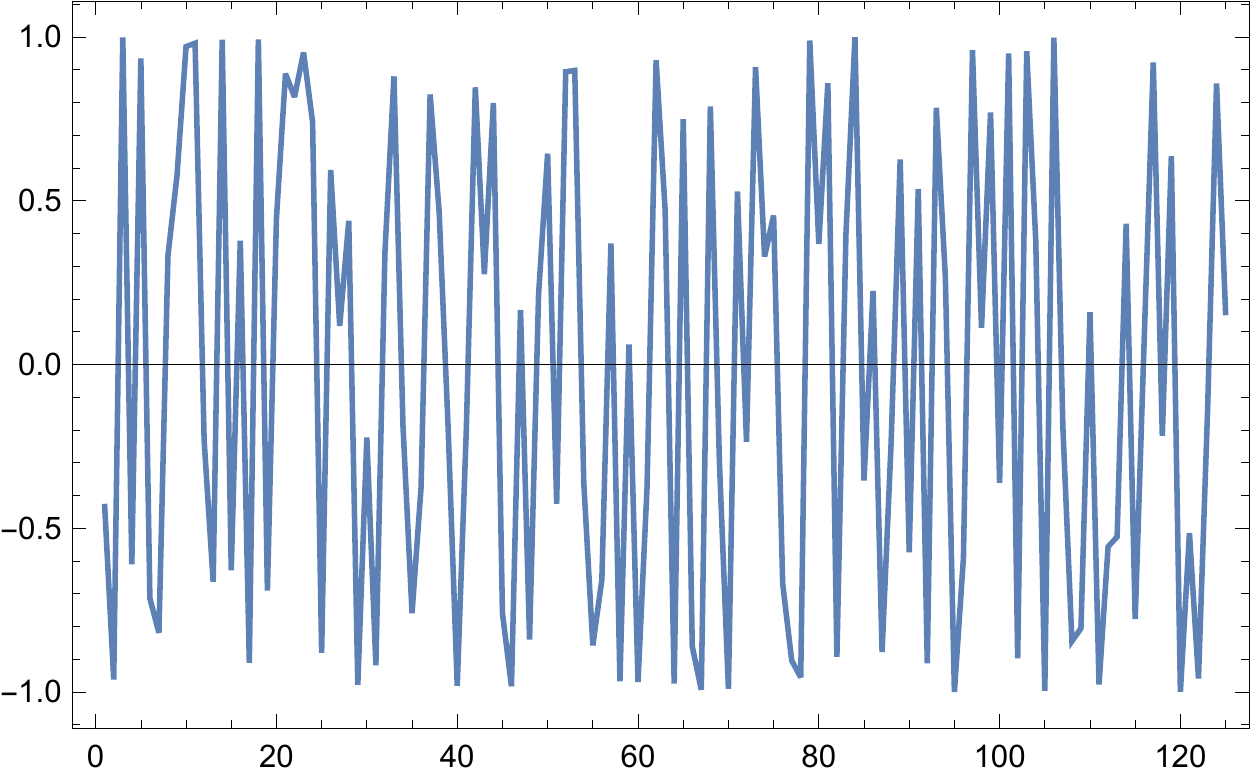}
\,\,\,
\includegraphics[width=0.45\textwidth]{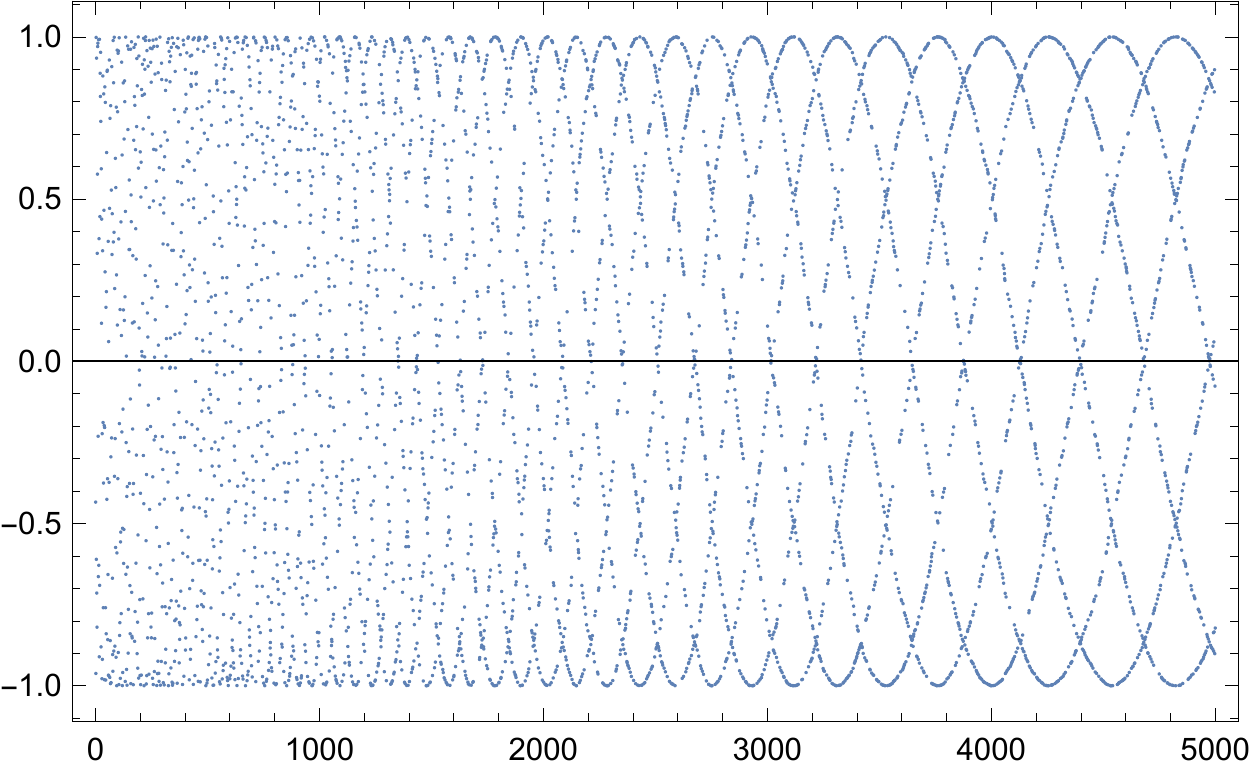}
\caption{
$b_n'(15)$  for the character $\chi_2$ mod q=7. (a) Left-hand side: short scale behavior of the sequence,  first 125 values. b) Righ-hand side: large scale behavior of the sequence, first 5000 values.} 
\label{chaos1}
\end{figure}

\begin{figure}[b]
\centering
\includegraphics[width=0.45\textwidth]{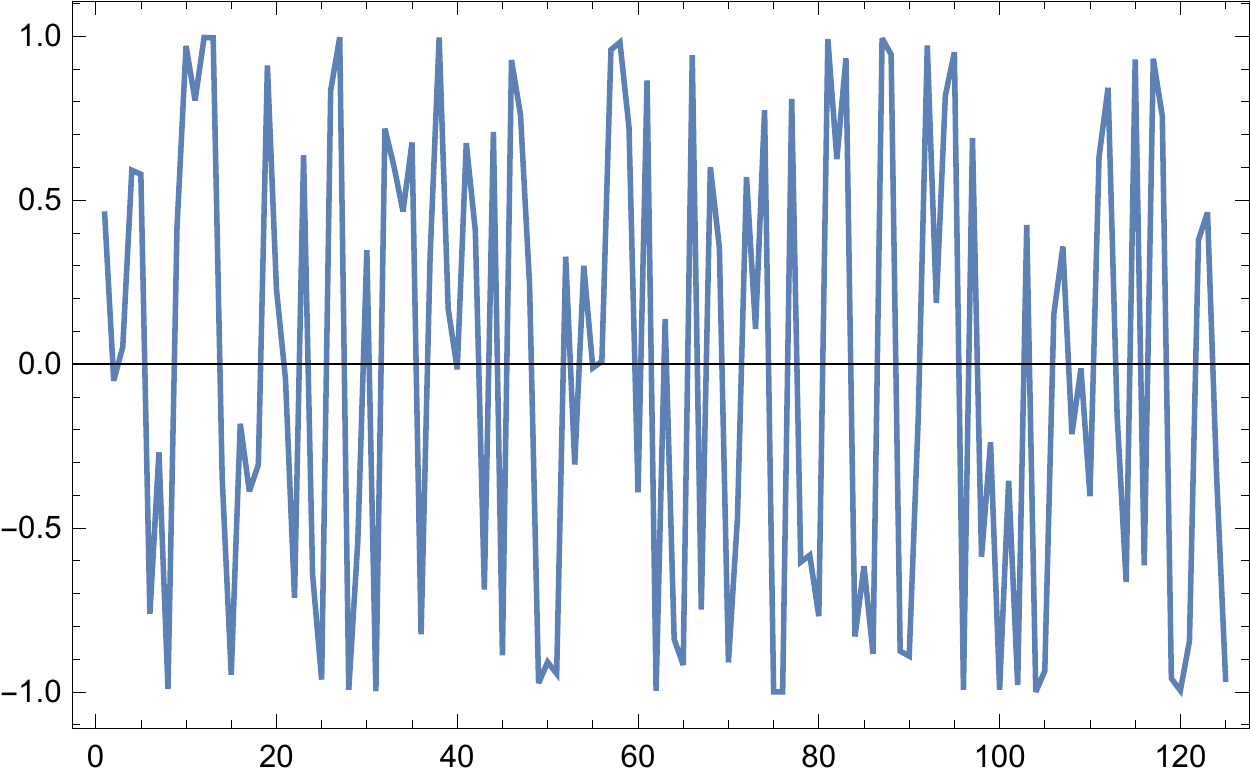}
\,\,\,
\includegraphics[width=0.45\textwidth]{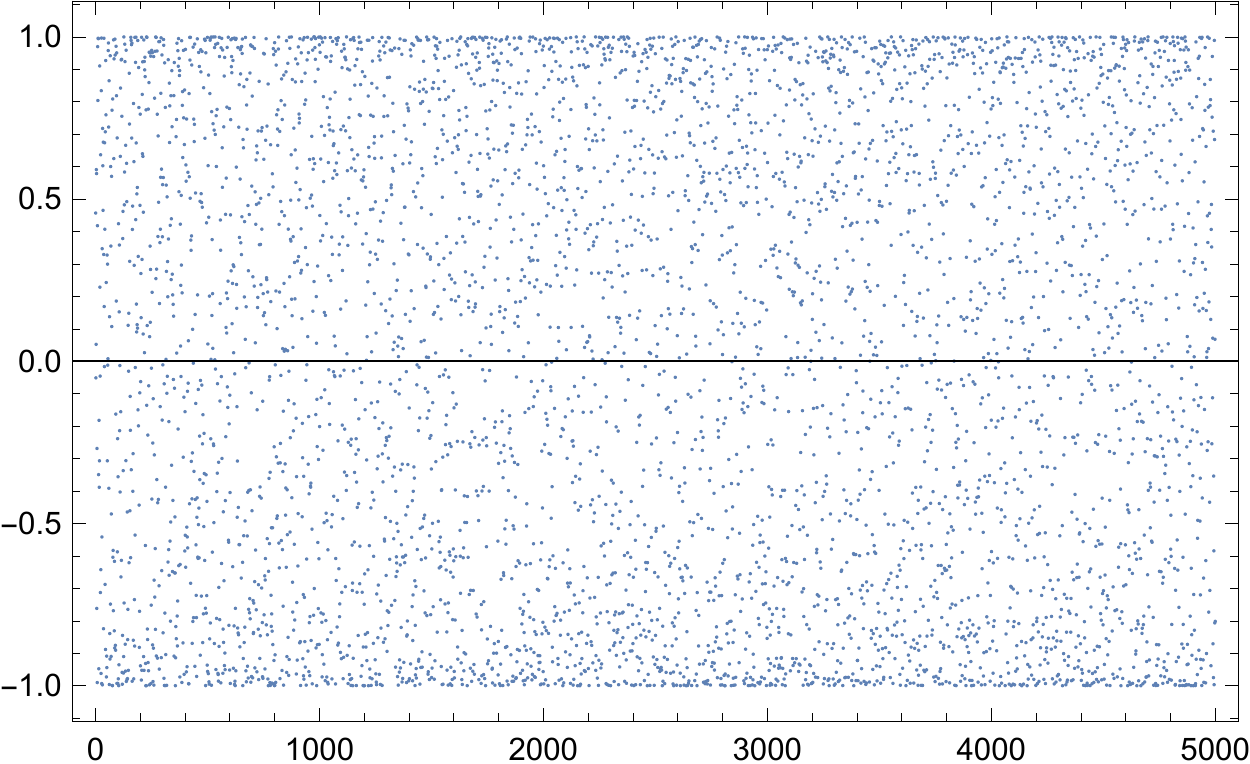}
\caption{$b_n'(500)$  for the character $\chi_2$ mod q=7. (a) Left-hand side: short scale behavior of the sequence,  first 125 values. b) Righ-hand side: large scale behavior of the sequence, first 5000 values.} 
\label{chaos}
\end{figure}

\item
When $t \neq 0$, no matter how small, the degeneracy of some of the previous straight lines is lifted and there are two mechanisms (although of quite different nature) which help in scrambling the angles $\psi_n'(t)$ and in giving rise to the apparent random behavior of the $b_n'(t)$'s: 
\begin{enumerate}
\item The first mechanism is due to the purely random term $t \,\log p'_n$. For the logarithm present in this expression, at a given $t$ this term changes slowly going from $n$ to $n+1$ and therefore it is necessary to arrive at an index $k$ such that $p_{n+k}' \simeq p_n' e^{\pi/t}$ to induce a change of phase equal to $\pi$ in the difference ($\psi_{n+k}(t) -\psi_{n}(t)$). Of course, the  larger the value of $t$,  the faster the change of the phase induced by this genuine random term. Imagine, in fact, that $t$ is large: to get a phase change equal to $\pi$ going from two consecutive indices $n$ and $n+1$, using the approximate formula $p_n \approx n \log n$, one determines that it is necessary to arrive to the index $n \approx t/\pi$. Although this mechanism may be considered a slow scrambling of the phase (especially for small values of $t$), it is nevertheless a mechanism present for any non-zero $t$. 
\item The second mechanism, which is definitely more effective and faster in scrambling the phase $\psi_n(t)$, is due to the previously discussed nearly chaotic jumps of the phase $\theta_{p_n}$ computed on the sequence of the primes. 
\end{enumerate}
As a result of these two scrambling mechanisms which, it is worth to underline, work for both the random sequence $b_n'(t)$ and 
the deterministic sequence $b_n(t)$, on large scales one observes  $\varphi(q)$ separate curves, statistically distributed in a symmetrical way with respect the vertical axis, as those shown for instance in Figures \ref{plotsequence2} and \ref{chaos1}, while on short scales, an erratic series of jumps among all possible values of these curves.   
\end{itemize}
It is also useful to notice that, increasing $t$, and in particular taking the limit $t \rightarrow \infty$, the scrambling of the values becomes more and more effective and there is indeed a smooth transition from a situation in which there is a set of $1 + \varphi(q)/2$ distinct curves to a situation in which there is a chaotic filling of the rectangle of sides $N \times 1$ in terms of the points of the sequence $\{b'_n(t)\}$, as shown in Figure \ref{chaos}. This is the reason which is behind a curve as the one in Figure \ref{symmcurve} for the histogram of the $b'_n(t)$'s.

\section{ A Central Limit Theorem for $B'_N (t)$}\label{centrallimittheorem}

In this section we prove a central limit theorem for the series $B'_N$ which has been the focus of our discussion thus far.  
Let us  first  recall Lyapunov's theorem which states the sufficient conditions under which the normal law for a set of random variables applies,  even if they are not equally distributed.

\vspace{3mm}
\begin{theorem}\label{Lyapunovtheorem1}(Lyapunov) 
{\it Let $x_n$, $n=1,2,\ldots, N$  be independent random variables with finite  mean $\mu_n$  and variance $\sigma_n^2$,  which are allowed to vary with $n$, and define the series $X_N = \sum_{n=1}^N  x_n$. Define  $m_N$ as the expectation value of $X_N$,
\beq
\label{mN}
m_N  =  \Ex \Bigl[X_N \Bigr]  =  \sum_{n=1}^{N}  \mu_n\,\,\,, 
\eeq 
and  $s_N^2$ the sum of variances
\beq
\label{sN2} 
s_N^2  =  \sum_{n=1}^{N}  \sigma_n^2 \,\,\,.
\eeq
If the Lyapunov condition is satisfied,  namely if for some $\delta >0$  
\beq
\label{LyCond}  
\lim_{N \to \infty}    \inv{s_N^{2 + \delta} }  \sum_{n= 1}^{N}   \Ex \Bigl[ | x_n - \mu_n |^{2 + \delta} \Bigr]   = 0,  
 \eeq 
then 
\beq
\label{Lyapunov}  
\inv{s_N}     \Bigl(  X_N  - m_N \Bigr)   ~ \dist ~  \CN(0,1)
\eeq
where $\CN(\mu, \sigma )$ is the normal distribution with mean $\mu$ and variance $\sigma$
\beq
\CN(\mu,\sigma) \,=\,
\frac{1}{\sqrt{2\pi\sigma^2}} \, e^{-\frac{(x-\mu)^2}{2\sigma^2}} \,\,\,.
\label{normaldisss}
\eeq
}  
\end{theorem}

\vspace{3mm}
\noindent
We are now in the position to establish the central limit theorem for the quantity $B_N'(t)$, as far as $t\neq 0$. 
With the values of $p'_n$ given in eq.\,(\ref{fake}), the central limit theorem involves both the quantities
\begin{eqnarray}
\label{munNP}
\mu_n & \, =\, & \Ex [b_n] = \inv{M+1} \sum_{m=0}^M \cos \[ t \log (p_n +m q) - \theta_{p_n} \],
\\ 
\label{sigmaNP}
\sigma_n^2 &\,=\, &\Ex [(b_n - \mu_n)^2] \,=\, \inv{M+1} 
\sum_{m=0}^M  \( \cos \[ t \log (p_n+ m q) - \theta_{p_n} \] - \mu_n \)^2 \\
&\,=\,&   
\inv{M+1} 
\sum_{m=0}^M \cos^2 \[ t \log (p_n+ m q) - \theta_{p_n} \] - \mu_n^2 \,\,\,. \nonumber
\end{eqnarray}
and their sums (\ref{mN}) and (\ref{sN2}). 
According to Lyapunov's theorem we then have   

\vspace{1mm}

\begin{theorem}
\label{BpCLTDir}  
For any $t\neq 0$, 
\beq
\label{CLTDirEq}
\inv{s_N}     \Bigl(  B'_N  - m_N \Bigr)   ~ \dist ~  \CN(0,1) \,\,\, ,
\eeq
where $m_N$ and $s_N$ are defined in \eqref{mN} and \eqref{sN2} along with \eqref{munNP} and \eqref{sigmaNP}.  
\end{theorem}

\bigskip

\begin{figure}[t]
\centering\includegraphics[width=.6\textwidth]{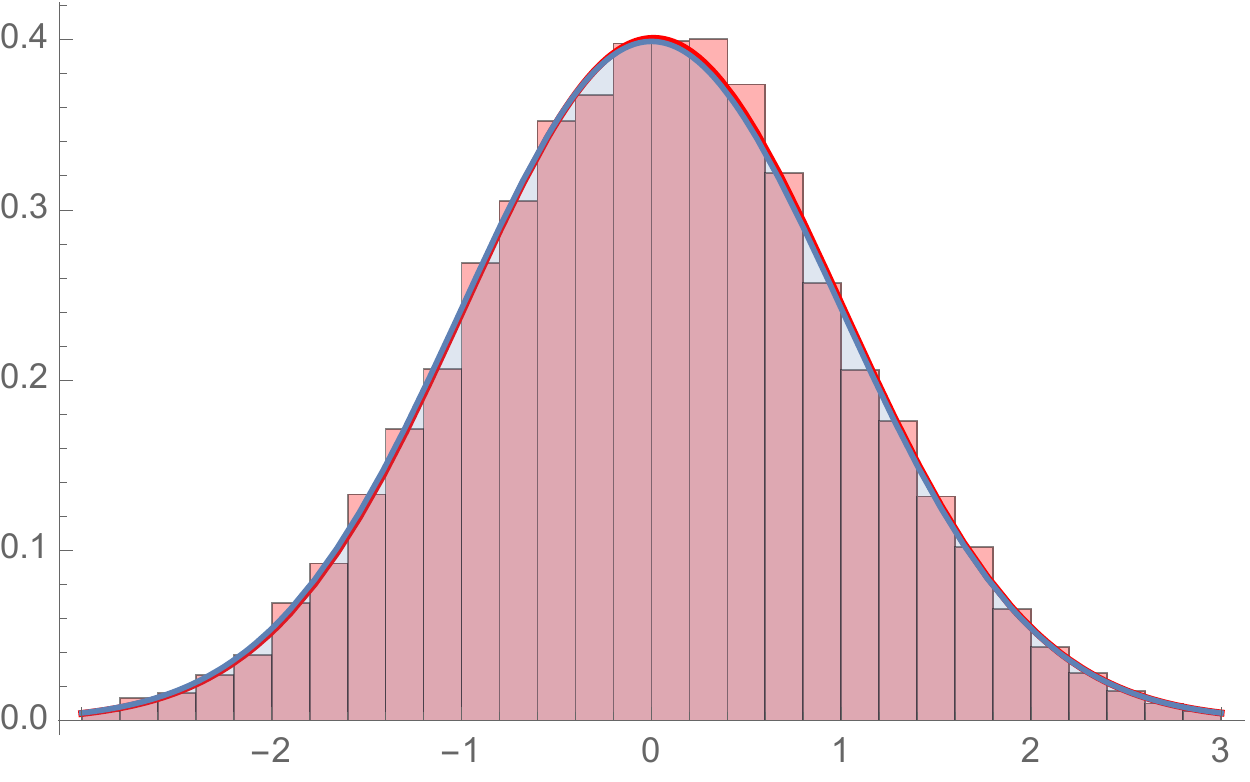}
\caption{Numerical evidence for Theorem \ref{BpCLTDir} based on the character \eqref{char}.      We fixed N=10,000, $M=3$,  and $t=100.$  
Displayed is a normalized histogram of the LHS of \eqref{CLTDirEq}  for $20,000$ states $\Primesp$.     The red curve is the  fit to the data,  which is the normal 
distribution $\CN(0.0108, 0.9955)$.   The  nearly indistinguishable  blue curve is the prediction $\CN(0,1)$.}
\label{CLT_Dir}
\end{figure}

Clearly,  the theorem only applies to $t\neq 0$ since for $t=0$,  $B'_N = m_N$.  Numerical evidence for Theorem \ref{BpCLTDir} can be found in Figure \ref{CLT_Dir}.

\section{The role of the mean}\label{mean}

For the purpose of  establishing  that $B'_N = O(\sqrt{N})$, the existence of a normal law distribution as the one discussed in the previous Section is unfortunately inconclusive since the Theorem \ref{BpCLTDir} concerns  the {\em difference} $(B'_N - m_N)$ rather than $B'_N$ itself! Yet, we can learn something from this analysis. 
First of all, for $t\neq 0$ we have $B'_N \neq m_M$.  Now $m_N$ is the average of $B'_N$,  and since the averaging smoothes out large fluctuations of $B'_N$, the growth of $m_N$  is either slower or the same as that of $B'_N$. Thus, unless there is some delicate cancellation in the difference $(B'_N - m_N)$ which occurs {\it for any t}, Theorem \ref{BpCLTDir} would imply that at worse both $B'_N$ and $m_N$ are  $O(\sqrt{N})$,  so that their difference is also $O(\sqrt{N})$. Unfortunately, we cannot rule out miraculous cancellations in the difference $(B'_N - m_N)$, thus we are back to studying the deterministic series $m_N$, which is quite similar to the original series $B_N$. In fact, for $M=0$, $B_N = m_N$.    
    
In hindsight, for $L$-functions of non-principal cases one can argue (discussed in Appendix \ref{growthBN})  that the asymptotic behavior in $N$ of the series $B_N(t)$ is entirely dictated by their behavior at $t=0$.    Of course,  a simpler argument just relies on the well known fact that the domain of convergence of $L$-functions are always half planes \cite{Apostol}.   
In light of this result, in the remainder of this paper we will only consider the series $C_N \equiv B_N (t=0)$.  
The focus is now on the following theorem:

\vspace{3mm}
\begin{theorem}\label{seriesinzero}
(Fran\c ca-LeClair \cite{EPFchi})
{\it Consider the sum on the primes 
\beq
\label{Gseries}
C_N \,=\, \sum_{\substack{n=1 \\ p_n \nmid q}}^{N}    \cos  \theta_{p_n}   
\,\,\,,
\eeq
and assume that for large $N$ it scales as 
\beq 
C_N  \simeq N^\alpha \,\,\,, 
\eeq
up to logarithms (see the discussion below Theorem \ref{BNtheorem}).
Then the GRH is true if 
\beq
\alpha = 1/2 + \epsilon \,\,\,
\eeq
 }
 for all  $\epsilon > 0$.   
\end{theorem}
Indeed, if $\alpha = 1/2$ then the function $B_N(t)$ grows  as $\sqrt{N}$ (up to logarithms) for {\em all} values of $t$ and therefore, using Theorem \ref{BNtheorem}, the convergence of the infinite product of the $L$-function can be safely extended down to  the critical line  $\Re (s) = \half$ without encountering any zeros.  
 
\def\s{\ell}
 
\section{The series $C_N$: Insights from random time series}\label{timeseriess}
\def\hatA{A} 

In light of Theorem \ref{seriesinzero}, the crucial quantity of our analysis has now become the series $C_N$  
\beq 
C_N \,=\, \sum_{n=1}^{N}    \cos  \theta_{p_n}   \equiv \sum_{n=1 }^{N}  c_n
\,\,\,, 
\label{defCNN}
\eeq
made of the sequence of the angles $\theta_{p_n}$ relative to the first $N$ primes
\beq
\hatA_{N}  \,=\, 
\{\theta_{p_n}\, ; ~  n = 1, 2, \ldots, N \} \,\,\,. 
\label{sequenceSNpart}
\eeq
For later use, let's also define the  {\it ordered} intervals of length $N$ starting at $\s$  
\beq
I_N(\s) =\{\s, \s+1, \s+2, \ldots, \s+N -1 \} \,\,\,, 
\,\,\,
\label{intervals}
\eeq 
and the associated sequence of angles $\hatA_{N}(\s) $ 
\beq
\hatA_{N} (\s)  \,=\, 
\{\theta_{p_n}\, ; ~ n \in I_N(\s)\}\,\,\,. 
\label{sequenceSNpart}
\eeq
Let's remind that the values of the angles $\theta_{p_n}$ belong to a finite and discrete set (see eqs. (\ref{thetan2} - \ref{PhiSet}))
\beq
\theta_{p_n} \in \Phi \,=\, \{ \phi_1,\phi_2,\ldots,\phi_r\}
\,\,\,\,\,\,\,\,
,
\,\,\,\,\,\,\,\,
r \leq \varphi(q) \,\,\,  .
\label{setangles}
\eeq
Hence the series $C_N$ is made of terms all of the same order, always smaller or equal to 1. Moreover, one could expect that the angles $\theta_{p_n}$ computed on the sequence of the primes are {\em equally distributed} among their possible $r$ values and, as a consequence, the values of the cosine of these angles are pairwise equal and opposite. If the $c_n$ were uncorrelated random variables with the properties just described, i.e. variables of average $\mu =0$ and finite variance $\sigma$, then the $\sqrt{N}$ behavior of the series $C_N$ will be simply guaranteed by the Lyapunov theorem recalled in Section \ref{randomprimeensemble}. 

To make any progress on the behavior of the series $C_N$ it is then necessary to study in more detail the statistical properties of the angles $\theta_{p_n}$ and their relative cosine. In the next section we will see that several properties of $\hatA_N$ are captured by the Dirichlet theorem \cite{Diric,SelbergD} and the Lemke Oliver-Soundararajan conjecture on the distribution of pairs of residues on consecutive primes \cite{OliverSoundararajan}. These two mathematical results will constitute the final and definite {\em theoretical} statements on the sequence of $\hatA_{N}$ on which we will base our future analysis. However, in this section we want to explore a different route, i.e. here we are going to study the sequence $\hatA_N$ from an {\em experimental} point of view. This means  that we are going to consider the angles $\theta_{p_n}$ as if they were the outputs of a random time sequence (of which we pretend to ignore the origin), with the role of discrete time played by the index $n$. From this point of view, $C_N$ assumes the meaning of a random time series and we can take advantage of several numerical methods developed to study these quantities \cite{timeseries1,timeseries2} to get some conclusions of pure statistical nature on our series $C_N$. Let's see what we can learn following these lines of thought, analyzing some significant examples. 

Let's choose for instance $q=7$: the maximal set of angles associated to the non-principal characters is shown 
in Table I and consists of 
\beq 
\Phi =\{
\alpha_1 = -2\pi/3 \,\,\,,
\,\,\, \alpha_2 =-\pi/3 \,\,\,,
\,\,\, \alpha_3 = 0 \,\,\,,
\,\,\, \alpha_4 = \pi/3 \,\,\,,
\,\,\, \alpha_5 = 2\pi/3, \,\,\,,
\,\,\, \alpha_6 = \pi \} 
\,\,\,.
\label{spettroangoli}
\eeq
Let's now take the character $\chi_2$ relative to this modulus and write the sequence of the corresponding angles relative to the increasing sequence of primes 
\beq
\hatA_{N} = \{ \alpha_5, \alpha_4, \alpha_2, \alpha_3, \alpha_1, \alpha_6, \alpha_4, \alpha_2, \alpha_5, \alpha_3, \alpha_4, \alpha_5, \alpha_6, \alpha_3, \alpha_2, \ldots \} 
\,\,\,.
\label{exampleseq}
\eeq
As we said, let's pretend that for this example and all the others we do not know where these sequences come from and therefore let's treat them as they were some random outputs associated to the rolling of a dice of $r$ faces. As for the rolling of a dice, it is then perfectly legitimate to enquire about the distribution of the various $\alpha$'s, the frequency of each of them, and whether the dice is biased or not, namely if the various outputs are correlated and how much they are correlated. 

\vspace{3mm}
\noindent
{\bf Relative probabilities}. Varying $N$, we can first study how many times the angles $\alpha_k$ appear in the sequence $\hatA_N$ and therefore define their relative probability $\Prob2_k$ as  
\beq
\Prob2_k \,=\, \frac{\# \,\{ \alpha_k \,\, \in \,\hatA_N\}}{N} \,\,\,. 
\label{frequencies}
\eeq
For instance, in the first 15 terms of the sequence (\ref{exampleseq}), the angle $\alpha_1$ appears only 1 time, $\alpha_2$ appears 3 times, $\alpha_3$ appears 2 times, $\alpha_4$ appears 3 times, $\alpha_5$ appears 3 times while $\alpha_6$ appears 2 times. Few examples will help to identify the trend of these relative probabilities.

\begin{table}[t]\label{tablefreq72}
\centering
\begin{tabular}{ |c | c | c | c | c | c | c | c | c | }
\hline\hline
N & $5^3$ & $5^4$ & $5 ^5$ & $5^6$ & $5^7$ & $5^8$ & $5^9$ & 
$5 ^{10}$ \\  
\hline
$\Prob2_1$ &
0.16000 & 0.16480 & 0.16544 & 0.16589 & 0.16644 & 0.16656 & 0.16664 & 0.16664 \\
\hline
$\Prob2_2$ & 
0.17600 & 0.16640 & 0.17024 & 0.16717 & 0.16696 & 0.16667 & 0.16670 &
0.16668 \\
\hline
$\Prob2_3$ & 0.16800 & 0.16160 & 0.16416 & 0.16659 & 0.16652 & 0.16662 & 0.16664 & 0.16668 \\ 
\hline
$\Prob2_4$ &
0.18400 & 0.17120 & 0.16640 & 0.16691 & 0.16698 & 0.16678 & 0.16668 & 0.16669 
 \\
 \hline
$\Prob2_5$ & 0.16000 & 0.16320 & 0.16512 & 0.16672 & 0.16646 & 0.16660 & 
0.16662 & 0.16666
 \\
 \hline
$\Prob2_6$ & 0.15200 & 0.17280 & 0.16864 & 0.16672 & 0.16664 & 0.16676 & 0.16673 & 0.16665
 \\
\hline
\end{tabular}
\caption{The probabilities $\Prob2_a$ of angles $\alpha_a$ ($a=1,\ldots,6$) relative to the 
character $\chi_2$ module 7 vs the length $N$ of the sequence $S_N$.}
\end{table}

\vspace{3mm}
\noindent
{\em First example}.
As a first example let's consider the statistics of the angles for the character $\chi_2$. Taking for $N$ the first $N=5^{10} = 9,765,625$ primes, it is evident that the probabilities of these angles tend to a common value equal to $1/6$, and the manner in which they reach these asymptotic values\footnote{The tiny differences between these probabilities at a finite $N$ (which become smaller and smaller with increasing $N$) can be traced to a well-known phenomenon, the so-called {\em Prime Number Races} (for a nice review on this subject, see \cite{GranvilleMartin}).} can be surmised by examining Table II. For $N=5^{10}$, the various relative probabilities differ each other for 
about $0.02\%$.

\vspace{3mm}
\noindent
{\em Second example}. As a second example, we consider the character $\chi_3$ (mod $7$). In this case there are only 3 angles which, adopting the same notation as before, are given by 
\beq 
\Phi =\{
\alpha_1 = -2\pi/3 \,\,\,,
\,\,\, \alpha_3 = 0 \,\,\,,
\,\,\, \alpha_5 = 2\pi/3 \} 
\,\,\,.
\label{spettroangoli2}
\eeq
As seen from Table III, the relative probability of these angles tend asymptotically to the common value $1/3$. The deviations from this value are in this case less than $0.003 \%$ for the first $N=5^{10}$ primes. 

\vspace{3mm}
\noindent
These examples and others give ample evidence of the equality of all relative probabilities of the appearance of the angles $\theta_{p_n}$. As a matter of fact, the equi-probability of each angle will be guaranteed by the Dirichlet theorem, as discussed in the next section.

\begin{table}[b]\label{tablefreq73}
\centering
\begin{tabular}{ |c | c | c | c | c | c | c | c | c | }
\hline\hline
N & $5^3$ & $5^4$ & $5 ^5$ & $5^6$ & $5^7$ & $5^8$ & $5^9$ & 
$5 ^{10}$ \\  
\hline
$\Prob2_1$ &
0.33600 & 0.32960 & 0.33536 & 0.33389 & 0.33343 & 0.33327 & 0.33332 & 0.33334 \\
\hline
$\Prob2_3$ & 
0.32000 & 0.33440 & 0.33280 & 0.33331 & 0.33316 & 0.33338 & 0.33337 & 0.33333
\\
\hline
$\Prob2_5$ & 
0.34400 & 0.33600 & 0.33184 & 0.33280 & 0.33341 & 0.33335 & 0.33331 & 
0.33333
 \\
\hline
\end{tabular}
\caption{The probabilities $\Prob2_a$ of angles $\alpha_a$ ($a=1, 3, 5$) relative to the 
character $\chi_3$ module 7 vs the length $N$ of the sequence $S_N$.}
\end{table}

\vspace{3mm}
\noindent
{\bf Stationarity}. It is also interesting to study the stationarity of the sequence $\hatA_{N}$. To this aim, let's consider the subsequences $A_N (\s) $ defined in (\ref{sequenceSNpart}) and  let's define the frequencies $\Prob2_k (\s) $ restricted, this time, only to these intervals  
\beq
\Prob2_k (\s) \
%prob_k 
\,=\, \frac{\# \,\{ \alpha_k \,\, \in \, A_N (\s) \}}{N} \,\,\,.
\label{frequenciesii}
\eeq
For large $N$, these frequencies seem to be reasonably ``translationally  invariant", i.e. largely independent of  the origin $\s$ of these intervals, since their relative variations of their values wrt the common asymptotic value are always order few percents, no matter how we change the origin of the intervals. An explicit example of this translation invariance of the frequencies is shown in Table IV for the angles of the character $\chi_2$ mod $q=7$.

\begin{table}[t]\label{tabble72part}
\centering
\begin{tabular}{ |c | c | c | c | c | c | c | c | c | }
\hline\hline
$\s$  & $10^5$ & $2 \times10^5$ & $3\times10^5$ & $4\times10^5$ & $5\times10^5$ & $6\times10^5$ & $7\times10^5$ & 
$8\times10^5$ \\  
\hline
${\cal P}_1(\s)$ &
0.1676 & 0.1665 & 0.1664 & 0.1655 & 0.1669 & 0.1670 & 0.1657 & 0.1668 \\
\hline
${\cal P}_2(\s)$ & 
0.1665 & 0.1665 & 0.1657 & 0.1659 & 0.1672 & 0.1677 & 0.1666 &
0.1674 \\
\hline
${\cal P}_3(\s)$ & 0.1659 & 0.1659 & 0.1659 & 0.1660 & 0.1664 & 0.1647 & 0.1667 & 0.1671 \\ 
\hline
${\cal P}_4(\s)$ &
0.1669 & 0.1667 & 0.1677 & 0.1683 & 0.1652 & 0.1679 & 0.1669 & 0.1668 
 \\
 \hline
${\cal P}_5(\s)$ & 0.1669 & 0.1670 & 0.1674 & 0.1672 & 0.1668 & 0.1673 & 
0.1661 & 0.1653
 \\
 \hline
${\cal P}_6(\s)$ & 0.1660 & 0.1658 & 0.1668 & 0.1670 & 0.1673 & 0.1657 & 0.1668 & 0.1675
 \\
\hline
\end{tabular}
\caption{The frequencies $\Prob2_k (\s)$ of angles $\alpha_k$ ($k=1,\ldots,6$) relative to the 
character $\chi_2$ module 7 for intervals of length $N = 10000$ varying the starting points $\s$ along the sequence  $A_{N} (\s) $.}
\end{table}

\vspace{3mm}
\noindent
{\bf Transition Probabilities}. Let's now make a step forward in the numerical analysis of the statistical properties of the sequence $\hatA_N$ by introducing the $k$-step probability $\Prob2_{ab}(k)$. This quantity can be defined as the number of pairs in which $\theta_{p_n} = a$ and $\theta_{p_{n+k}}=b$ divided by the number $N_a$ of instances $\theta_{p_k}=a$ that are present in the sequence $\hatA_n$. This definition implies that not necessarily $\Prob2_{ab}(k) = \Prob2_{ba}(k)$ although it is always true that $\sum_{b} \Prob2_{ab}(k) 
=1$.

\vspace{3mm}
\noindent
{\em One-step Probability Transition}. The one-step transition probability $P_{ab}(1)$ is the simplest and refers to the statistics of the next-neighboor pairs of values $(\theta_{p_i}, \theta_{p_{i+1}})$. Let's consider once again the case $q=7$ and the angles coming from the character $\chi_2$. Taking  the first $N=5^9 = 1,953,125$ primes, these are the corresponding values of $\Prob2_{ab}(1)$     
  \beq
  \label{onestepppp}
  \Prob2_{ab}(1) = 
  \left(
\begin{array}{cccccc}
 0.086860 & 0.13015 & 0.19755 & 0.22431 & 0.15683 & 0.20430 \\
 0.25019 & 0.091781 & 0.14963 & 0.15665 & 0.20558 & 0.14616 \\
 0.18018 & 0.20487 & 0.092175 & 0.20381 & 0.14614 & 0.17282 \\
 0.13538 & 0.19748 & 0.15025 & 0.087599 & 0.24923 & 0.18006 \\
 0.19660 & 0.14888 & 0.22663 & 0.13026 & 0.091846 & 0.20578 \\
 0.15061 & 0.22704 & 0.18358 & 0.19742 & 0.15004 & 0.091308 \\
\end{array}
\right)
\eeq
Observe that the entries of this matrix are not equal. Consider, for instance, the first row: this means that if at a certain point $i$ of the sequence $\hatA_n$ we have $\theta_{p_i} =\alpha_1$, there is only a $8.6 \%$ probability that the next value $\theta_{p_{i+1}}$ is still $\alpha_1$, while the most probable next angle following $\alpha_1$ is $\theta_{p_{i+1}} = \alpha_4$, whose relative probability is equal to $22.4\%$. 
\begin{figure}[b]
\centering
\includegraphics[width=0.4\textwidth]{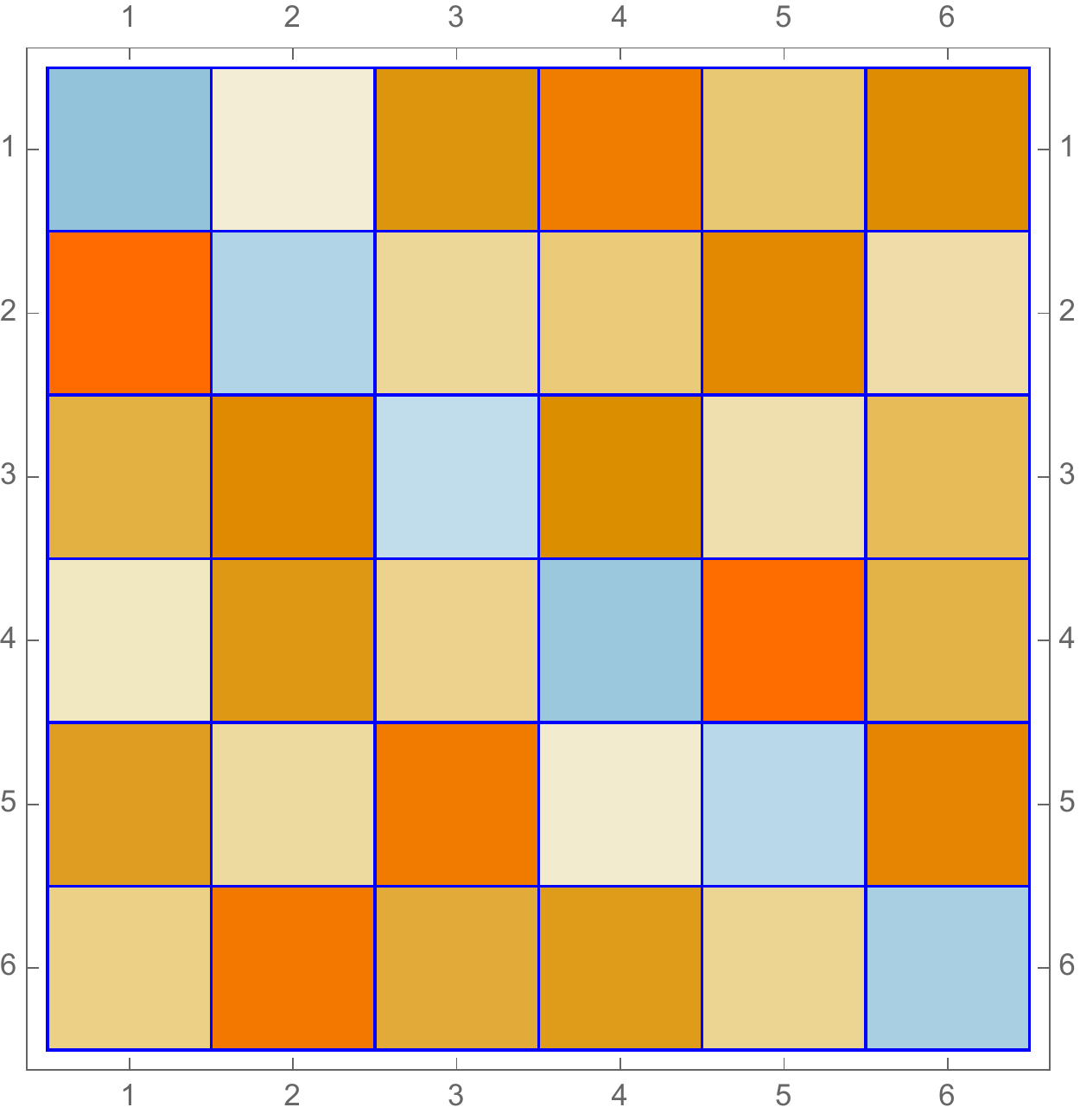}
\caption{One-step transition probability matrix for the character $\chi_2$ of modulus $q$ for $N=5^9$. Cool colors stay for low values of the matrix while warm colors for higher values.} 
\label{matrixplot1}
\end{figure}

As a general feature of this matrix (which holds for all non-principal characters of modulus $q$) we have the phenomenon of {\em anti-correlation} of equal angles, in the sense that consecutive pairs of equal angles $(\theta_{p_k} = \alpha_i, \theta_{p_{k+1}} = \alpha_i)$ are always the less probable output: correspondingly, the lowest entries of the matrix $\Prob2_{ab}(1)$ are always along the diagonal. A graphical way to show the information encoded in this matrix is shown in Figure \ref{matrixplot1} where we use cool colors for low values of the probabilities and warm colors for higher values.

\vspace{3mm}
\noindent
{\em Non-Markovian property}. Even though the entries of $P_{ab}(1)$ are different, if one takes enough large powers of this matrix one gets a {\bf constant matrix} with entries approximatively equal to $1/r$, where $r$ is the order of the character. Taking once again $q=7 $ and $\chi_2$ as example, it is enough for instance to take the $6$-th power of the matrix (\ref{onestepppp}) to get
\beq
(\Prob2(1))^6_{ab} = 
\left(
\begin{array}{cccccc}
 0.16664 & 0.16670 & 0.16664 & 0.16670 & 0.16661 & 0.16670 \\
 0.16664 & 0.16670 & 0.16663 & 0.16670 & 0.16662 & 0.16670 \\
 0.16664 & 0.16670 & 0.16664 & 0.16670 & 0.16661 & 0.16670 \\
 0.16664 & 0.16670 & 0.16663 & 0.16670 & 0.16662 & 0.16670 \\
 0.16664 & 0.16670 & 0.16664 & 0.16670 & 0.16662 & 0.16670 \\
 0.16664 & 0.16670 & 0.16664 & 0.16670 & 0.16661 & 0.16670 \\
\end{array}
\right)
\eeq
This result means that if variables $\theta_{p_i}$'s were only one-step correlated, this correlation would be essentially lost after $6$ steps, where each value becomes once again equiprobable, independent of the value of the angle assumed $6$ steps earlier. On the other hand, we can directly compute  the $6$-step transition probability and compare the two expressions. In the example at hand, 
this $6$-step probability is given by  
\beq
\label{prob6666}
\Prob2_{ab}(6) = 
\left(
\begin{array}{cccccc}
 0.16091 & 0.16817 & 0.17293 & 0.17255 & 0.16241 & 0.16304 \\
 0.17283 & 0.15798 & 0.16797 & 0.16280 & 0.17009 & 0.16833 \\
 0.16948 & 0.17138 & 0.16108 & 0.16333 & 0.16825 & 0.16647 \\
 0.16935 & 0.16235 & 0.16380 & 0.16203 & 0.17251 & 0.16996 \\
 0.16203 & 0.16786 & 0.17138 & 0.16780 & 0.15926 & 0.17167 \\
 0.16523 & 0.17248 & 0.16266 & 0.17155 & 0.16716 & 0.16091 \\
\end{array}
\right)
\eeq
Comparing $\Prob2(6)$ with $(\Prob2(1))^6$, one sees that the entries of these matrices, although quite close, are nevertheless different. This implies that, at least for a finite length $N$ of the sequence $\hatA_N$ that is sampled, the transition probabilities do not have markovian properties. Notice that also for $\Prob2(6)$  there persists the anti-correlation effect for equal values, i.e. the entries along the diagonal are smaller than the other entries. These, and other properties of the $k$-step transition probabilities, will be the content of the Lemke Oliver-Soundararajan conjecture discussed in the next section.

\vspace{3mm}
\noindent
{\bf Correlations}. The previous analysis has shown that the angles $\theta_n$ are equally distributed in the sequence $\hatA_N$ and moreover that there is an interesting pattern of correlation among the terms of this sequence. An important further insight is whether these correlations are weak or strong. To address numerically  such a question, we are going to study the correlation function at lag $j$ of the variables $c_n = \cos\theta_{p_n}$. Let's remind that for a generic time series with variables $y_k$ ($k=1,2,\ldots, N$), the correlation function at lag $j$ is defined as 
\beq
C(j) \,=\, \frac{\sum_{i=1}^{N-1} (y_i - \mu)(y_{i+j} - \mu)}{\sum_{i=1}^N (y_i - \mu)^2 }\,\,\,,
\eeq
where $\mu$ is the arithmetic mean of the time series 
\beq
\mu \,=\,\frac{1}{N} \,\sum_{i=1}^N y_i \,\,\,.
\eeq
Notice that $C(0)=1$. The {\em spectral density} of the time series is given by 
\beq
{\cal F}(k) \,=\,| \tilde C(k) |^2 \,\,\,,
\label{spectraldensity}
\eeq
where $\tilde C(k)$ is expressed by the discrete Fourier transform of the correlation function 
\beq
\tilde C(k) \,=\,\frac{1}{\sqrt{n}}  \,\sum_{j=0}^{N-1} C(j) \, e^{2 \pi i k/N} \,\,\,.
\label{spedtral density} 
\eeq

\begin{figure}[t]
\centering
\includegraphics[width=0.45\textwidth]{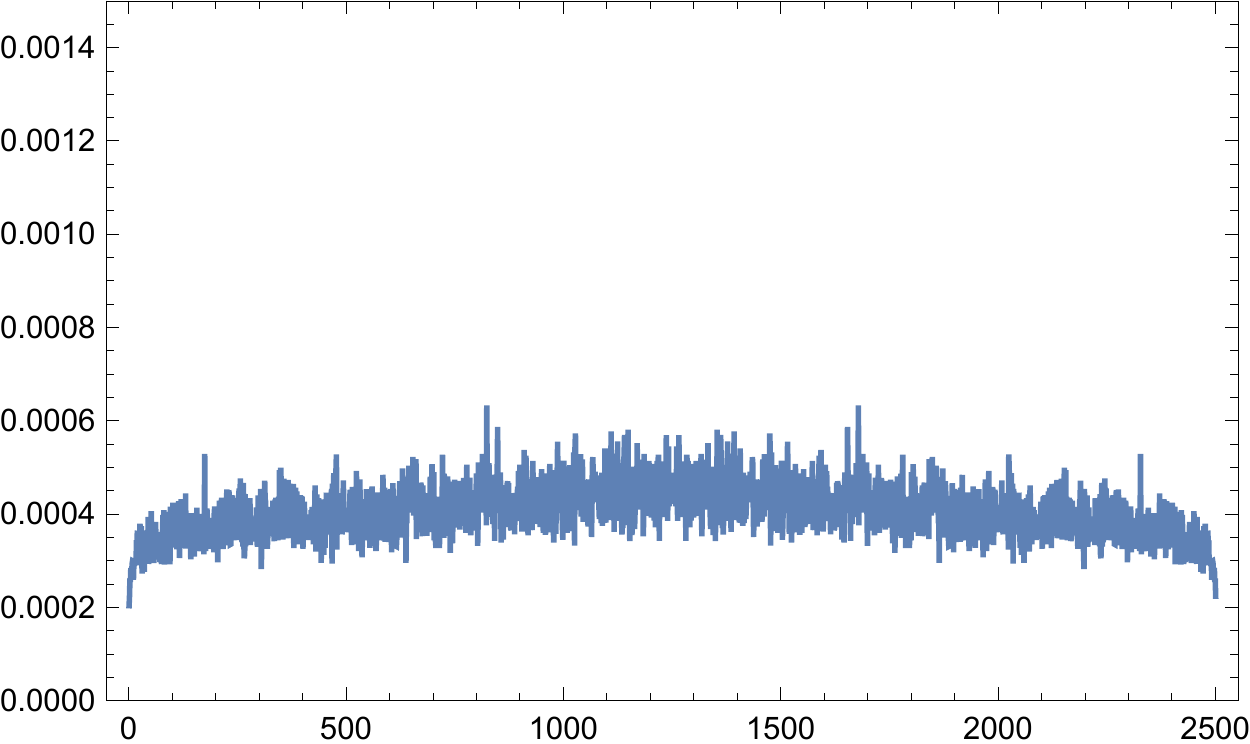}
\caption{Spectral density of the first 2500 frequencies relative to the character $\chi_2$ (mod 7) for 
the variables $c_n$ with $N=5 \times 10^5$. The flatness of this spectral density shows that the variables $c_n$ are weakly correlated.} 
\label{spectraldensityFig}
\end{figure}

\begin{figure}[b]
\centering
\includegraphics[width=0.45\textwidth]{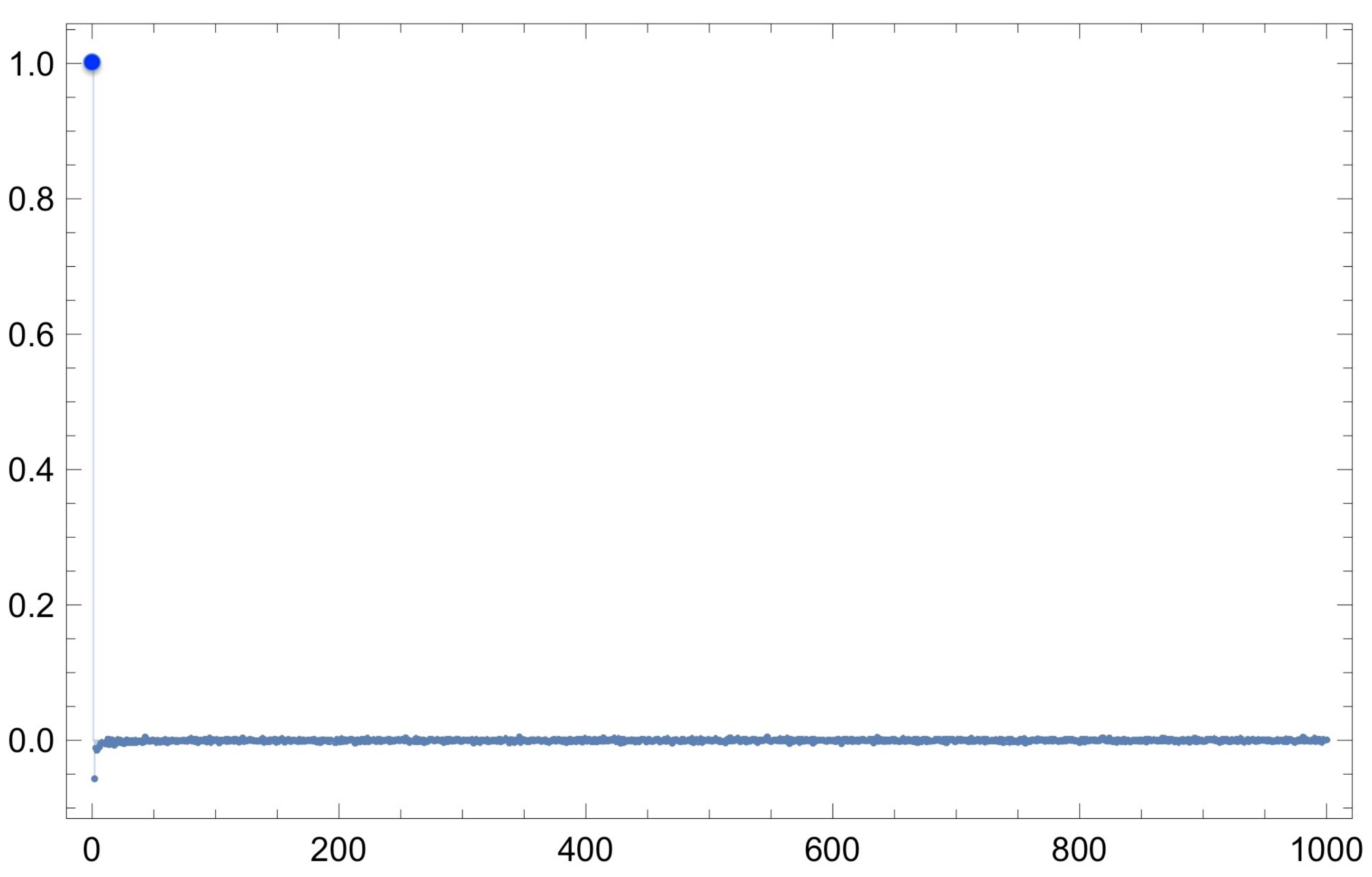}
\,\,\,
\includegraphics[width=0.47\textwidth]{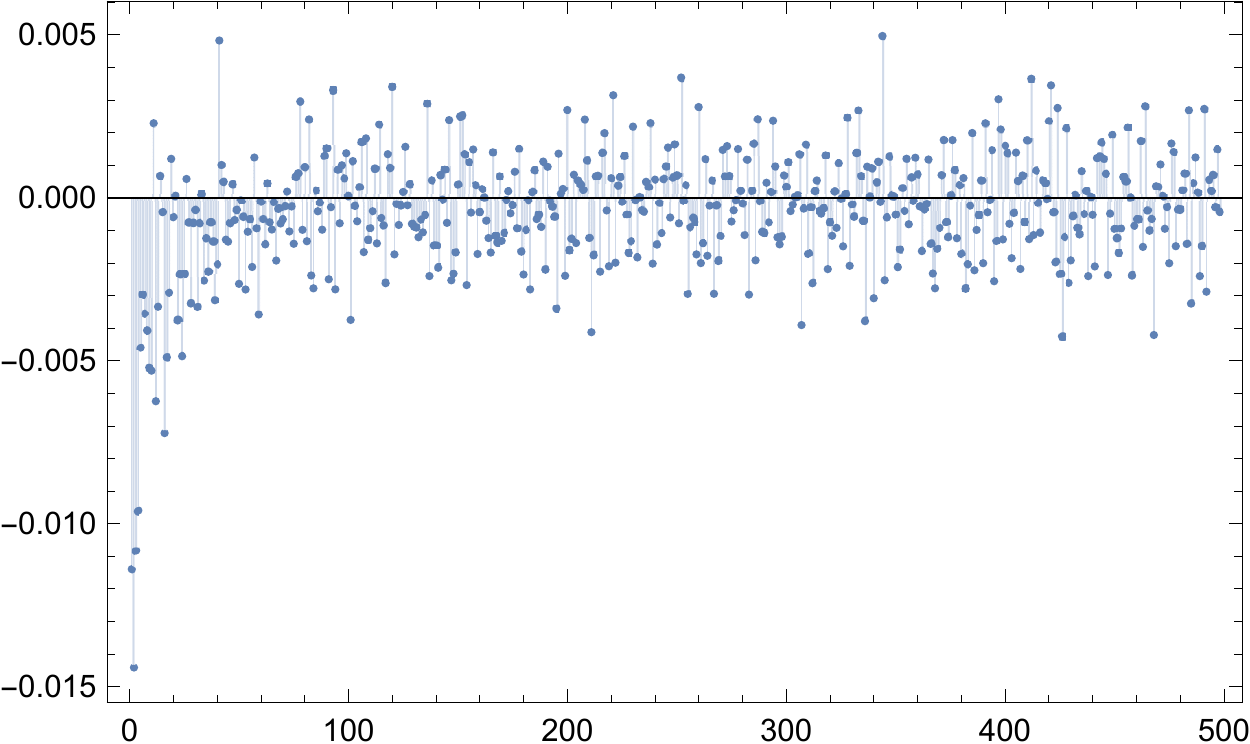}
\caption{Left hand side: correlation function $C(j)$ for the first $1000$ lags for the variables $c_n$ relative to the character $\chi_2$ (mod 7). Notice that $C(0) =1$. Right hand side: zoom on the values of $C(j)$ for $j \geq 1$. 
The small values of this function for $j \geq 1$ shows that the variables $c_n$ are weakly uncorrelated also at small separations.
For these figures, we consider $N=5 \times 10^5$ values of the $c_n$'s. } 
\label{corrlag}
\end{figure}

For an uncorrelated set of variables $y_i$, the correlation function is essentially zero for $j \neq 0$ and its 
spectral density is flat: as a rule of thumb, the  flatter the spectral density,  the less correlated are its variables.
Figure \ref{spectraldensityFig} shows the spectral density of the original variables $c_n$ for the sequence $\hatA_N$ relative to a particular character $\chi$ : such a curve is relatively flat and, correspondingly, the plot of the correlation function $C(j)$ versus the lag $j$ of the variables $c_n$ shows that, apart from  the value $C(0)=1$, for all other lags $j \geq 1$ the correlation function is extremely small (see Figure \ref{corrlag}). This result holds in general for all other sequences $\hatA_N$ relative to other characters.


\vspace{3mm}
\noindent
{\bf Summary}. Let's collect the main indications obtained by our experimental statistical analysis of the sequence $\hatA_N$. 
\begin{enumerate}
\item Moving through the sequence of the primes, the angles $\theta_{p_n}$ vary in a complicated and irregular way and 
there is an obvious analogy with the rolling of a dice with $r$ faces. 
\item As in the case of a dice, these angles seem to be equi-distributed along the sequence $\hatA_N$. 
\item There are however indications that the outputs $\theta_{p_n}$ are correlated, although weakly. The transition matrices relative to pairs of the angles $\theta_{p_n}$ separated by $k$ steps highlight an anti-correlation effect for equal angles and a non-markovian property. 
\end{enumerate}

In the next section we will see that the items $2$ and $3$ will be the content of the Dirichlet theorem and Lemke Oliver-Soundararajan conjecture respectively.

\section{Dirichlet theorem and Lemke Oliver-Soundararajan conjecture}\label{stprth}

We present here two important results which capture the statistical properties of the angles $\theta_{p_n}$. The first result concerns a theorem by Dirichlet, which originates from the interesting question whether there are an infinite number of primes in arithmetic progressions such as 
\beq
S_m \,=\, q \, m + a \,\,\,\,\,,\,\,\,\, m=0, 1, 2, \ldots 
\,\,\,\,\,\,\,\,\,\,\,\, 
q, a  \in \mathbb{N} 
\label{sequence}
\eeq
The number $q$ is the {\em modulus} while the number $a$ as the {\em residue}. As already mentioned in Section \ref{CHH}, to find a prime among the values of $S_m$ necessarily the two natural numbers $q$ and $a$ must have no common divisors, namely they should be {\em coprime}, a condition expressed as $(q,a) = 1$. Dirichlet proved that such a condition is also sufficient  \cite{Diric} and, as a consequence, there is the analog of the prime number theorem for arithmetic progressions.  Namely, define 
\beq
\label{pntq}
\pi_a (x; q )  = \# \{ p < x\, : p = a ~ {\rm mod} ~ q \}\,\,\,,
\eeq
and 
\beq
\label{nprimes}
\pi (x )  = \# \{ p < x\,\}\,\,\,.
\eeq
Then, for $x \rightarrow \infty$, Dirichlet proved that 
\beq
\label{pntq2}
\pi_a (x; q ) \sim \pi(x)/\varphi (q)\,\,\,.
\eeq
Since  
\beq
\pi(x) \sim {\rm Li}(x)
\,\,\,\,\,\,\,\,\,\,
,
\,\,\,\,\,\,\,\,\,\,
x \rightarrow \infty 
\label{asymptoticpi}
\eeq
where ${\rm Li}(x) = \int_1^x dt/\log t \sim x/\log (x)$ is the log integral function, eq.\,(\ref{pntq2}) can be also written as 
\beq 
\pi_a(x; q) \sim {\rm Li}(x)/\varphi(q) \,\,\,.
\label{asympiaq}
\eeq
Since the angles $\theta_{p}$ are functions of the residue $a$ of the prime $p$ mod $q$, Dirichlet's theorem is equivalent to  the statement that the angles $\theta_{p_n}$ are {\em equally distributed} among their possible $r$ values:   

\vspace{3mm}
\begin{theorem}\label{Dirichlettheorem} (Dirichlet)  
{\it Let $\chi(p_n) = e^{i \theta_{p_n}} \neq 0$ be a non-principal Dirichlet character modulo $q$ and $\pi(x)$ the number of primes  less than $x$. These distinct roots of unity form a finite and discrete set, i.e. $\theta_{p_n} \in \Phi = \{ \phi_1,\phi_2,\ldots,\phi_r\}$ with $r \leq \varphi(q)$ and we have
\beq 
f_i = \lim_{x\rightarrow \infty} 
\frac{\#\{ p\leq x \, : \theta_p = \phi_i\}}{\pi(x)}
= \frac{1}{r} 
\label{equipp}
\eeq
for all $i = 1, 2, \ldots, r$ where $p$ denotes a prime while $f_i$ denotes the frequency of the event $\theta_p = \phi_i$ occurring.
}
\end{theorem}

\vspace{3mm}

It is important to notice that the Dirichlet theorem does not say anything about the possible correlations of the angles in the sequence $A_N$.  For example, as shown in the previous section, correlations of these variables is probed by how many times the pairs $(\phi_a,\phi_b)$ appear as values of two consecutive angles $\theta_{p_n}$ and $\theta_{p_{n+1}}$, or angles separated by $k$ steps in the sequence $A_{N} $, i.e. $\theta_{p_n} $ and $\theta_{p_{n+k}}$. The theoretical formulation of this problem has been recently addressed by Lemke Oliver and Soundararajan on the basis of the Hardy-Littlewood prime k-tuples conjecture. 
Let's notice that in the paper \cite{OliverSoundararajan}, instead of the angles $\theta_{p_n}$, Lemke Oliver and Soundararajan were directly concerned with the patterns of residues mod $q$ among the sequences of consecutive primes less than an integer $x$ (on this subject see also \cite{Shiu,Ash} and references therein). For our purposes, this is equivalent to the correlations among the angles $\theta_{p_n}$ since these quantities are just functions on the residues.  We will refer to the residues as ``$a$" (or ``$b$") in accordance to \eqref{sequence}:  
\beq
\label{resid}
p_n = a ~ {\rm mod} ~ q, ~~~a \in \{ 0,1,2, \ldots , \varphi(q)  \} \,\,\,. 
\eeq
If $q$ is not a prime, then not all values of $a$ in the above set are realized. If $q$ is instead a prime, then 
for $p_n > q$,  there are $\varphi (q) = q-1$ possible values of the residue $a$ and only for the special case when $p=q$ is the residue equal to $0$. Hereafter we focus our attention to $q$ equal to a prime and to the counting functions 
\beq
\pi_{ab} (x; q, k) \,= \#\{ p_n < x  : 
p_n \equiv a\,\,({\rm mod}\,\,q)\,\,,\,\, p_{n+k} \equiv b \,\,({\rm mod}\,\,q) \}\,\,\,. 
\eeq
For instance,  for $q=3$ and $k=1$,  $\pi_{ab}$ counts the number of consecutive primes whose residues have the patterns $(a,b) = (1,1), (1,2), (2,1), (2,2)$. Based on the pseudo-randomness of the primes, for $x \rightarrow \infty$ one would expect that the primes counted by $\pi_{ab}(x; q, k)$ go as 
\beq 
\pi_{ab}(x;q,k) \sim \pi(x)/(\varphi(q))^2\,\,\,,
\label{largeasymptotic}
\eeq
independent of the separation $k$ of the two residues. However, as shown by  Lemke Oliver and Soundararajan, for {\em finite values of $x$} there are potentially large corrections in the expected asymptotic behavior which create biases toward certain patterns of residues. In the following, in particular, we focus our attention on the $\varphi (q) \times  \varphi (q) $ matrices\footnote{LOS define $f_{ab}$ as above but with $\pi (x)$ replaced by the log integral ${\rm Li} (x)$.  The latter is simply the leading approximation to $\pi (x)$ based on the prime number theorem, thus our definition is actually more meaningful. In the large $x$ limit,  the results are the same whether one uses $\pi (x)$ or ${\rm Li }(x)$.}
\beq
f_{ab}(x,q,k) \,=\, \frac{\#\{p_n \leq x : p_n \equiv a \,\,({\rm mod}\,\,q)\,\,,\,\, p_{n+k} \equiv b \,\,({\rm mod}\,\,q) 
\}}{\pi(x)} \,\,\,, 
\label{definitionfreqOliver}
\eeq
which, for nearby $x$, give the {\em local densities} of pairs of primes, in which $p_n \equiv  a \,\,({\rm mod}\,\,q) $ will be followed, after $k$ steps, by a prime $p_{n+k} \equiv b \,\,({\rm mod}\,\,q)$. Here we quote the large $x$ behavior of these quantities \cite{OliverSoundararajan}: 

\vspace{1mm}

\begin{conjecture}(Lemke Oliver-Soundararajan). 
{\it For large values of $x$ we have\footnote{It is possible to express $f_{ab}(x,q,1)$ and $f_{ba}(x,q,1)$ individually (and they are not equal) but their expression is rather complicated, see \cite{OliverSoundararajan}. Moreover, the expressions (\ref{firstconjab}) and (\ref{firstconjaa}) given here are those of LOS  but specialised to the modulus $q$ being a prime.} 
\beq 
f_{ab}(x,q,1) + f_{ba}(x,q,1) \,=\,\frac{2}{(\varphi(q))^2} \left[1 + \frac{\log\log x}{2 \log x} - \log\frac{q}{2\pi} 
\frac{1}{2 \log x} + O\left(\frac{1}{(\log x)^{7/4}}\right)\right]\,\,\,,
\label{firstconjab}
\eeq
whereas 
\beq
f_{aa}(x,q,1) \,=\, \frac{1}{(\varphi(q))^2} \left[1 - \frac{(\varphi(q) -1)}{2} \, \frac{\log\log x}{\log x} + 
(\varphi(q) -1) \,\log\frac{q}{2\pi} \,\frac{1}{2 \log x} +
O\left(\frac{1}{(\log x)^{7/4}}\right)\right]\,\,\,.
\label{firstconjaa}
\eeq
} 
\end{conjecture}

\vspace{1mm}

\begin{conjecture}(Lemke Oliver-Soundararajan). 
{\it  For  $k \geq 2$,  then for large values of $x$ we have 
\beq 
f_{ab}(x,q,k)  \,=\,\frac{1}{(\varphi(q))^2} \left(1 + \frac{1}{2 (k-1)} \,\frac{1}{\log x} + O\left(\frac{1}{(\log x)^{7/4}}\right)\right)\,\,\,,
\label{secconjab}
\eeq
whereas 
\beq
f_{aa}(x,q,k) \,=\, \frac{1}{(\varphi(q))^2} \left(1 - \frac{(\varphi(q) -1)}{2(k-1)} \, \frac{1}{\log x} + O\left(\frac{1}{(\log x)^{7/4}}\right)\right)\,\,\,.
\label{secconjaa}
\eeq
} 
\end{conjecture}
 
\vspace{1.5mm}

The opposite signs in the second term in \eqref{firstconjaa} versus \eqref{secconjab} are responsible for the bias and the anti-correlations that we saw from the numerical studies of the previous section. Notice that the formulas above present a permutation symmetry ${\mathcal S}_{\varphi(q)}$ (since the only thing that matters is whether the residues are equal or different) which, for a matrix as $\Prob2(6)$ computed with $N = 5^9$ and reported in eq.\,(\ref{prob6666}), was already verified with a precision of the order $2\%$. Notice that the counting functions $ \pi_{ab}(x;q,k)$ of the pairs of primes in $\hatA_N $ relative to various residues are given by
\beq
\label{generaldistrrr}
 \pi_{ab}(x;q,k)  =  \
\pi (x) \, f_{ab}(x,q,k) \, \sim \frac{x}{\log x}  f_{ab} (x, q, k)
\eeq
Let's further comment on some important features which emerge from these functions $\pi_{ab}(x;q,k)$. For $x\to \infty$, these formulas state that all pairs of residues in $\hatA_N  $ are {\em equally probable} (both for consecutive primes and primes separated by $k$ steps) and their probability is given by $1/(\varphi(q))^2$. This means that, in the limit $x \rightarrow \infty$, the angles $\theta_{p_n}$ in any subsequence $\hatA_N $ are completely {\em uncorrelated} and this is the most important property for  the aim of establishing the GRH! However, at any {\em finite} value of $x$, the next neighbor variables in the subsequence $\hatA_N $ tend to be {\em anti-correlated}, as we already noticed: the occurrence of pairs of equal residues $(a,a)$ for next neighbor primes are always less probable than the occurrence of pairs of different residues $(a,b)$ although this may be considered a finite-size effect, since it 
vanishes as $\sim \log\log x/\log x$. At any {\em finite} $x$, this anti-correlation phenomenon also persists for primes which are separated by $k$ steps and the matrices $f(x,q,k)$ are {\em not} equal to $(f(x,q,1))^k$, i.e. these probabilities do not satisfy the Markovian property, as also we noticed earlier. This correlation decreases as $1/k$ with the separation $k$ of the two primes but it is also a finite size effect since the coefficient in front of this $1/k$ correlation vanishes as $1/\log x$ when $x \rightarrow \infty$. The Markovian property of these matrices is of course restored in the $x \rightarrow \infty$ limit.

\def\chat{\hat{c}}
\def\Chat{\hat{C}}
\def\Ghat{\hat{G}}

\def\hatA{A} 

\section{A normal distribution for the series $C_N$}\label{Timeseries}

Let's recall that our aim is to estimate how the series $C_N$ grows with $N$. If we want to view it as a random time series, we have to face the problem of defining an ensemble $\CE$ relative to the possible values of $C_N$ together with their relative probabilities. An obvious obstacle is that, for any given character, there is of course {\em one and only one} series $C_N$. This, however, is a common problem in many time series, in particular for all those that refer to situations for which it is impossible to \textquotedblleft turn back time\textquotedblright. Indeed, in these cases it is impossible to have access to all possible outputs and therefore equally impossible to define the relative probabilities. In the literature, this is known as the {\em Single Brownian Trajectory Problem} (see, for istance \cite{brow1,brow2,brow3,brow4} and references therein). 

In order to deal with this problem, we can consider an arbitrarily long time series and, in order to sample it, take \textquotedblleft stroboscopic\textquotedblright \,snapshots of it, in the following way. Define the  {\it ordered} intervals of length $N$ starting at $\s$  
\beq
I_N(\s) =\{\s, \s+1, \s+2, \ldots, \s+N -1 \} \,\,\,, 
\,\,\,
\label{intervals}
\eeq 
and  the associated angles $\hatA_{N}(\s) $ 
\beq
\hatA_{N} (\s)  \,=\, 
\{\theta_{p_n}\, ; ~ n \in I_N(\s)\}\,\,\,. 
\label{sequenceSNpart}
\eeq
We then define {\it block variables} $C_N (\s)$ based on the above intervals:
\beq
C_N(\s) \,=\,\sum_{k\in I_N(\s)}  c_{k} = \sum_{k=\s}^{\s + N -1}  \cos \theta_{p_n}  \,\,\,. 
\label{groupvariables}
\eeq
For reasons that will become clear,  it will be convenient to also define
\beq 
A_{N_1,N_2} \,  \equiv   A_{N_2 - N_1 +1} (N_1)  = \, \{\theta_{p_n}\, ; ~ n=N_1, N_1+1,\ldots , N_2 \}\,\,\,,
\label{sequenceSN}
\eeq
relative to primes between $p_{N_1}$ and $p_{N_2}$. Choosing 
\beq
N_1 < \s < \s + N -1 < N_2\,\,\,\,,
\label{appartenenza}
\eeq
$\hatA_{N}(\s)$ is of course a subset  of $A_{N_1,N_2}$.  Imagine we fix a very large value of $N_1$ and then vary $N_2$: in this way we can consider arbitrarily long sequences $A_{N_1,N_2}$, out of which many and well separated block variables $C_N(\s)$ of the {\em same length} $N$ can be defined  and used as members of the ensemble to which belongs the original sequence $C_N(1)$! This is equivalent to the stroboscopic snapshots behind the solution of the the Single Brownian Trajectory Problem (see the forthcoming subsection).  The validity of this self-averaging procedure relies on two aspects of the corresponding time series: its ergodicity and stationarity. Let's discuss these two aspects separately. 

In the case of our sequences $A_{N_1,N_2}$, their ergodicity is simply guaranteed by the presence of {\em all} possible outputs of the angles $\theta_{p_n}$ along the sequence of the primes. Their stationarity is an issue more subtle which can be settled however on the basis of the following considerations. According to the formulas of LOS , there are correlations which explicitly depend on the point $x$ along the sequence of the primes and therefore, for arbitrary values of the extrema $N_1$ and $N_2$, they break -- strictly speaking -- the stationarity of the sequences $A_{N_1,N_2}$. There are however two facts which help in solving this issue: the first is that, as we already commented, these are {\em finite size effects} which vanish when $x\rightarrow \infty$; the second is the equivalence between the series 
\beq
\label{importantequivalence}
\sum_{n=1}^N c_n \sim \sum_{n=\s}^{N} c_n, ~~~~~{\rm for} ~ N \to \infty
\eeq
which holds since we are interested in their behavior only for $N\rightarrow \infty$ and which implies that we have always the freedom to drop a finite number $l$ of the first terms of the series $C_N$. Thanks to this equivalence, even at finite $x$ we can focus our attention on sequences whose extrema $N_1$ and $N_2$ are such that the correlations are both weak and sufficiently uniform along the entire length of these intervals. Intervals $(N_1,N_2)$ which satisfy this property will be called {\em inertial intervals} and sequences based on these intervals can be made as stationary as one desires. For instance, choosing $N_1 = 10^{200}$ and $N_2 = 10^{250}$, the correction to a uniform background $1/(\varphi(q))^2$ distribution is only of the order $0.20 \%$ and $0.17\%$ respectively at the beginning and at the end of the sequence  $A_{N_1,N_2}$, therefore with a breaking of the stationarity that can be quantified of the order of $0.03 \%$. These values come from the correction $1/\log x$ present in the LOS with respect to the constant values of the correlations, computed for $x =N_1$ and $x=N_2$. Of course we can choose arbitrarily higher values of $N_1$ and $N_2$ (since the primes are infinite) and make the corresponding sequence $A_{N_1,N_2}$ stationary with arbitrarily higher degree of confidence. By the same token, namely enlarging the size of the sequences $A_{N_1,N_2}$, we can always set up a proper ensemble for $C_N(1)$ for {\em any} $N$, no matter how large. Notice that as $N_{1}$ and $N_2 \to \infty$, also $\s \to \infty$. Moreover, we are going to assume the inequalities  
\beq
\label{ineq}
1 \ll N \ll \s,  
\eeq
so that  $p_\s \approx p_{\s+N}$.    

\subsection{Statistical Ensemble $\CE$ for  the series $C_N$}\label{subesemble and its variance}
\begin{figure}[t]
\centering
\includegraphics[width=0.40\textwidth]{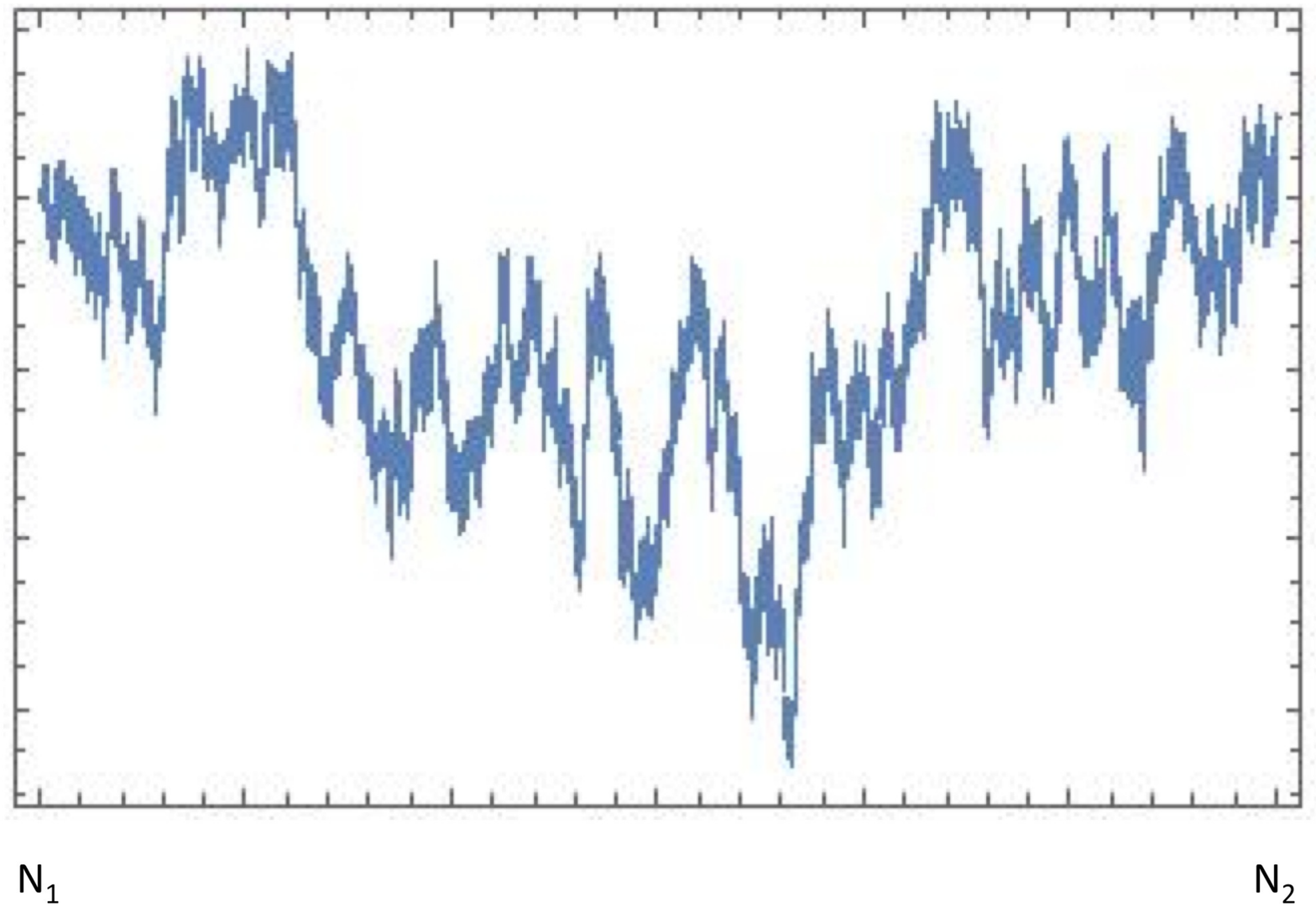}
\,\,\,
\includegraphics[width=0.42\textwidth]{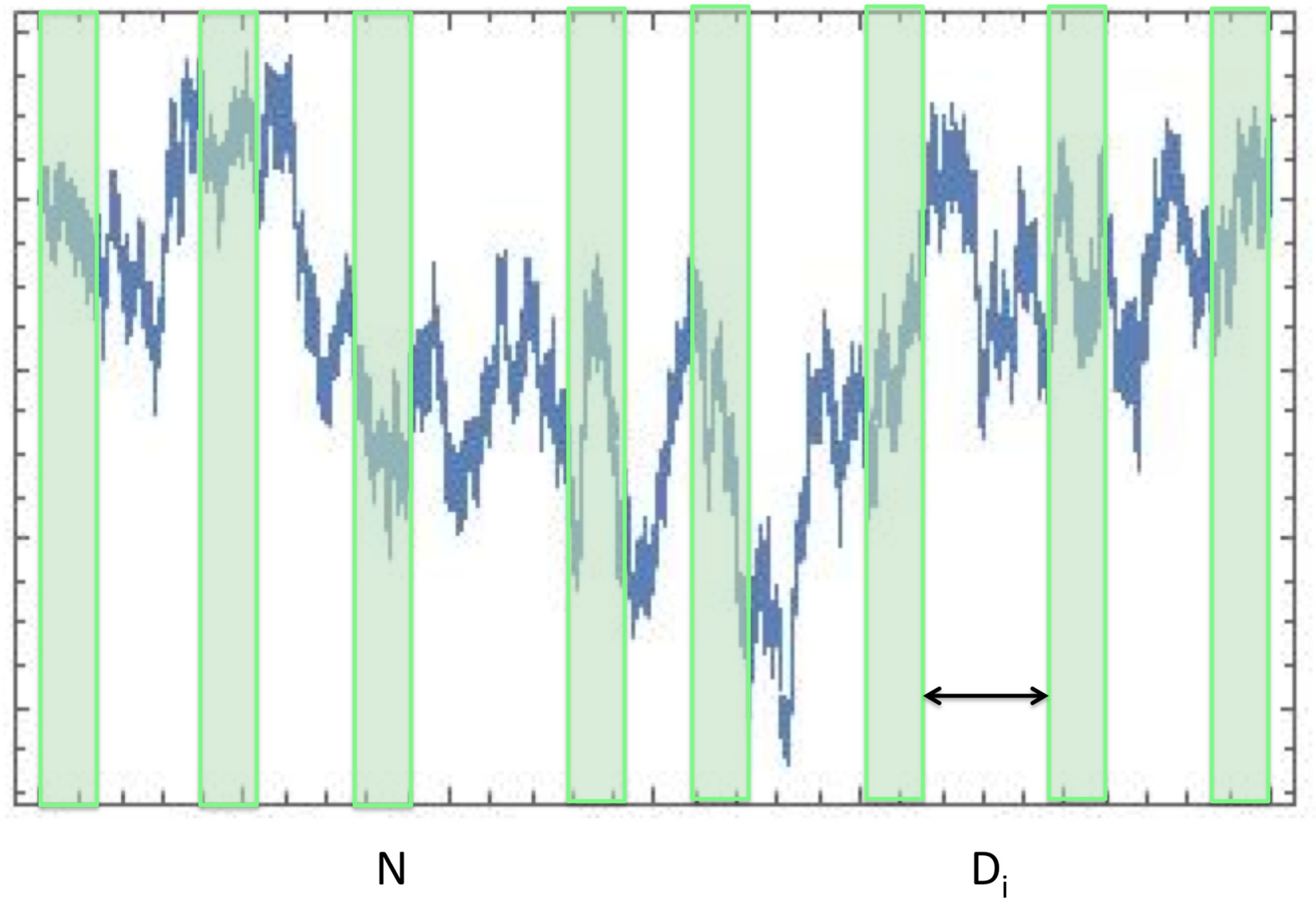}
\caption{Left hand side: series $\sum_{n=N_1}^{N} c_n$ vs $N$, in the inertial interval $(N_1,N_2)$. 
Right hand side: sampling of the time series done in terms of block variables $C_N (\s)$ of length $N \gg 1$ relative to the green intervals separated by distances $D_i$. Under the hypothesis of stationarity of the sequence $A_{N_1,N_2}$, the values of these blocks define a probabilistic ensemble $\CE$ for the quantity $C_N$ relative to the sum of the first $N$ values $c_n$. } 
\label{partitiontimeseriesGN}
\end{figure}

The block variables $C_N(\s)$ are the equivalent of the  
\textquotedblleft stroboscopic\textquotedblright \,images of length $N$ of a single Brownian trajectory (see Figure \ref{partitiontimeseriesGN}) and they allow us to control the irregular behavior of the original series $C_N(1)$ by proliferating it into a collection of sums of the same length $N$. It is this collection of sums that forms the {\em set of events}, i.e. the ensemble $\CE$ relative to the sums of $N$ consecutive terms $c_n$.  As discussed originally in \cite{shortpaper}, this ensemble is defined as follows:

\begin{enumerate}
\item Consider two very large integers $N_1$ and $N_2$  (which eventually we will send to infinity), with $N_1\gg 1$, $N_2\gg1$  but also $L \equiv (N_2 - N_1) \gg 1$ such that, for a given character $\chi$ of modulus $q$,  
the sequence $A_{N_1,N_2}$ is inertial. 
\item For any fixed integer $N$, with $1 \ll N \ll L$, consider the union of  sets 
\beq 
{\cal S}_M \,=\, \bigcup_{i=1,\ldots M} I_N(i) 
\,\,\,\,\,\,\,\,
,
\,\,\,\,\,\,\,\,
N_1 \leq i < N_2 
\,\,\,\,\,\,,\,\,\,\,\, I_N(i) \cap I_N(j) = 0 
\,\,\,,\,\,\,  i \neq j 
\eeq
made of $M$ {\em non-overlapping} and also well separated intervals of length $N$ whose origin is between the two large numbers $N_1$ and $N_2$ (see Figure \ref{partitiontimeseriesGN}). The integer $M$ is the cardinality of the set $S_M$.    These conditions ensure that the block variables $C_N(i)$ computed on such disjoint intervals are very weakly correlated and therefore we can assume that we are dealing statistically with $M$ separated copies of the original series $C_N(1)$. 
\item At any given $N_1$ and $N_2$, the cardinality ${\rm card} ({\cal S}_M) = M$ of these sets cannot be larger of course than $L/N$. There is however a large freedom in generating  them: 
\begin{enumerate}[label={\bf \alph*}]
\item We can take, for instance, $M$ intervals $I_N(\s)$ separated by a fixed distance $D$, with the condition that $M(N + D) = L$; 
\item Alternatively,  we can take,  $M$ intervals $I_N(\s)$ separated by random distances $D_i$ such that $M N + \sum_{i=1}^M D_i = L$. 
\end{enumerate}
\item The ensemble $\CE$ is then defined as  the set of the $M$ block variables $C_N(\s)$ relative to the intervals $I_N(\s) \in {\cal S}_M$:
\beq
\label{CEdef}
\CE = \{ C_N (\s) \},  ~~~~{\rm with~~}  I_N (\s) \in {\cal S}_M
\eeq
\end{enumerate}
\noindent 
In summary, choosing two very large and well separated  integers $N_1$ and $N_2$, we can generate a large number of sets of intervals ${\cal S}_M$ and use the corresponding block variables of length $N$ to sample the typical values taken by a series consisting of a sum of $N$ consecutive terms $c_n$. In view of the ergodicity and stationarity of the sequence $A_{N_1,N_2}$ for $N_1\rightarrow \infty$ and $N_2\rightarrow \infty$, this is equivalent to determining  the statistical properties of the original series $C_N$.

\vspace{2.4mm}

The most important quantities of the series $C_N(\s)$ are its mean and variance. In particular, the value of the mean is a simple consequence of the Dirichlet theorem, as show hereafter.

\vspace{3mm}
\noindent
{\bf Mean of $C_N$}. In the limit $N\rightarrow \infty$, the series $C_N$ has zero mean 
\beq 
\mu\, \equiv \, \lim_{N\rightarrow \infty} \frac{1}{N} 
\sum_{n=1}^N \cos\theta_{p_n} \,=\, 0 \,\,\,.
\label{meanC_N}
\eeq
The proof is quite simple.  Consider the case when the cardinality $r$ of the set $\Phi$ of the angles coincides with $\varphi(q)$, i.e. $r = \varphi(q)$ (recall the definition 
of the angles $\alpha_k$ given in \eqref{notationangles}).     We can use then eq.\,(\ref{pairwiseangles}) to group pairwise the terms of the sum and since  
\beq
\cos(\alpha_{\varphi/2 +k}) = - \cos\alpha_{k}
\,\,\,\,\,\,
,\,\,\,\, 
k=1,\ldots ,\varphi(q)/2 \,\,\,,  
\label{pairwisecosine}
\eeq
we have 
\beq 
\mu \,=\,\, \lim_{N\rightarrow \infty} \frac{1}{N} 
\sum_{n=1}^N \cos\theta_{p_n} \,=\, 
 \sum_{k=1}^{\varphi(q)/2} \cos\theta_k \,
 \left(f_k - f_{\varphi/2+k}\right)
 \,=\, 0\,\,\,
\label{meanC_N11}
\eeq
since, in the $N\rightarrow \infty$ limit, from  the Dirichlet theorem all frequencies $f_{n}$ are equal. Analogous results can be easily obtained also when $r < \varphi(q)$. In the double limit $N_1 \rightarrow \infty$ and $N_2 \rightarrow \infty$ (so that also $N \rightarrow \infty$), from the stationarity properties of the sequence $A_{N_1,N_2}$ the same is true for the ensemble average of the large $N$ block variables $C_N(\s)$ 
\beq 
\Ex \Bigl[C_N  (\s )\Bigr] \,=\, 0 \,\,\,.
\eeq
In conclusion, the ensemble $\CE$ consists of block variables $C_N(\s)$ equally distributed among positive and negative values.

\vspace{3mm}
\noindent
{\bf Variance of the block variables $C_N(\s)$}. The block variables $C_N(\s)$ are defined in eq.\,(\ref{groupvariables}).
Let us  first define the variance $b^2$ of the cosine on the set of the $r$ angles 
\beq
\label{variancecosine}
b^2 \,\equiv\,\frac{1}{r} \sum_{k=1}^r \cos^2\phi_k \,=\, 
\left\{
\begin{array}{lll} 
1 & \,\,\,\,,\,\,\,\, & {\rm \,\,}   {\rm if ~ \chi ~ is ~ real} \\
1/2 &  \,\,\,\,,\,\,\,\, & {\rm  \,\,} {\rm if ~ \chi ~ is ~ complex } 
\end{array}
\right.
\eeq 
If  $\chi$ is real, then the only values of the character are $\chi = \pm 1$. 

If the terms $c_n$ entering the block variables $C_N(\s)$ were uncorrelated, the probability distribution of these block variables could be computed in terms of the characteristic function $\hat P(k)$ of the variable $c\equiv \cos\theta$ given by 
\beq 
P (x) \,=\,\frac{1}{2\pi} \int_{-\infty}^{\infty} dk (\hat P(k))^N \, e^{-i k x} \,\,\,.  
\label{PX}
\eeq
This expression would have led immediately to the gaussian behavior relative to the central limit theorem, since 
\beq 
\hat P(k) \,\simeq 1 - b^2\,\frac{k^2}{2} + \cdots \,\,\,,
\eeq
with $b^2$ given in eq.\,(\ref{variancecosine}), and for large $N$ 
\beq
P (x) \,=\,\frac{1}{2\pi}\,\int_{-\infty}^{\infty} dk \,e^{N \log \hat P(k)} \, e^{-i k x} \,
\simeq \frac{1}{2\pi}\,\int_{-\infty}^{\infty} dk \,e^{-N b^2 \,k^2/2}\ \, e^{-i k x} \,
\simeq e^{-x^2/(2 N b^2)}\,\,\,. 
\eeq
So, if the $c_n$'s were uncorrelated, the block variables $C_N(\s)$ would be certainly  gaussian distributed with a variance equal to $N$ times the variance $b^2$ of the $c_n$'s. 

However, for any sequence $A_{N_1, N_2}$, the variables $c_n$ {\em are} weakly correlated,  as indicated by the LOS conjecture. A priori, these correlations do not prevent to have a central limit theorem, as we are going to show. We can use the correlations between the variables $c_n$  to compute the variance $\sigma_N^2$ of the block variables\footnote{Here we present the argument relative to the case $r=\varphi(q)$ but the final expression of the variance, eq.\,(\ref{mostimportantformula}), holds for all cases. Moreover, in the following we will consider block variables nearby the position ${\bf \s}$, with $N_1 \leq \s < N_2$. 
} $C_N(\s)$. To this aim, consider the block variable $C_N(\s)$ belonging to the ensemble $\CE$ defined above and take the ensemble average of its square 
\begin{eqnarray}
\sigma_N^2  (\s) & \,=\, & \Ex \Bigl[ (C_N(\s))^2 \Bigr] \,=\,\, \sum_{l=0}^{N-1}  \sum_{m=0}^{N-1} 
\Ex \Bigl[ c_{\s + l } c_{\s + m } \Bigr]\,=\,\sum_{m=0}^{N-1} \Ex \Bigl[ c^2_{\s+m}\Bigr] + 
2 \, \sum_{m=1}^{N-1} [N - m] \,\Ex \Bigl[ c_{\s} c_{\s+ m} \Bigr] \,\,\,\nonumber \\
&\,=\,&\, N \, \Ex \Bigl[ c^2_\s \Bigr] + 2 \, \sum_{m=0}^{N-1} (N - m) \Ex \Bigl[ c_{\s } c_{\s + m} \Bigr] \,\,\,,
\end{eqnarray}
where we used the stationarity of the ensemble to group the contributions of the various pairs separated by $k$ steps (there are $(N-m)$ of them). Isolating further the term $m=1$ in the second quantity of the expression above, we have that the variance can be expressed as 
\beq
\Ex \Bigl[ (C_N(\s))^2 \Bigr] \,=\, {\cal D}_0 + {\cal D}_1 + {\cal D}_2 
\eeq
where 
\begin{eqnarray}
{\cal D}_0 \,& = & \,  N \, \Ex \Bigl[ c^2_\s\Bigr] \,\,\,,\nonumber \\
{\cal D}_1 \,& = & \,  2  (N - 1) \,\Ex \Bigl[ c_{\s} c_{\s + 1} \Bigr] \,\,\,,\\
{\cal D}_2 \,& = & \, 2 \sum_{m=2}^{N-1}  [N - m] \,\Ex \Bigl[ c_{\s} c_{\s + m} \Bigr] \,\,\,.\nonumber
\end{eqnarray}
The variables $c^2$ are statistically equi-distributed on the $r$ angles and their variance $b^2$ on these angles was given in eq.\,(\ref{variancecosine}), so ${\cal D}_0$ is expressed as  
\beq
{\cal D}_0 \,=\,b^2 \,  N \,\,\,. 
\label{d00}
\eeq
In order to compute ${\cal D}_1$, we need the formulas (\ref{firstconjab}) and (\ref{firstconjaa}) relative to the residues of two next neighbor primes. In light of eqs.\,(\ref{firstconjab}) and (\ref{pairwisecosine}), notice that in  the average of the product of the two cosines on the ensemble, keeping initially the $\cos\theta_{p_i}$ fixed, there are only {\em two terms} which contribute to the average: the first when $\theta_{p_{i+1}}  \,=\,\theta_{p_i}$ (with weight $f_{aa}(p_\s,q,1)$), the second when $\theta_{p_{i+1}} \,=\,\theta_{p_i} + \pi$ (with weight $f_{ab}(p_\s,q,1)$), while all other terms cancel out pairwise. Summing now on the $\varphi(q)$ values taken by $\theta_{p_i}$, we have 
\beq
{\cal D}_1 \,=\, - b^2\, (N-1)  \,\frac{\log\log p_\s}{\log p_\s} +b^2\, \log\frac{q}{2\pi} \frac{1}{\log p_\s}\,\,\,. 
\label{d11}
\eeq

The calculation is essentially similar for the other term ${\cal D}_2$, the only difference between the dependence of the separation $m$ of the two cosines
\beq
{\cal D}_2 \,=\, - b^2 \,\frac{1}{\log p_\s}\,\sum_{m=2}^{N-1} (N - m) \frac{1}{m-1}
\,\,\,.
\eeq
Putting together the three terms, we arrive to the following theorem: 
\vspace{3mm}
\begin{theorem}\label{GM1}   
{\it Assuming the validity of the LOS  conjectures, in the inertial intervals 
the variance $\sigma_N^2$ of the the block variables of length $N$ is given by  
\beq
\sigma_N^2(\s)/b^2\,=\, \Ex \Bigl[ (C_N(\s))^2 \Bigr]/b^2 \,=\,  N \,\lambda(N,\s) \,  + \,\rho(N,\s) \,\,\,,
\label{mostimportantformula}
\eeq
where 
\begin{eqnarray}
\lambda(N,\s) &\,=\,& \left[1 +\frac{1}{\log p_\s} \left(1 - \sum_{m=1}^{N-2} \frac{1}{m}\right) - \frac{\log\log p_\s}{\log p_\s}\right] \,\,\,,
\label{correctionfactorlambda}\\
\rho(N,\s) &\,=\,& \, \frac{1}{\log p_\s} \left[\log \left(\frac{q \log p_\s}{2 \pi e^2}\right) + \sum_{m=1}^{N-2} \frac{1}{m}
\right] 
\label{interceptrho}\,\,\,.
\end{eqnarray}
}
\end{theorem}

\vspace{3mm}
This theorem states that, in all the inertial intervals, the variance of block variables $C_N(\s)$ of length $N$ scales {\em linearly} with N, up to a correction factor $\lambda(N,\s)$ which is {\em independent} of the modulus $q$, but depends on the prime $p_{\s}$ of the inertial interval around which we consider the block variables. Notice that, for $\s\rightarrow \infty$, we recover a purely gaussian expression for  the variance 
\beq
\lim_{\s \rightarrow \infty} \sigma^2_N (\s)  \,=\, b^2 \, N \,\,\,.  
\label{limitxinfinity}
\eeq
Notice that we have also used the inequality \eqref{ineq} which implies $p_{N+\s} \approx p_\s$. Keeping instead $\s$ finite and considering the large $N$ asymptotic of this formula, the factor $\lambda(N,\s)$ introduces a logarithmic correction in the variance since 
\beq
\sum_{m=1}^{N-1} \frac{1}{m} \simeq \log N + \gamma_{E} \,\,\,,
\eeq
where $\gamma_E$ is the Euler-Mascheroni constant. Notice that, as far as $\s$ is finite, for the anti-correlation of the residues of consecutive primes, we have $\lambda(N,\s) < 1$ and therefore the variance of the block variables $C_N(\s)$ at a finite $\s$ is always {\em smaller} than the variance of $N$ uncorrelated variables.

As anticipated, Theorem \ref{GM1},  along with  \eqref{limitxinfinity}, implies that in the limit $\s \to \infty$ the properly normalized block variables are gaussian distributed: 
\beq
\frac{C_N(\s)}{\sigma_N(\s)} ~ \dist ~  \CN(0,1)\,\,\,, 
\label{finallynormal}
\eeq
where finite $\s$  corrections to  $\sigma_N(\s)$  are given in eq.\,(\ref{mostimportantformula}).  For finite $\s$, the distribution is not purely gaussian, and this non-gaussianity captures the existence of correlations between the primes for a given character.

\vspace{1.5mm}

In light of this result, taking larger and larger inertial intervals, and therefore correspondingly larger and larger values of $N$,  the block variables $C_N(\s)$ of length $N$ always scale as 
\beq
C_N(\s) \,= O(N^{1/2+\epsilon}) \,\,\,, 
\label{fffff}
\eeq
for arbitrarily small $\epsilon >0$.  
Notice that, in probabilistic language, in the limit $N \rightarrow \infty$ this behavior occurs  with probability equal to 1. Indeed, since $\sigma_N(\s) \leq \sqrt{N}$, using the normal law distribution (\ref{finallynormal}), in the limit $N\rightarrow\infty$ we have 
 \barray
 \nonumber 
 \prob \[  |C_N(\s)  |    <   d  \, \sqrt{N} \]  & > &
 \prob \[  |C_N(\s)   |    <   d  \,  \sigma_N(\s)  \]   \,=\,  
 \inv{\sqrt{2 \pi}} \int_{-d}^d  dx\,  e^{-x^2/2}  
 \\ 
 \label{PrBN} 
 &=& 1 -  \frac{e^{- d^2 /2}} {\sqrt{ 2 \pi } } \(  \frac{2}{d}  +  O\(  \inv{d^2 }  \)  \).
 \earray
  Chose $d  = \kappa N^\epsilon$ for any $\kappa >0$.   
 Then for any $\epsilon > 0 $,   
 \beq
 \label{ProbOne}
 \lim_{N \to \infty}  \prob \[  C_N(\s)  = O(N^{1/2+ \epsilon}) \]  = 1    \,\,\,.
 \eeq

It is important to stress that our result is actually stronger than this probabilistic argument.  As in all probabilistic arguments,  one is concerned about  \textquotedblleft rare events\textquotedblright  \,of measure zero.   For example, in flipping of a coin, a sequence of $10^{100}$ heads has probability equal to essentially zero,  but it is still possible!  To eliminate in our case the possibility that rare events spoil the asymptotic behavior of the series, let us use a reductio ad absurdum argument, namely let's assume that, in view of some rare events\footnote{For our series this would be an infinitely long series of the same residues for successive primes.}, the series $C_N$ for $N \rightarrow \infty$ rather than going as in eq.\,(\ref{fffff}) would instead behave  as $C_N \simeq N^{\alpha}$ (up to logarithmic corrections) with $\alpha \neq \half$. It this were true, such a behavior of the series $C_N$ should hold for any neighborhood of infinity, namely for all $N$ which satisfies $N > N^*$, for any arbitrarily large $N^*$. But, in turn, this fact would imply that the variance of $C_N$ should always go as $N^{\alpha}$ for {\em all} the infinitely many ensembles $\CE$ of the inertial sequences $A_{N_1,N_2}$ with $N_1 > N^*$. The infinite occurrence of such behavior would contradict firstly the notion itself of \textquotedblleft rare events\textquotedblright \, and, secondly, it would be in clear contrast with the explicit expression (\ref{mostimportantformula}) of the variance computed on all the infinitely many ensembles $\CE$ of the inertial intervals. In other words, we cannot exclude that {\em some} $C_N$ for some specific starting point $\tilde \s$ of the series (\ref{importantequivalence}) and even for long values of $N$ may grow as $N^\alpha$ with $\alpha \geq  \half$, but, if this is the case, using the equivalence of the various series related to $C_N$, we can always change at our will $\tilde \s$ and also take larger and larger values of $N$. Choosing the new $\tilde \s$ to be inside any of the inertial intervals, a behavior of the block variable as $N^\alpha$ would then disappear in favour of the only stable behavior of the series $C_N$ under any possible translation of the inertial intervals and any possible ensemble $\CE$ set up in these intervals, and this is precisely the scaling law of the random walk given by $N^{1/2}$ (again up to logarithms).

\subsection{For the curious reader}

\label{LinearCN}

\def\Ctilde{{\tilde{C}}} 
Eq.\,(\ref{mostimportantformula}) states that, asymptotically, the variance of the series $C_N$ grows linearly in $N$ apart from a factor $\lambda(N,l)$ which takes into 
account finite size corrections in the analysis of the block variables. Such a corrective term depends on the prime $p_{\s}$ of the inertial interval around which we consider the block variables and goes to $1$ when $p_{\s} \rightarrow \infty$ while, at finite $p_{\s}$, it introduces  at most a logarithmic correction in $N$ which -- we know -- is harmless for what concerns 
the implications of Theorem \ref{BNtheorem}. 

\begin{figure}[t]
\centering
\includegraphics[width=0.6\textwidth]{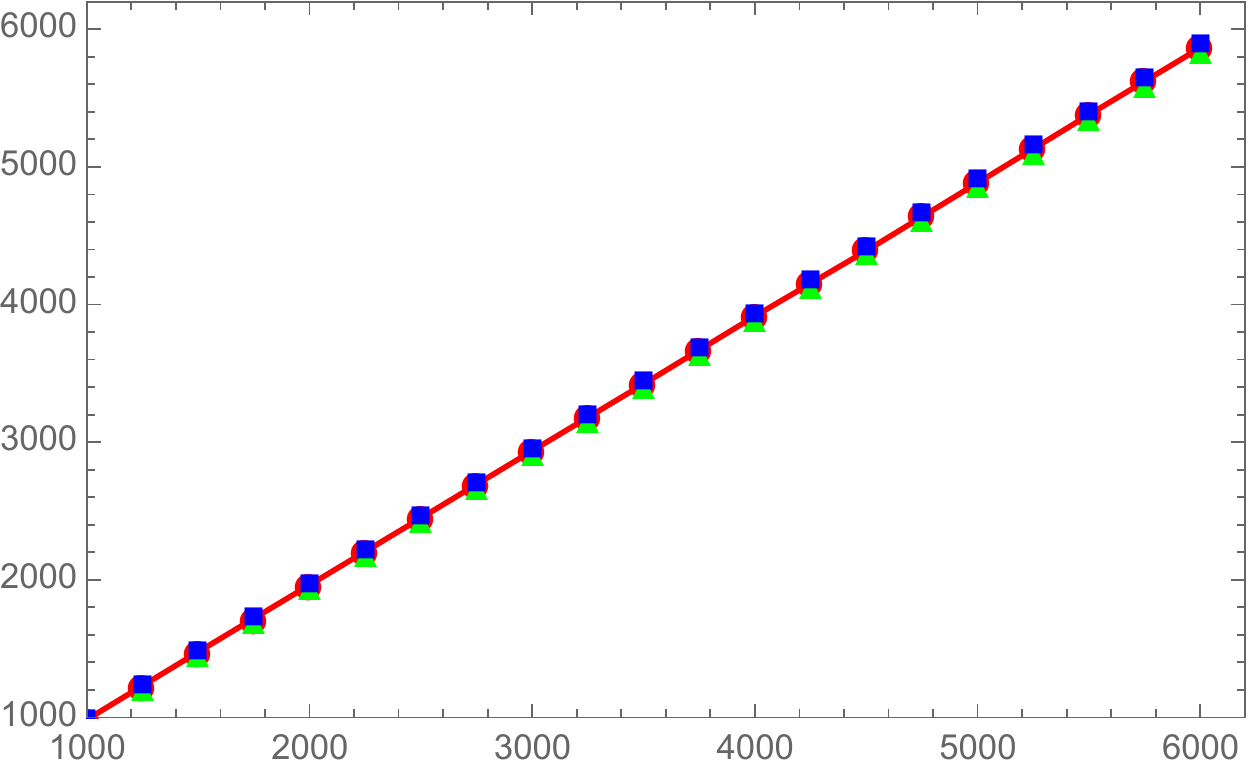}
\caption{$\Ex \[(\Ctilde_N)^2\]$ versus $N$ in the range $(1000,6000)$ in steps of 250 relative to the three non-principal characters with modulus $q=5$, with the legend 
$\chi_2 \rightarrow {\red \bullet}$; $\chi_3 \rightarrow \begingroup
    \color{green}
    \blacktriangle
    \endgroup
$ and $\chi_4 \rightarrow  \begingroup
\color{blue} \blacksquare\endgroup$. The red line is the theoretical prediction with the slope equal to 1. For all characters 
$\Ex \[(\Ctilde_N)^2\]$ grows linearly with $N$ with a slope close to $1$ with few percent of approximation (see Table \ref{tablechq=5C_N}). 
} 
\label{finalfitq55}
\end{figure}

\begin{figure}[t]
\centering
\includegraphics[width=0.6\textwidth]{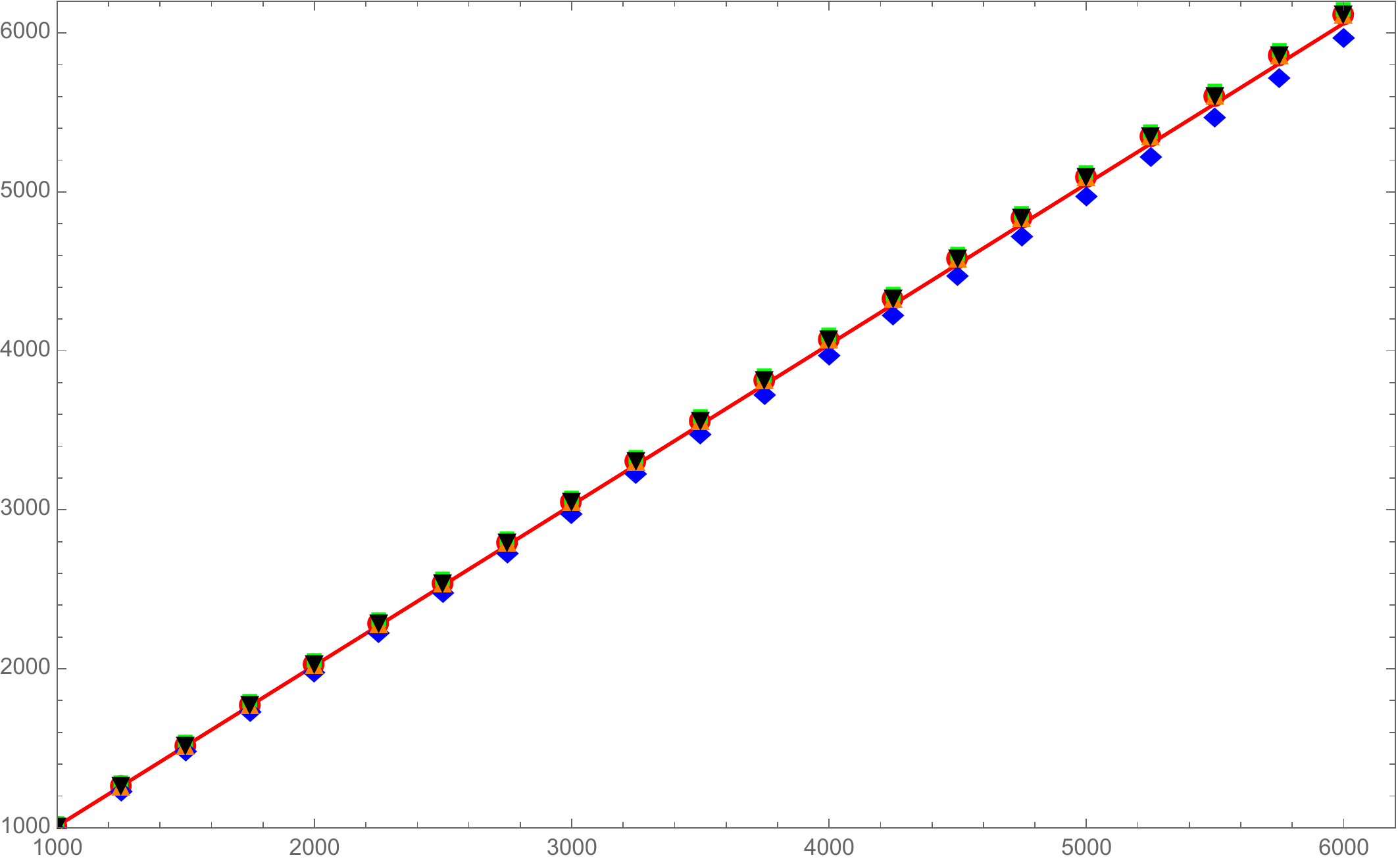}
\caption{$\Ex \[(\Ctilde_N)^2\]$ versus $N$ in the range $(1000,6000)$ in steps of 250 relative to the five non-principal characters with modulus $q=7$, with the legend 
$\chi_2 \rightarrow {\red \bullet}$; $\chi_3 \rightarrow \begingroup
    \color{green}
    \blacktriangle
    \endgroup
$; $\chi_4 \rightarrow  \begingroup
\color{blue} \diamond 
\endgroup$; $\chi_5 \rightarrow \begingroup 
\color{orange} \bigtriangleup \endgroup$ and 
$\chi_6 \rightarrow \begingroup 
\color{black} \bigtriangledown\endgroup$. The red line is the theoretical prediction with the slope equal to 1. For all characters 
$\Ex \[(\Ctilde_N)^2\]$ grows linearly with $N$ with a slope close to $1$ with few percent of approximation (see  
Table \ref{tablechq=7C_N}). 
} 
\label{finalfitq77} 
\end{figure}

One may be curious about the robustness of the expression (\ref{mostimportantformula}) for capturing finite size effects also for intervals $A_{N_1,N_2}$ which are large enough but nevertheless finite: to be definite, let's say the sets of angles $\theta_{p_n}$ relative to the first $10^7$ primes - sets which are pretty simple to generate on a laptop without the need to use a dedicated computer to number theory. As shown by some examples below, it is quite remarkable that the variance $\sigma_N(\s)$ closely follows the behavior predicted by eq.\,(\ref{mostimportantformula}) already for these sets. The examples presented here involve all non-principal characters relative to two different modulus, $q=5$ and $q=7$, and 
are obtained according to the following protocol: 

\begin{enumerate}
\item We have chosen $N_1= 10^5$ and $N_2 =10^7$, therefore with $L = N_2- N_1=9.9 \times 10^6$.
\item We have chosen the length $N$ of the block variables varying in the range $(1000,6000)$ and separated either by a constant interval of length $D=800$ or by interval of random length around this range of value (there was always not much difference between the two cases).
\item At a given $N$, we have computed $C_N(\s)$ (relative to the block variables located at position $\s$) and we have divided this quantity by  $\sqrt{b^2 \lambda(N,\s)}$, defining in this way the new variable 
\beq
\Ctilde_N (\s)  \equiv \frac{C_N(\s)}{\sqrt{b^2 \lambda(N,\s)}} \,\,\,.
\label{defX}
\eeq 
Notice that in this new normalization of the variables $C_N(\s)$ we have not taken into account the term $\rho(N,\s)$ given in eq.\,(\ref{interceptrho}) for the reason that this term is a sub-leading correction in $N$ of the variance and therefore it is expected that should not significantly affect the data. This implies, however, that the variance of $\Ctilde_N$ may have a slope in $N$ slightly different from $1$. 
\item For any given $N$, the cardinality of the set ${\cal S}_M$ is $M=L/(N+D)$, a value that in our examples has been always larger than $10^3$, i.e. large enough to have a reasonable sampling of the quantities $C_N$.  
\item 
In all cases examined, relative to various characters of different modulus $q$, we have always observed a  {\em linear} plot of $\Ex \[(\Ctilde_N)^2 \] $ versus $N$ with a slope remarkably equal to $1$ within very few percent  of approximation. The plot relative to the behavior of $\Ex \[(\Ctilde_N)^2 \]$ versus $N$ for all non-principal characters of $q=5$ 
is in Figure \ref{finalfitq55}, with the relative data given in Table \ref{tablechq=5C_N}. The plot relative to the behavior of $\Ex \[(\Ctilde_N)^2 \]$ versus $N$ for all non-principal characters of $q=7$ is instead in Figure \ref{finalfitq77}, with the relative data given in Table \ref{tablechq=7C_N}.  
\end{enumerate}
These results give evidence of the robust nature of the formula (\ref{mostimportantformula}) which indeed seems to be able to capture efficiently finite size effects of the block variables already for samples of the order of the first $10^7$ primes and it is expected to get better and better going up in the number of primes considered.

.

\newpage

\vspace{5cm}
.
\vspace{5cm}

\begin{table}[bh]
\hspace{-55mm}
\begin{minipage}{.6\textwidth}
\centering
\begin{tabular}{ c |c | c | c | c | c | c | c | c | c | c | c | c | c |}
\cline{2-14}
\cline{2-14}
 & $N$ & 1000 & 1500 & 2000 & 2500 & 3000 & 3500 & 4000 & 4500 & 5000 & 5500 & 6000 & Fit\\ [0.5ex] % inserts table %heading
\cline{2-14}\cline{2-14}
$\chi_2 \rightarrow$ & $\Ex[\Ctilde_N^2]$  &
983.3 & 1474.2 & 1964.2 & 2461.5 & 2950.6 & 3425.7 & 3924.4 & 
4395.6 & 4902.2 & 5374.8 & 5870.8 & $\Ctilde_N^2 = 0.98 N $\\
\cline{2-14}\cline{2-14}
$\chi_3\rightarrow$ & $\Ex[\Ctilde_N]^2$ & 
985.7 & 1476.8 & 1949.4 & 2432.4 & 2917.6 & 3410.8 & 3902.1 & 
4391.7 & 4855.9 & 5361.6 & 5840.4 & $\Ctilde_N^2= 0.97 N $\\
\cline{2-14}\cline{2-14}
$\chi_4\rightarrow$ &$\Ex[\Ctilde_N]^2$ &983.7 & 1475.4 & 1973.1 & 2451.2 & 2943.2 & 3429.7 & 3920.2 & 
4383.3 & 4894.0 & 5399.4 & 5883.8  & $\Ctilde_N^2 = 0.98 N $\\
\cline{2-14}\cline{2-14}
\end{tabular}
\end{minipage}
\caption{$\Ex[(\tilde C_N)^2]$ versus $N$ (here reported in steps of 500) relative to the non-principal characters with modulus $q=5$ (see Table I for their definition) in the range $(1000,6000)$. The slightly different values of $\Ex[(\tilde C_N)^2]$ for the characters $\chi_2$ and $\chi_4$ which are complex conjugate one to other is for the random generation of the block variables. 
In the last column of the table the best fit of the linear growth of $(\tilde C_N)^2$ vs $N$ for each character: in all cases, the slope determined by the best fit differs from the asymptotic value $1$ at most by  $3 \%$. }
\label{tablechq=5C_N}
\end{table}

\newpage
\newpage

\vspace{5cm}
.
\vspace{5cm}

\begin{table}[bh]
\hspace{-55mm}
\begin{minipage}{.6\textwidth}
\centering
\begin{tabular}{ c |c | c | c | c | c | c | c | c | c | c | c | c | c |}
\cline{2-14}
\cline{2-14}
 & $N$ & 1000 & 1500 & 2000 & 2500 & 3000 & 3500 & 4000 & 4500 & 5000 & 5500 & 6000 & Fit\\ [0.5ex] % inserts table %heading
\cline{2-14}\cline{2-14}
$\chi_2 \rightarrow$ & $\Ex[\Ctilde_N^2]$  &
1020.6 & 1533.0 & 2042.8 & 2559.1 & 3070.3 & 3576.0 & 4108.2 &
4598.7 & 5125.0 & 5604.8 & 6113.2 &$\Ctilde_N^2 = 1.02 N $\\
\cline{2-14}\cline{2-14}
$\chi_3\rightarrow$ & $\Ex[\Ctilde_N^2]$ & 
1016.6 & 1527.2 & 2032.6 & 2548.8 & 3082.3 & 3584.1 & 4100.8 & 
4617.3 & 5100.8 & 5626.9 & 6153.7 & $\Ctilde_N^2 = 1.02 N $\\
\cline{2-14}\cline{2-14}
$\chi_4\rightarrow$ &$\Ex[\Ctilde_N^2]$ &
1010.7 & 1507.5 & 2005.4 & 2493.4 & 2985.5 & 3483.3 & 3987.7 & 
4485.8 & 4999.7 & 5483.7 & 5981.59   & $\Ctilde_N^2 = 0.99 N $\\
\cline{2-14}\cline{2-14}
$\chi_5\rightarrow$ &$\Ex[\Ctilde_N^2]$ &
1019.1 & 1527.7 & 2031.7 & 2549.1 & 3073.1 & 3587.7 & 4099.4 &
4584.8 & 5122.8 & 5625.3 & 6135.0 & $\Ctilde_N^2 = 1.02 N $\\
\cline{2-14}\cline{2-14}
$\chi_6\rightarrow$ &$\Ex[\Ctilde_N^2]$ &
1021.1 & 1530.2 & 2043.1 & 2558.3 & 3071.2 & 3575.7 & 4109.5 &
4599.4 & 5123.2 & 5602.2 & 6115.8  & $\Ctilde_N^2 = 1.02 N $\\
\cline{2-14}\cline{2-14}
\end{tabular}
\end{minipage}
\caption{$\Ex[(\tilde C_N)^2]$ versus $N$ (here reported in steps of 500) relative to the non-principal characters with modulus $q=5$ (see Table I for their definition )in the range $(1000,6000)$. The slightly different values of $\Ex[(\tilde C_N)^2]$ for the characters $\chi_2$ and $\chi_6$ which are complex conjugate one to other is for the random generation of the block variables. In the last column of the table the best fit of the linear growth of $(\tilde C_N)^2 $ vs $N$ for each character: in all cases, the slope determined by the best fit differs from the asymptotic value $1$ by at most  $2\%$.}
\label{tablechq=7C_N}
\end{table}

\newpage

\begin{figure}[b]
\centering\includegraphics[width=.6\textwidth]{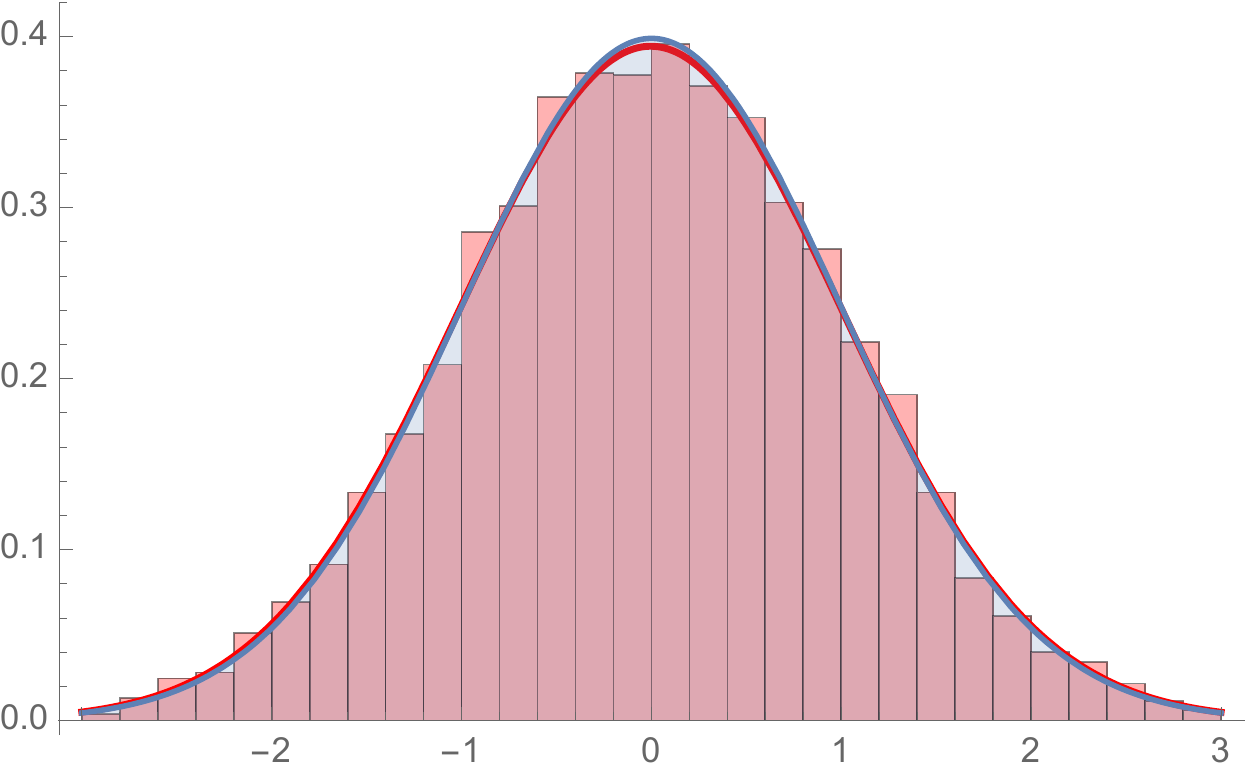}
\caption{Numerical evidence for the normal distribution proposed in \eqref{finallynormal}.  What is shown is a histogram of the LHS of  \eqref{finallynormal} which are properly normalized block variables $C_N  (\s)$ for the character $\chi_2$ mod $7$ in Table \ref{tablech}.  The ensemble $\CE$ corresponds to $N=6000$,  $D=100$,  with 
$M=10.000$ states. The red curve is the fit to the data, which is the normal distribution $\CN(-0.004, 1.01)$.  The nearly indistinguishable  blue curve is the prediction $\CN(0,1)$.}
\label{FigureCLT1}
\end{figure}

For a given character and a given $N$, we have also made the histogram of the $M$ values relative to the ensemble $\CE$ 
of the quantities $\Ctilde_N(\s)$, normalized to their variance  $\Ex \[(\Ctilde_N)^2\]$. 
As expected, the distribution of the $M$ values of $\Ctilde_N (\s) $, for $N$ large enough, is gaussian distributed with a very high level of confidence  (see Figure \ref{FigureCLT1} for 
one of such examples).

\section{Conclusions}\label{conclusions}

In this paper we have addressed the Generalized Riemann Hypothesis for the Dirichlet $L$-functions of non-principal characters  based on studying an enlarged region of convergence of their 
Euler infinite product representation. We have shown that the convergence of the Euler product is controlled by the large $N$ behavior of the series $C_N$ defined in eq.\,(\ref{defCNN}): a purely diffusive random walk behavior as $N^{1/2 +\epsilon}$ of this series, for arbitrarily small $\epsilon >0$, signifies that all zeros of these functions are along the critical line $\Re(s) = \half$. We have established this result through a series of steps which have enlightened various aspects of the problem. 

First we have considered a random set of $L$-functions which have the virtue of having exactly the same zeros as the original $L$-function and we have established a normal law for the
analog of the series $C_N$ in Theorem \ref{BpCLTDir}.  However the implications for the GRH were inconclusive because such a normal law rules the fluctuation of the series 
$C_N$ but with respect to their mean $m_N$, and we showed that estimating the behavior of $m_N$ was tantamount to proving the validity itself of the GRH. 

However we have subsequently shown that there is a natural explanation of such a diffusive behavior of the series $C_N$, which can be established using the Dirichlet theorem on the equidistribution of reduced residue classes modulo $q$ and the Lemke Oliver-Soundararajan conjecture on the distribution of pairs of residues on consecutive primes. As a matter of fact, the series $C_N$ are amenable of a probabilistic approach albeit they are deterministic quantities: from this point of view, they share several properties with random series encountered in other scientific fields. As for other random series, however, in order to control the growth of the series $C_N$ by varying $N$ one has to face the so-called {\em Single Brownian Trajectory Problem}. Such a problem in our case can be solved by defining an ensemble $\CE$ which involves block variables of the original series $C_N$: these block variables provide 
\textquotedblleft stroboscopic\textquotedblright \,snapshots of the original series and realize its sampling. The mean of this variance vanishes  by  virtue of the Dirichlet theorem while the variance, as shown in eq.\,(\ref{mostimportantformula}), goes linearly in $N$ (up to logarithmic corrections). This leads to the normal distribution \eqref{finallynormal} for  the series $C_N$.  As discussed in the text, there are some correlations among consecutive angles $\theta_{p_n}$  which however do {\em not} spoil the diffusive behavior of the series $C_N$.

In summary, based on Theorem \ref{GM1}, which assumes the LOS conjectures, we can establish a purely diffusive random walk behavior of the series $C_N$,  and this implies 
that all zeros of the Dirichlet $L$-functions of non-principal characters are along the same critical line since the Euler product converges for $\Re (s) > \half$.  A natural question is how strongly this conclusion depends on the LOS conjectures? We would answer that the most important property of the pair correlation is its asymptotic uncorrelated 
behavior given in eq.\,(\ref{largeasymptotic}), while the details of the LOS formula are essential only for controlling finite $\s$ effects.

\vspace{5mm}
 
\section*{Acknowledgments}

AL would like to thank Steve Gonek  for discussions and for pointing out the reference \cite{GrosswaldSchnitzer}.  GM would like to thank Don Zagier, Karma Dajani, Giorgio Parisi, Gianni Dal Maso, Andrea Gambassi and Satya Majumdar for interesting discussions and Robert Lemke Oliver for useful email correspondence.  AL would like to thank SISSA in Trieste, Italy, where this work was begun while GM would like to thank the Institute for Theoretical Physics in Utrecht where the great majority of this work was done. GM would like also to thank the Simons Center in Stony Brook and the International Institute of Physics in Natal for the warm hospitality and support during the initial and final parts of this work respectively.

\newpage

\appendix

\section{$L$-functions as grand canonical partition functions of non-interacting particles}\label{AA}

In this section we show that the Dirichlet $L$-functions can be interpreted as generalized grand canonical partition functions of an infinite set of non-interacting bosonic particles. The idea is not new (see \cite{Julia,Spector}) but it is worth recalling  it. 

Consider an Hilbert space whose Fock basis is given by an infinite countable set of  bosonic  creation operators $a^\dagger_1, a^\dagger_2, a^\dagger_3, \ldots $ associated to the increasing sequence $p_1, p_2, p_3, \ldots$ of prime numbers, with the energy of each mode given by 
\beq
\epsilon_n \,=\,\log p_n \,\,\,,
\label{energymode}
\eeq
A generic state of such an Hilbert space can be expressed as 
\beq
| M \rangle \,=\, \left(\prod_{k=1}^r (a_{i_k})^{\sigma_k} \right) \,| 0 \,\rangle 
\,\,\,, 
\label{multiparticle}
\eeq
where the integer $M$ is given by 
\beq
M \,=\, p_{i_1}^{\sigma_1} \, p_{i_2}^{\sigma_2} \cdots p_{i_k}^{\sigma_k} \,\,\,. 
\eeq
From  the unique factorization of the integers in terms of the primes, $M$ is uniquely specified in terms of the bosonic creation operators and clearly their order does not matter. With this notation, the primes give rise to one-particle states while composite numbers are given in terms of multi-particle states. Assuming no interaction among the modes, the energy of this states is equal to the sum of the energies of its constituents 
\beq
E_M \,=\,\sum_{k=1}^r \sigma_k \,\epsilon_{i_k} \, =\,   \sum_{k=1}^r \sigma_k \,\log p_{i_k}\, 
\,=\, \log M \,\,\,.
\label{total energy}
\eeq
We call such a non-interacting system the  {\em prime number gas}. We would like now to compute the generalized grand canonical partition functions of this system by eventually filtering some of its states: for instance, taken a set ${\mathcal A}$ of $k$ primes $\left\{p_{a_1}, p_{a_2}, \ldots, p_{a_k}\right\}$ we can decide to keep all the states (\ref{multiparticle}) in which there never appears any of the corresponding creation operators. At the same time, defining the integer number $q = p_{a_1} \cdot p_{a_2} \cdots p_{a_k}$ given by the product of all the primes in the set ${\mathcal A}$, we can assign to the remaining states $| \tilde M \rangle$ with non-zero residue mod $q$ 
a complex weight $\chi(\tilde M) \equiv e^{i \theta_{\tilde M}}$ of unit modulus  according to the rules of the characters mod $q$ already recalled in Section II, which we repeat here for convenience:  
\begin{enumerate}
\item 
$\chi(m+q) \,=\, \chi(m) $. 
\item 
$\chi(1) =1 $ and $\chi(0) = 0$. 
\item 
$\chi( m n ) \,=\, \chi(m) \, \chi(n)$. 
\item 
$\chi(m) = 0$ if $(m,q) > 1$ and $\chi(m) \neq 0$ if $(m,q) =1$. 
\item 
If $(m,q) =1$ then $(\chi(m))^{\varphi(q)} =1$, namely $\chi(m)$ have to be $\varphi(q)$-roots of unity. 
\end{enumerate}
Let's call $\theta_{\tilde M}$ the {\em abelian charge} assigned to the state $ | \tilde M \rangle$:  according to the Rule 3, it is given by the sum of the charges $\theta_{p_i}$ assigned to the prime factors $p_i$ present in $\tilde M$, weighted with their relative multiplicities 
\beq
\tilde M  \,=\,p_{i_1}^{\sigma_1} \, p_{i_2}^{\sigma_2} \cdots p_{i_h}^{\sigma_h} 
\,\,\,\,\,\,\,
\longrightarrow 
\,\,\,\,\,\,\,
\theta_{\tilde M} \,=\, \sum_{i=1}^h \sigma_i \theta_{p_i}
\eeq
From the general theory of the characters, we know that there $\varphi(q)$ consistent ways of assigning the charges to the states $| \tilde M \rangle$, keeping into account the period $q$ stated by the Rule 1. States $| M \rangle$ which have zero residue mod $q$ have a weight $\chi(M) = 0$. 

Let's now consider the generalized grand canonical partition function associated to a given set of weights $\chi(M)$ for the states $| M \rangle$, according to the Rule 4 above, and at the inverse temperature $s = 1/T$, 
\beq
\Omega_{\chi}(s) \,=\,
\sum_{M=1}^{\infty} \chi(M) \, e^{-s E_M} \,=\,
\sum_{M=1}^{\infty} \frac{\chi(M)}{M^s} \,\,\,.
\eeq
However, from  the non-interactive nature of the system and the additivity of the charges 
assigned to the states, the grand canonical partition function is just given by the infinite product of the constituent partition functions 
\beq
\Omega_{\chi}(s) \,=\,\prod_{n=1}^\infty \left(
\sum_{\sigma_n=0}^\infty \left(\chi(p_n) \, e^{-s \epsilon_n}\right)^{\sigma_n} 
\right)
\,=\,\prod_{n=1}^\infty \frac{1}{1- \chi(p_n) e^{-s \epsilon_n}} 
\,=\,\prod_{n=1}^\infty \frac{1}{1- \frac{\chi(p_n)}{p_n^s}} 
\eeq
Clearly $\Omega_{\chi}(s) = L(s,\chi)$ and the Euler identity  
\beq
\sum_{n=1}^{\infty} \frac{\chi(n)}{n^s} \,=\,
\prod_{p} \frac{1}{1- \frac{\chi(p)}{p^s}}
\,\,\,,
\eeq
can be then interpreted as an equivalence of the microcanonical and grand canonical ensemble of the non-interacting prime number gas.   

%\newpage

\section{Poles and Fisher zeros of the generalized partition functions}
Given the periodicity of the character $\chi(m)$, the infinite sum on the states for the generalized 
partition function 
\beq
L(s,\chi) \,=\,\sum_{m=1}^\infty \frac{\chi(m)}{m^s} \,\,\,,
\eeq
can be organized as a finite sum of the $\varphi(q)$ different characters as (see 
eq.\,(\ref{LHurwitz}) in the text) 
\beq
L(s,\chi) \,=\,\frac{1}{q^s} \, \sum_{r=1}^q \chi(r) \, \zeta\left(s, \frac{r}{q}\right) 
\,\,\,, \label{LHurwitz2}
\eeq
where 
\beq
\zeta(s,a) \,=\, \sum_{n=0}^\infty \frac{1}{(n+a)^s} \,\,\,,
\label{Hurwitz22}
\eeq 
Notice that each partition function associated to the Hurwitz function $\zeta(s,a)$ 
is divergent at $s=1$, where there is a pole. From a statistical mechanics point of view, such a singularity, known as {\em Hagedorn temperature}, is due to the exponential divergence of the density of states of the system: indeed, with a spectrum of energy given by $E_n = \log (n+a)$, the density of states is given by 
\beq
\omega(E) \,=\, \frac{dn}{dE} \,=\, e^E \,\,\,,
\eeq
and, making the change of variable $n \rightarrow E$ in the Hurwitz function, we 
have 
\beq
\zeta(s,a) \,\simeq\, \int dE \omega(E) \, e^{-s E} \,=\, 
\int dE \, e^{-(s-1) E} \,=\, \frac{1}{s-1} + \cdots 
\eeq
Hence, at sufficiently  high temperature, the exponentially  decreasing Boltzmann factor $e^{-s E}$ is no longer able to compensate for the exponential growth  of the number of states and therefore the partition function of the system explodes. 

Looking at eq.\,(\ref{LHurwitz2}), the partition function of each charge sector is individually divergent but for all characters but the principal one it holds that 
\beq
\sum_{r=0}^q \chi(r)\,=\, 
\left\{
\begin{array}{cll}
\frac{\varphi(q)}{q} & & {\rm if }\, \, \chi = \chi_1 \\
0 & & {\rm if }\,\, \chi \neq \chi_1 \,\,\,.
\end{array}
\right.
\eeq
and therefore the singularity at $s=1$ cancels. Therefore, for all non principal characters, the 
corresponding $L$-function is an entire function of $s$, fully characterized by its zeros in $s$. In statistical mechanics those zeros are known as {\em Fisher zeros} \cite{Fisher} and they can help in clarifying the nature of the physical system under study. We can focus the 
attention only on the non-trivial zeros of the $L$-function by considering the completed $L$-function $\hat L(s,\chi)$ for primitive characters
\beq
\hat L(s,\chi) \equiv \left(\frac{q}{\pi}\right)^{(s+\delta)/2} \, 
\Gamma\left(\frac{s+\delta}{2}\right) \, L(s,\chi)\,\,.  
\label{completedLfunctionaa}
\eeq
Denoting by $\rho_\chi$ the non-trivial zeros in the critical strip $0 < \sigma < 1$, the completed $L$-function admits the Hadamard infinite product \cite{Steuding} 
\beq
\hat L(s,\chi) \,=\, \exp\left(A_{\chi} + B_{\chi} s\right) \, 
\prod_{\rho_{\chi}} \left(1-\frac{s}{\rho_{\chi}}\right) \, e^{s/\rho_{\chi}} \,\,\,,
\eeq
where the series $\sum_{\rho_{\chi}} |\rho_{\chi}|^{-1}$ diverges while $\sum_{\rho_{\chi}} |\rho_{\chi}|^{-1 - \epsilon}$ converges for any positive $\epsilon$. $A_\chi$ and $B_\chi$ are two constants which depends on the primitive character $\chi$. 

The logarithm of the completed $L$-function has the meaning of the free energy $F(s,\chi)$ of the prime number gas and it admits the high-temperature series expansion 
\beq
F(s,\chi) \,=\, A_\chi + B_\chi s + \sum_{n=2}^\infty \alpha_n s^n \,\,\,,
\label{generalisedFreeEnergy}
\eeq
where the coefficients $\alpha_n$ are nothing else but the $n$-th inverse moment of the zeros  
\beq
\alpha_n \,=\, \sum_{\rho_\chi} \frac{1}{\rho_\chi^n} \,\,\,.
\label{momentzeros}
\eeq
As in the case of the Riemann zeta function, also for the $L$-function it is possible to compute the number $N(T,\chi)$ of non-trivial zeros $\rho_{\chi} = \beta_\chi 
+ i \gamma_\chi$ with $| \gamma_\chi| \leq T$ and the leading behavior of this function 
is given by 
\beq
N(T,\chi) \,=\,\frac{T}{\pi} \,\log\frac{q T}{2 \pi e} + {\mathcal O}(\log q T) 
\,\,\,.
\label{asymptoticzeros}
\eeq
In this expression zeros from the lower half-plane are counted as well for the lack of a
symmetry with respect to the real axis in case of non-real characters. For the real characters, this formula implies an average density of zeros given by 
\beq
d(T) \simeq \frac{1}{2\pi q} \log\frac{q T}{2\pi} \,\,\,.
\label{densityof zeros}\,
\eeq

\section{Proof of the theorems by Grosswald-Schnitzer and Chernoff.}\label{proofstheorems}

In this Appendix we briefly discuss the proof of the two theorems presented in Section III of the text. 

Let's start first with the proof of the Grosswald-Schnitzer theorem. Consider the set of integers $p'_n$ 
which satisfy the conditions given in eq.\,(\ref{ppn}) and, in terms of them, define for 
$\Re(s) > 1$ the function $L'(s,\chi)$ by means of the absolutely convergent infinite product 
\beq
L'(s,\chi) \,=\,\prod_{n}\left(1-\frac{\chi(p'_n)}{(p'_n)^s}\right) \,\,\,.
\eeq 
Let's now set 
\beq
\theta(s) \,=\,\prod_{n} \frac{\left(1 - \chi(p_n)\, p_n^{-s}\right)}{\left(1 - \chi(p'_n)\,p_n'^{-s}\right)} 
\,\,\,, 
\eeq 

The  function $\theta (s)$   can be proven  to converge absolutely for  $\Re(s) > 0$ and with no zeros in this region. Since 
\beq
L'(s,\chi) \,=\, \theta(s) \, L(s,\chi) \,\,\,,
\label{crucialidentity}
\eeq
it is clear that $L'(s,\chi)$ inherits the analytic structure of the Dirichlet $L$-function.   In particular it can be analytically continued into the whole half plane $\Re(s) >0$ and has exactly the same zeros (including multiplicities) as the original $L$-function in the critical strip. 

\vspace{3mm}

Let us now consider Proof of Chernoff's theorem. To this aim let's take 
the infinite product representation of the Riemann $\zeta$-function 
\beq
\zeta(s) \,=\,\prod_{n=1}^\infty  \( 1 -  \frac{1}{p_n^s} \)^{-1} \,\,\,,
\eeq
and  its logarithm
\begin{eqnarray}
\log \zeta(s) &\,=\, &- \sum_{n=1}^\infty \log\left(1 - \frac{1}{p_n}\right)
\,=\,\sum_{n=1}^\infty \sum_{k=1}^\infty \frac{p_n^{- s k}}{k} \\
& = & \sum_{k=1}^\infty \frac{1}{k} \sum_{n=1}^\infty p_n^{-s k} 
\nonumber \,\,\,.
\end{eqnarray}
It is easy to see that the divergence of this expression is controlled just by the first term, namely by the series\footnote{The function $\eta (s)$ is commonly referred to as the
``prime zeta-function".}
\beq
\eta(s) \,=\,\sum_{n=1}^\infty \frac{1}{p_n^s} \,\,\,.
\eeq
If we now substitute into  this expression the continuous approximation of the $n$-th prime based on the prime number theorem, 
i.e. $p_n \simeq n \,\log n$, we end up with the series 
\beq
\eta'(s) \,=\,\sum_{n=2}^\infty \frac{1}{(n \,\log n)^s}
\,\,\,. 
\eeq  
Writing it as Stieltjes integral 
\beq
\eta'(s) \,=\,\int_2^\infty (x \,\log x)^{-s} \, d[x] \,\,\,
\eeq
and proceeding first through an integration by parts  and some  additional steps, 
one ends up with the final expression \cite{Chernoff}
\beq
\eta'(x) \,=\, I(s)  - \int_2^\infty \frac{\{x\} \, (\log x +1)}{(x \,\log x)^{s+1}} \,\,\,,
\label{psi}
\eeq
where $\{x\} = x - [x]$ is the fractional part of $x$ and 
\beq
I(s) \,=\, \frac{2 (2 \log 2)^{-s}}{s} + \frac{1}{s} \,\left[
\frac{(2 \log 2)^{-s-1}}{s-1} + (s-1)^{-(s-1)} \, \Gamma(1-s) +
 \int_0^{\log 2} u^{-(s-1)} \,e^{-(s-1) u}\, du \right]
\,\,\,.
\eeq
Apart from the explicit singularities at $s=0$ and $s=1$, the quantity $I(s)$ has an analytic continuation into the physical strip and has no singularities there. Moreover, the last term in eq.\,(\ref{psi}) defines an analytic function for $\Re(s) > 0$. Hence, altogether the function $\eta's)$ has an analytic continuation into the critical strip $0 < \Re(s) < 1$ in which it has no singularities. Hence the infinite product 
\beq
\zeta'(s) \,=\,\prod_{n=2}^\infty \left(1 - (n \,\log n)^{-s}\right)^{-1} \,\,\,,
\eeq
has an analytic continuation into the physical strip and has no zeros  there. 

\section{Kac's central limit theorem}\label{Kactheor}

One of the remarkable results of Mark Kac regards the behavior of deterministic trigonometric series with linearly independent frequencies \cite{Kac}, as the one 
given in eq.\,(\ref{noshiftfunction}) in the text. Let's recall that real numbers $\lambda_1, \lambda_2, \ldots, \lambda_n$ are called independent on the field of the rationals if the only solution $(k_1,k_2,\ldots, k_n)$ of the equation 
\beq 
k_1 \,\lambda_1 + k_2 \, \lambda_2 + \cdots + k_n\, \lambda_n \,=\,0 
\eeq
is 
\beq 
k_1 = k_2 = k_3 = \cdots = 0 \,\,\,.
\eeq 
Notice that the sequence $\lambda_n = \log p_n$ consists indeed of linearly independent numbers. Let's now state the Kac's theorem. 

\begin{figure}[t]
\centering
\includegraphics[width=0.40\textwidth]{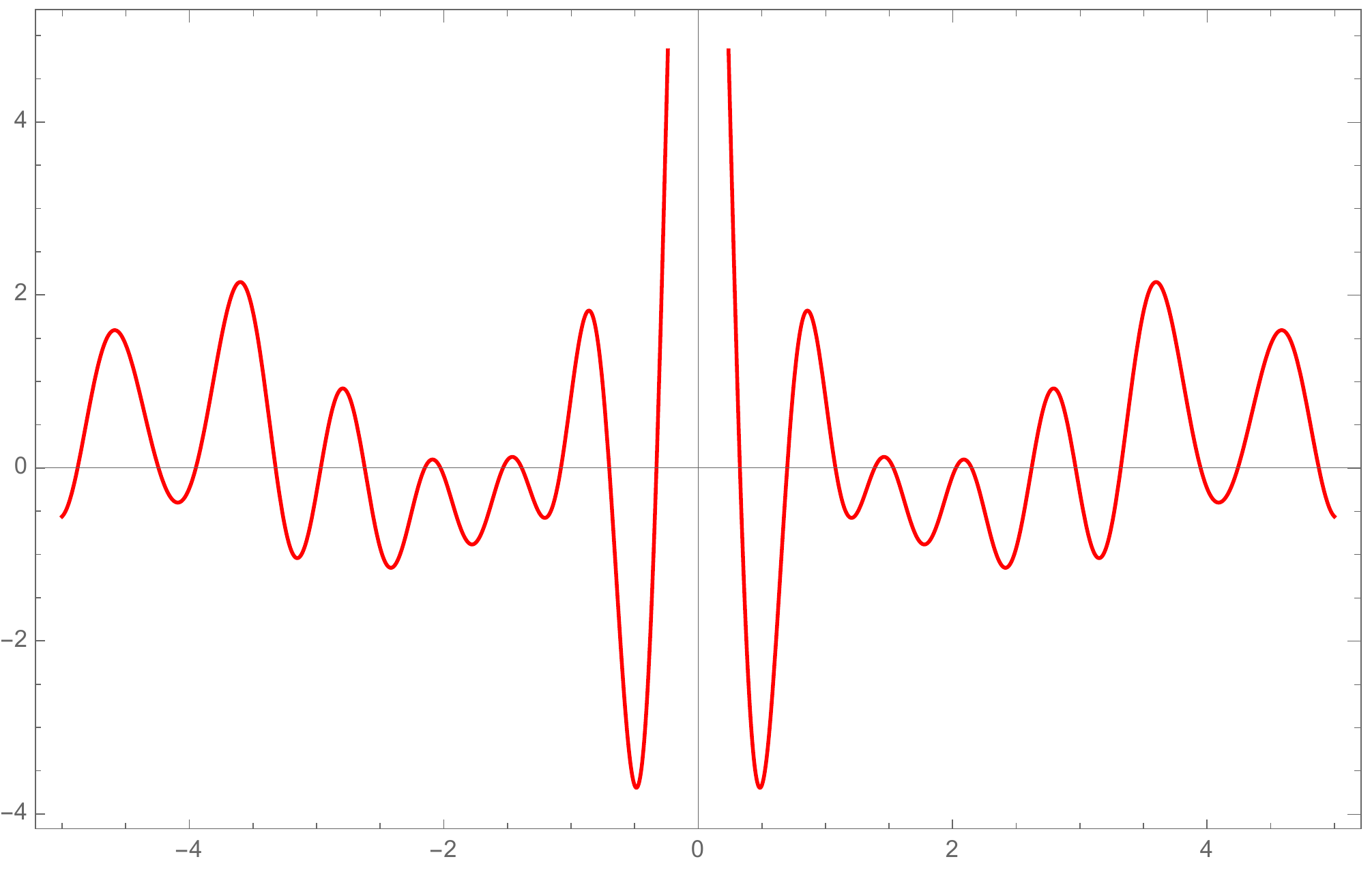}
\,\,\,
\includegraphics[width=0.40\textwidth]{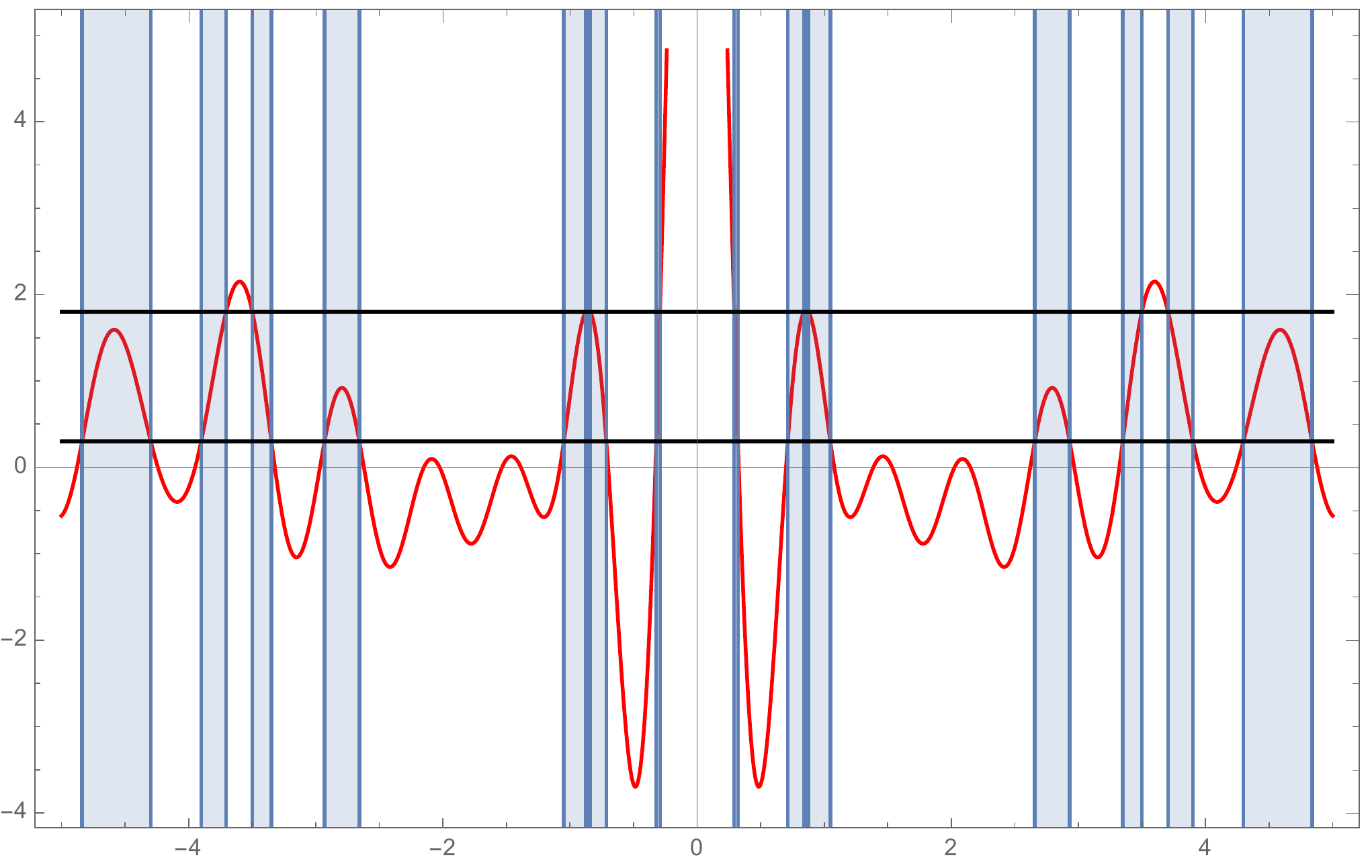}
\caption{Content of the Kac's theorem. Left hand side: plot of the function $F_N(t)$ on an interval $T$. Right-hand side: 
intervals (coloured in blue) where the function $F_N(t)$ is between two values $a$ and $b$. Kac's theorem states that the sum of the coloured intervals divided by $T$ is a well-defined quantity in the limit $T\rightarrow \infty$ which, when $N \rightarrow \infty$, is given in terms of the normal distribution. } 
\label{Kacgraphical}
\end{figure}

\begin{figure}[t]
\centering
\includegraphics[width=0.40\textwidth]{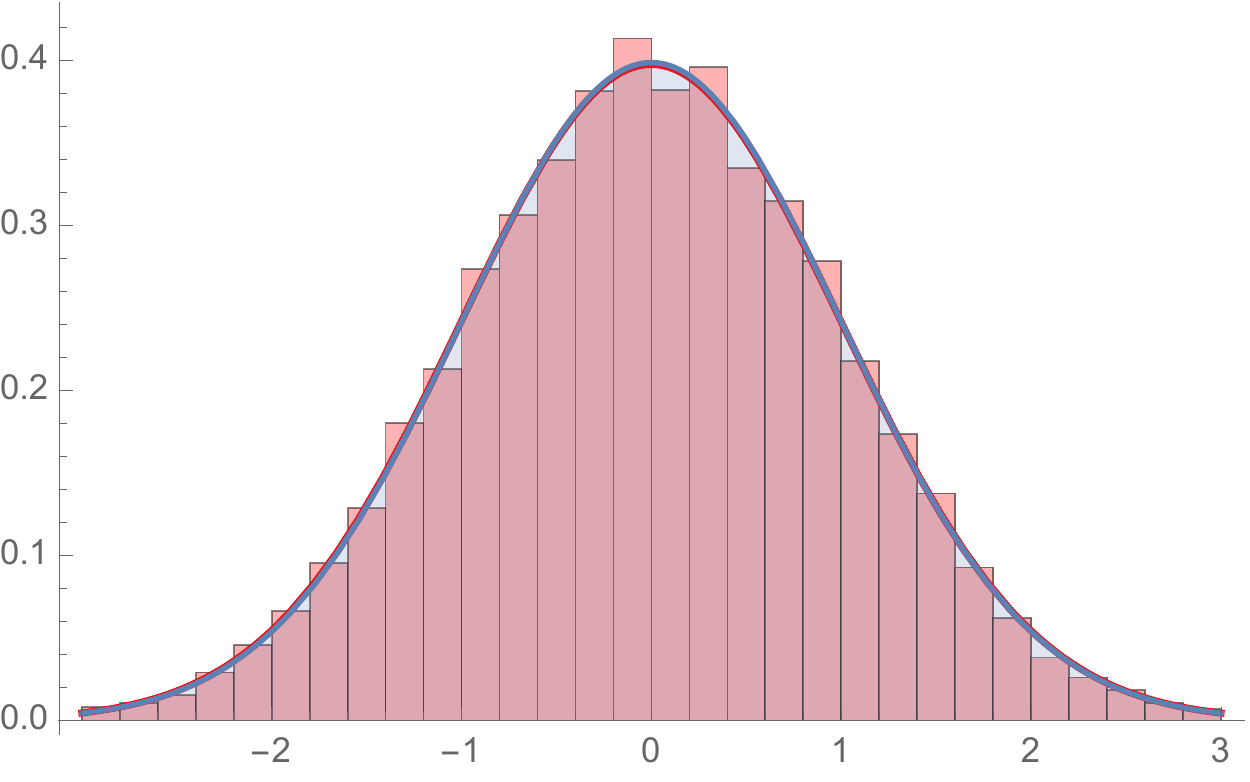}
\,\,\,\hspace{5mm}
\includegraphics[width=0.30\textwidth]{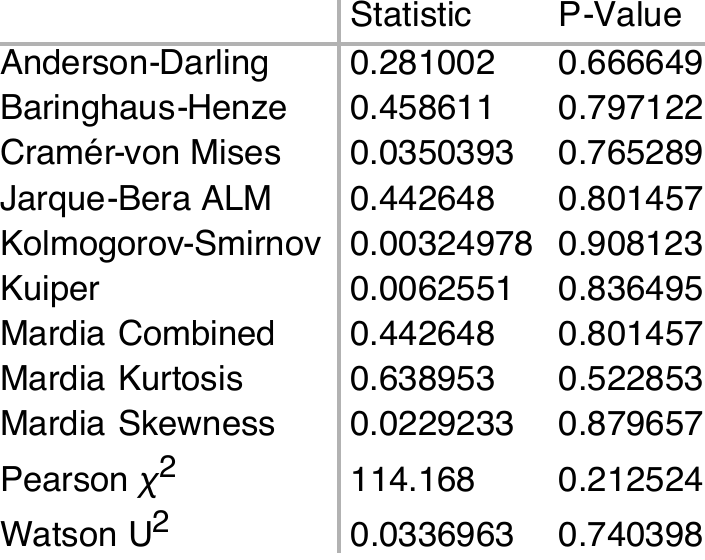}
\caption{Numerical analysis of the Kac's theorem. In this example $T =  10^5$ while $N=10^3$. For the frequencies $\lambda_i$ we have chosen an increasing sequence of incommensurate irrational numbers such that $\Delta\lambda_i = \lambda_i - \lambda_{i-1} 
= {\mathcal O}(1)$. We have generated $M=2\times 10^4$ random points $t_i \in (-T,T)$ and made an histogram of the corresponding values $F_N(t_i)$. The result is the normal distribution on the left hand side, as confirmed with a high level of confidence of the various indicators of the fit shown on the right hand side. 
} 
\label{normalllKac}
\end{figure}

\vspace{3mm}
\begin{theorem}(Kac) 
{\it Let $\lambda_n$, $n=1,2,\ldots, N$  be a sequence of linearly independent numbers on the field of the rational and consider the function $F_N(t)$ defined as
\beq 
F_N(t) \,=\, \sqrt{2} \,\frac{\cos\lambda_1t +  \cos\lambda_2t + \cdots + \cos\lambda_Nt}{\sqrt{N}} \,\,\,.
\label{FNNNN}
\eeq
Let $\mu_R[S]$ denote the Lebesgue measure of a set $S$ on $\mathbb{R}$. Then 
\beq
\label{Kaccc}
\lim_{N\rightarrow \infty} \lim_{T\rightarrow \infty} \frac{1}{2T} \, \mu\, \left[-T \leq t \leq T \, : a \leq F_N(t) \leq b \right] = 
 \frac{1}{\sqrt{2\pi}} \int_a^b \exp\left[-\frac{x^2}{2} \right] dx \,\,\,.
 \eeq 
 }
\end{theorem}
The content of the Kac's theorem is illustrated in Figure \ref{Kacgraphical}. 
\noindent 
It is crucial to stress that it is important to take the time average and the limit $T\rightarrow \infty$ {\em before} the limit $N \rightarrow \infty$ and that these two limits do not commute. In other words, the Kac's theorem concerns with the infinite time averages of the family of the partial sums $F_N$ rather than the infinite time average of the limit function $F(t) = \lim_{N\rightarrow \infty} F_N(t)$. In order to 
appreciate this point, let's present the main points of the proof of this theorem \cite{Kac}. Define 
\beq
g(x) \,=\, \left\{ 
\begin{array}{lll}
1 &, & a < x < b \\
0 &, & {\rm otherwise}
\end{array}
\right.
\eeq
which, in terms of Fourier transform, can be expressed as
\beq 
g(x) \,=\, \frac{1}{2\pi} \int_{-\infty}^{+\infty} G(\xi) e^{-i x \xi} \, d\xi
 \,\,\,, 
\eeq
where 
\beq 
G(\xi) \,=\, 
\frac{e^{i b \xi} - e^{i a \xi}}{i \xi}\,\,\,.
\eeq 
Let's also define the time average of a function $A(t)$ as 
\beq
\langle A \rangle \,=\, \lim_{T\rightarrow \infty} \frac{1}{2 T} \, \int_{-T}^{T} A(t) \, dt\,\,\,.
\label{timeaveragge}
\eeq
Let $S(a,b)$ be the set of points $t$ on the real axis for which 
\beq
a \leq  \sqrt{2} \,\frac{\cos\lambda_1t +  \cos\lambda_2t + \cdots + \cos\lambda_Nt}{\sqrt{N}} \leq b \,\,\,.
\eeq
and let's consider the fraction of time in the interval $(-T,T)$ where the function $F_N(t)$ is between these two values. Such a fraction can be computed in terms of the function $g(x)$ as 
\begin{eqnarray}
&& \frac{1}{2 T} \int_{-T}^T g\left(\sqrt{2} \,\frac{\cos\lambda_1t +  \cos\lambda_2t + \cdots + \cos\lambda_Nt}{\sqrt{N}}\right)\, dt \,=
\\ 
&& = \frac{1}{2\pi} \int_{-\infty}^{\infty} G(\xi) \,\left[\frac{1}{2 T} \int_{-T}^T 
\exp\left(i \xi \sqrt{2} \,\frac{\cos\lambda_1t +  \cos\lambda_2t + \cdots + \cos\lambda_Nt}{\sqrt{N}}\right) \, dt \right] d\xi \nonumber \,\,\,.
\end{eqnarray}
The crucial point now is that, only taking the limit $T \rightarrow \infty$, we have a decoupling of the various terms, namely
\beq
\lim_{T\rightarrow \infty} \frac{1}{2 T} \int_{-T}^T 
\exp\left(i \xi \sqrt{2} \,\frac{\cos\lambda_1t +  \cos\lambda_2t + \cdots + \cos\lambda_Nt}{\sqrt{N}}\right) \, dt \,=\, 
\left[J_0\left(\sqrt{2} \frac{\xi}{\sqrt{N}}\right)\right]^N \,\,\,,
\label{arbN}
\eeq
where $J_0(x)$ is the Bessel function. To see how this happens consider the case $N=2$ for which 
\beq
 \frac{1}{2 T} \int_{-T}^T 
\exp\left(i \xi \,(\cos\lambda_1t +  \cos\lambda_2t) \right) \, dt \,=\,
\sum_{k,l=0}^\infty \frac{(i \xi)^k (i \xi)^l}{k! l!}  \frac{1}{2 T} \,\int_{-T}^T \cos^k \lambda_1 t \, \cos^l \lambda_2 t\, dt\,\,\,. 
\eeq
The integrand can be written as linear combination of exponentials 
\beq
\cos^k \lambda_1 t \, \cos^l \lambda_2 t \,=\, 
\frac{1}{2^k} \,\frac{1}{2^l} \sum_{r=0}^k \sum_{s=0}^l 
\binom{k}{r} \,\binom{l}{s} \, e^{i \left[2 r - k) \lambda_1 + (2 s - l) \lambda_2 \right] t} 
\eeq
Since 
\beq
 \frac{1}{2 T} \int_{-T}^T e^{i a t}
 \,=\,
 \left\{
\begin{array}{cll}
1 & , & \alpha =0 \\
\frac{\sin\alpha T}{T}& , & \alpha \neq 0  
\end{array}
\right.
\eeq
we have
\beq
\langle e^{i \alpha t}\rangle \,=\,
\lim_{T\rightarrow \infty}  \frac{1}{2 T} \int_{-T}^T e^{i a t} \,=\, 
\left\{
\begin{array}{lll}
1 & , & \alpha =0 \\
0 & , & \alpha \neq 0  
\end{array}
\right.
\eeq
Given the linear independence of the frequencies, the only solution of the equation 
$$
(2 r - k) \lambda_1 + (2 s- l) \lambda_2 \,=\, 0 
$$
is given by $k = 2 r$ and $ l = 2 s$ and therefore 
\beq
\langle 
\cos^k \lambda_1 t \, \cos^l \lambda_2 t \rangle 
\,=\, \langle \cos^k \lambda_1 t \rangle \, \langle \cos^l \lambda_2 t \rangle\,\,\,. 
\eeq 
Hence 
\beq
\langle e^{i \xi (\cos\lambda_1 t + \cos\lambda_2 t)} \rangle \,=\, 
\langle e^{i \xi (\cos\lambda_1 t )} \rangle \, \langle e^{i \xi (\cos\lambda_2 t)} \rangle
\eeq
and 
\beq
\langle e^{i \xi (\cos\lambda t )} \rangle \,=\, \frac{1}{2 \pi} \int_0^{2 \pi} e^{i \xi \cos \theta} d\theta \,=\, J_0(\xi) \,\,\,.
\eeq 
The calculation presented for $N=2$ can be generalized to arbitrary $N$ and this leads to the equation (\ref{arbN}). Notice it was crucial to take the limit $T \rightarrow \infty$ to have the factorized expression (\ref{arbN}) in terms of the $N$ terms of the original sum $F_N(t)$.  
Such a factorization expresses the statistical independence of each variable and this property leads directly to a central limit theorem. 
In fact, since 
\beq
\lim_{N\rightarrow \infty} \left[J_0\left(\sqrt{2} \frac{\xi}{\sqrt{N}}\right)\right]^N \,=\, e^{-\xi^2/2} 
\eeq
we arrive to the result (\ref{Kaccc}). It is quite simple to have a numerical confirmation of the normal distribution implied by the Kac's theorem. To this aim, one chooses an interval $T$ sufficiently large and a large integer $N$. Then one generates an uniform random distribution of $M$ points $t_i \in (-T,T)$ on which one computes the function $F_N(t)$. The histogram of the $M$ values $F_N(t_i)$ gives rise to a curve very close to a normal distribution (see Figure \ref{normalllKac}), as confirmed by the various indicators of the fit. 

\vspace{3mm}
\noindent
{\bf What happens at a given $t$?}
As we have seen, the Kac's theorem concerns with the fraction of time (in the interval $(-T,T)$ and in the limit $T \rightarrow \infty$) that the function $F_N(t)$ spends in a given range of values $(a,b)$. Let's now pose a different question. Suppose that we fix a given instant of time $t$: can we say how the series $F_N(t)$ goes with $N$? It is worth stressing that the answer to this question does {\em not} follow from the Kac's theorem and it is not straightforward: it depends both on the nature of the sequence $\{\lambda_k\}$, in particular how these frequencies grow with the index $k$, and also the value of $t$ chosen. It is easy to identify some simple instances which show that this is indeed the situation. 

Let's consider an arbitrary sequence of admissible frequencies, i.e. linearly independent on the field of the rational numbers, and let's first study the behavior of the function $F_N(t)$ 
given in eq.\,(\ref{FNNNN}), as function of $N$ at some particular values of $t$. For instance, if $t=0$, the numerator of $F_N(0)$ goes as $N$ and therefore the series $F_N(0)$ grows as $\sqrt{N}$ when $N \rightarrow \infty$. On the other hand, for $t \rightarrow \infty$, the rapidly oscillating angles of the various cosines average to $0$ and therefore in this limit we are essentially in the condition of the Kac's theorem, so that we can conclude that, for $t\rightarrow \infty$, $F_N(t) \sim {\cal O}(1)$ when $N \rightarrow \infty$. 
However, beside these simple cases, it is quite difficult to draw some general conclusions on the behavior of $F_N(t)$ at a generic value of $t$ when $N \rightarrow \infty$ for a generic sequence of admissible frequencies.  

Deepening the analysis, to have $F_N(t) \sim {\cal O}(1)$ for $N \rightarrow \infty$ at a given $t$ it is of course necessary that the sequence $\{\cos\lambda_k t\}$ ($k=1,2,\ldots,N)$ has zero average: this condition is guaranteed if the sequence of the angles $\{\theta_k\}$ associated by the fractional part\footnote{It only matters the fractional part for the periodicity of the cosine.}  of $\frac{\lambda_k t}{2\pi}$, i.e. 
\beq
\theta_k \,=\, \frac{\lambda_k t}{2\pi} - \left[\frac{\lambda_k t}{2\pi}\right] \,\,\,, 
\label{fractionalparrt}
\eeq
is equidistributed on the interval $(0,1)$ or, at least, symmetric distributed under the transformation $\theta \rightarrow 1/2 - \theta$. This is the case, for instance, for a sequence of admissible frequencies $\lambda_k$ that grow linearly with the index $k$ as $\lambda_k \sim k$: indeed, for any finite value of $t > t_*$ (where $t_* \sim 1$), the series $F_N(t)$ 
associated to these frequencies is always bounded for $N \rightarrow \infty$. In other words, the sum of the cosines in the numerator of $F_N(t)$ behaves in this case as a random walk. But, taking instead another admissible sequence of frequencies, e.g. the sequence $\{\log k\}$ of the logarithm of the integers, the corresponding sequence (\ref{fractionalparrt}) is {\em not} uniformly distributed in the interval $(0,1)$ \cite{Steundingunif} and this implies that the behavior of the series $F_N(t)$ associated to this sequence as a function of $N$ may be different for $t = 2 \pi$ and for $t \rightarrow \infty$. To show that this sequence is not uniformly distributed, we can use the Weyl criterion that states that a sequence $\{a_n\}$ is equidistributed modulo $1$ if and only if for all non-zero integers $m$
\beq
\lim_{n\rightarrow \infty} \frac{1}{n} \sum_{k=1}^n e^{2 \pi i m a_k} \,=\, 0 \,\,\,.
\label{Weyl}
\eeq
With $a_k = \log k$ and $m=1$, we have 
\beq
 \sum_{k=1}^n e^{2 \pi i  \log k} \,=\,\sum_{k=1}^N k^{2 \pi i} \,=\,
 \sum_{k=1}^N \left(\frac{k}{N}\right)^{2 \pi i} \, N^{2 \pi i} \sim 
 N^{1+2 \pi i} \,\int_0^1 u^{2\pi i} du \,=\,
 \frac{N^{1+2\pi i}}{1+2\pi i} 
 \eeq
 and therefore the limit (\ref{Weyl}) does not go to zero since it continues to oscillate.  

Notice that if we use the approximation $p_n \sim n \log n$ for the primes, the Weyl criterion seems also to suggests that the sequence $\{\lambda_k = \log p_k\}$ associated to the primes is not uniformly distributed since 
\beq
\sum_{k=1}^n e^{2 \pi i  \log p_k} \,=\,\sum_{k=1}^N p_k^{2 \pi i} \,=\,
 \sum_{k=1}^N \left(\frac{p_k}{p_N}\right)^{2 \pi i} \, p_N^{2 \pi i} \sim 
 N^{1+2 \pi i} \,\int_0^1 (u (\log u -1))^{2\pi i} du \,=\,{\cal A} \,
 N^{1+2\pi i}
 \eeq
where ${\cal A}$ is the finite value of the integral and therefore the limit (\ref{Weyl}) does not go to zero.

\section{Growth of the series $B_N(t)$}\label{growthBN}
In this Appendix we present a simple argument  showing that the scaling law for large $N$ behavior of the series $B_N(t)$ 
for $L$-functions of non-principal cases is {\em independent} of  $t$ and it is completely fixed by the large $N$ behavior of the series at $t=0$. For what follows, it is important that $\theta_{p_n} \neq 0$ and this is the crucial condition that distinguishes this case from the one analysed at the end of the previous Appendix (see the paragraph {\em What happens at a given $t$?}). Let $\alpha$ be the exponent associated to the asymptotic behavior in $N$ of the series $B_N(t)$ evaluated at $t=0$  
\beq
C_N \equiv B_N(0) \,=\, \sum_{\substack{n=1 \\ p_n \nmid q}}^{N}    \cos \( \theta_{p_n} \)   \simeq N^\alpha 
\,\,\,\,\,\,\,\,
,
\,\,\,\,\,\,\,\,
N \rightarrow \infty
\,\,\,.
\label{GseriesAp}
\eeq
We want to show that this asymptotic behavior in $N$ of the series is the {\em same} for any finite value of $t$ and that $t$ can eventually only affect the prefactor ${\cal A}(t)$ 
\beq
B_N(t) \simeq  {\cal A}_\chi(t) \, N^{\alpha}
\,\,\,\,\,\,\,\,
,
\,\,\,\,\,\,\,\,
N \rightarrow \infty
\,\,\,
.
\eeq
We will use a {\em reductio ad absurdum} argument, initially assuming that there exists a finite value of $t$, say $t=t^*$, for which, {\em up to a given} $N$, the series $B_N(t)$ grows with a different exponent\footnote{The exponent $\beta$ cannot be less than $1/2$ because, as already commented in the text, it is already known that at least a certain number of zeros of the $L$-functions are on the critical line, see for instance \cite{Conrey,Hughes,Conrey2}: from the convergence properties of the integral (\ref{importantintegral}), this implies that the exponent $\beta$ must satisfy $1/2 \leq \beta \leq 1$.} $\beta$ as
\beq
B_N(t^*) \, =\, \sum_{\substack{n=1 \\ p_n \nmid q}}^{N}    \cos \(t^* \log p_n - \theta_{p_n} \) \simeq 
N^\beta 
\label{betaexponent}
\,\,\,.
\eeq
If $t^*$ is finite, it certainly exists an integer $M$ such that
\beq 
t^* \simeq M\,\,\,,
\label{posing}
\eeq
and we can have two cases: (a) $ M < N$ or (b) $M > N$. Since both cases finally lead to the same conclusions, we can assume that $M > N$. Once we have identified such an integer $M$, let's study how the series $B_N(t^*)$ behaves if we now change the upper extreme $N \rightarrow \tilde N$, with $\tilde N \gg N$. We can split the new sum $B_{\tilde N}(t^*)$ into three terms 
\beq 
B_{\tilde N}(t^*) \,=\, \sum_{n=1}^{N}[\cdots] + 
 \sum_{n=N+1}^{\tilde M}[\cdots] 
 + \sum_{n=\tilde M+1}^{\tilde N}[\cdots]  
\label{threeterms}
 \eeq
where $ N < M  < \tilde M$, with $ \tilde M \gg M$, and in the three cases, $[\cdots]$ stays for $\cos \(t^* \log p_n - \theta_{p_n} \)$. 

\vspace{1mm}
 
From eq.\,(\ref{betaexponent}), the first term goes as $N^\beta$ but since we are now interested in how the series goes with the {\em new} upper extremum $\tilde N$, this term is just a constant value for the new series $B_{\tilde N}(t)$ and therefore can be safely discarded. The same is also true for the second term in (\ref{threeterms}). The key term is then the last one, on which we now focus our attention. Let's now divide the interval $(\tilde M+1,\tilde N)$ into $k$ intervals, with $k$ sufficiently large: the length $d$ of these intervals is then   
\beq
d = \frac{(\tilde N - \tilde M+1)}{k}\simeq \frac{\tilde N}{k} + {\mathcal O}\left(\frac{1}{\tilde N}\right)
 \,\,\,, 
\eeq
since we can also take $\tilde N \gg \tilde M$. In this way, $\tilde N \simeq k d$. In the following we assume for simplicity that $d$ is an integer. We can now show that in each of these intervals the change of the phases $\psi_n(t^*)= t^* \log p_n -\theta_{p_n}$, varying $n$, depends essentially only on the angles $\theta_{p_n}$ and not on $t^*$:  indeed, going from the prime $p_q$ to $p_{q+m}$, where $q \gg M$ and $1 \leq m \leq d$, we have 
\beq
\psi_{p_{q+m}} - \psi_{p_q} \equiv \Delta \psi_{p_{q + m},p_q} \,=\, t^* \,\log\left(\frac{p_{q+m}}{p_q}\right) -  
\Delta\theta_{p_{q+m},p_q} 
\simeq \frac{m t^*}{q} - \Delta\theta_{p_{q+m},p_q}\,\,\,,
\eeq
where $\Delta\theta_{p_{q+m},p_q} = \theta_{p_{q+m}} - \theta_{p_q}$ and we have used $p_n \sim n \,\log n$ to evaluate the change of phase due to the first term. In view of eq.\,(\ref{posing}) and the condition $q \gg M$, one can see that the first term can be made  infinitesimal and completely negligible with respect to the second one, in particular it can be made more and more negligible simply increasing $\tilde N$ and $\tilde M$. Therefore, varying $m$, the change of phases is only due to the angles $\theta$'s and does not depend on $t^*$.  So
\beq
\cos\psi_{p_q+m} \,=\,\cos(\psi_{p_q} + \Delta\psi_{p_q+m,p_q})\simeq 
\cos(\psi_{p_q} - \Delta\theta_{p_{q+m},p_q})\,\,\,,
\eeq  
where $\psi_{p_q}$ is a fixed quantity. With 
\beq 
\cos(\psi_{p_q} - \Delta\theta_{p_{q+m},p_q}) \,=\,\cos\psi_{p_q} \, \cos\Delta\theta_{p_{q+m},p_q} 
+ \sin\psi_{p_q} \, \sin\Delta\theta_{p_{q+m},p_q}\,\,\,,
\eeq
summing on a sufficient large number $d$ of these terms and using the scaling law (\ref{Gseries}) (which also 
holds for the sum on the sin's), their sum in each interval goes as 
\beq
\sum_{n=q}^{q+d} [\cdots] \,\simeq \,A_q d^\alpha\,\,\,, 
\eeq
where $A_q$ is a constant. Since there are $k$ of these contributions, then the last term in (\ref{threeterms})  goes as 
\beq
 \sum_{n=\tilde M+1}^{\tilde N}[\cdots]  
\simeq (A_1+ \cdots A_k)  \, d^\alpha \,=\, \simeq {\cal A} \,\left(\frac{\tilde N}{k}\right)^\alpha \simeq \tilde N^\alpha
\,\,\,.
\eeq

\vspace{3mm}

In conclusion, even if one would assume the existence of finite value of $t^*$ for which, up to a certain $N$, the series $B_N(t^*)$ scales as $N^\beta$, with $\beta \neq \alpha$, going to larger values of $N$, $N \rightarrow \tilde N$, the series will be driven to the scaling behavior $\tilde N^\alpha$, where the exponent $\alpha$ is defined in eq.\,(\ref{Gseries}). Hence the asymptotic behavior of the series $B_N(t)$ for large values of $N$ is ruled by the asymptotic behavior of the series at $t=0$ and therefore the GRH relies only on the Theorem \ref{seriesinzero} stated in the text.

\end{document}